\PassOptionsToPackage{unicode}{hyperref}
\PassOptionsToPackage{naturalnames}{hyperref}
\documentclass[a4paper,11pt]{article}
\usepackage[affil-it]{authblk}
\usepackage{etex}
\usepackage{graphicx}

\usepackage{caption,subcaption} 

\usepackage{amsmath}
\usepackage{amssymb}
\usepackage{amsthm}
\usepackage{empheq}
\usepackage{multirow}
\usepackage{makecell}
\usepackage{pict2e}
\usepackage{booktabs}
\usepackage[figuresright]{rotating}
\usepackage{bbm}
\usepackage[all,light]{draftcopy}
\usepackage{tikz}

\usetikzlibrary{calc,trees,positioning,arrows,chains,shapes.geometric,%
    decorations.pathreplacing,decorations.pathmorphing,shapes,%
    matrix,shapes.symbols,external,fadings,patterns}

\usepackage{algorithm}      
\usepackage{pgfplots}
\usepackage{listings}		
\usepackage[noend]{algpseudocode}
\usepackage{changepage}     
\usepackage{color}              
\definecolor{myred}{rgb}{1,0,0}

\usepackage{epstopdf}
\usepackage{graphics} 
\usepackage{epsfig} 
\usepackage[ plainpages = false, pdfpagelabels,
	     pdfpagelayout = useoutlines,
	     bookmarks,
	     bookmarksopen = true,
	     bookmarksnumbered = true,
	     breaklinks = true,
	     linktocpage,
	     pagebackref,         
	     colorlinks = false,  
	     linkcolor = blue,
	     urlcolor  = blue,
	     citecolor = red,
	     anchorcolor = green,
	     hyperindex = true,
	     hyperfigures
	     ]{hyperref}
\usepackage{pdflscape}

\usepackage{todonotes} 

\usepackage{lineno}
\DeclareMathOperator\supp{supp}
\DeclareMathOperator\err{err}

\newcommand{\inftyint}{\int_{-\infty}^{\infty}}

\newcommand{\maboxinormal}{\frac{1}{\sqrt{2 \pi}}e^{-\xi^2 /2}}
\newcommand{\maboxinormalflat}{1/\sqrt{2 \pi}\;e^{-\xi^2 /2}}
\newcommand{\dx}{\partial_x}
\newcommand{\dt}{\partial_t}
\newcommand{\dxi}{\partial_{\xi}}
\newcommand{\kappas}{(\kappa_1, \ldots, \kappa_{n-3})}

\newcommand{\seta}{\sqrt{\theta}}
\newcommand{\asys}{A_{sys}}
\newcommand{\reals}{\mathbb{R}}
\newcommand{\naturals}{\mathbb{N}}

\newtheorem{theorem}{Theorem}[section]

\theoremstyle{remark}
\newtheorem{remark}[theorem]{Remark}

\theoremstyle{definition}
\newtheorem{definition}[theorem]{Definition}

\newcommand{\vect}[1]{\textrm{\boldmath${#1}$}} 

\newcommand\Eqn{Equation }

\newcommand\vel{u}          

\let\originalleft\left
\let\originalright\right
\renewcommand{\left}{\mathopen{}\mathclose\bgroup\originalleft}
\renewcommand{\right}{\aftergroup\egroup\originalright}

\binoppenalty=\maxdimen 
\relpenalty=\maxdimen   

\date{\today}             

\title{Spline Moment Models for the one-dimensional Boltzmann-Bhatnagar-Gross-Krook equation}

\author{Julian Koellermeier%
  \thanks{\texttt{koellermeier@zedat.fu-berlin.de}; Corresponding author}}
\affil{Institut f\"ur Mathematik, Freie Universit\"at Berlin,\\ 14195 Berlin, Germany, \\and\\ School of Mathematical Sciences, Peking University, \\ 100871 Beijing, China}
\author{Ullika Scholz%
  \thanks{\texttt{ullika.scholz@fu-berlin.de}}}
\affil{Institut f\"ur Mathematik, Freie Universit\"at Berlin,\\ 14195 Berlin, Germany}

\begin{document}

\maketitle


\begin{abstract}
We introduce Spline Moment Equations (SME) for kinetic equations using a new weighted spline ansatz of the distribution function and investigate the ansatz, the model, and its performance by simulating the one-dimensional Boltzmann-Bhatnagar-Gross-Krook equation.
The new basis is composed of weighted constrained splines for the approximation of distribution functions that preserves mass, momentum, and energy.
This basis is then used to derive moment equations using a Galerkin approach for a shifted and scaled Boltzmann-Bhatnagar-Gross-Krook equation, to allow for an accurate and efficient discretization in velocity space with an adaptive grid. The equations are given in compact analytical form and we show that the hyperbolicity properties are similar to the well-known Grad moment model.
The model is investigated numerically using the shock tube, the symmetric two-beam test and a stationary shock structure test case. All tests reveal the good approximation properties of the new SME model when the parameters of the spline basis functions are chosen properly. The new SME model outperforms existing moment models and results in a smaller error while using a small number of variables for efficient computations.

\vspace{4mm}
\end{abstract}

\section{Introduction}
\label{sec:1}

Kinetic theory is used in many applications, especially to derive simplified model equations for systems consisting of small, colliding particles that can no longer be described by standard continuum equations \cite{cercignani1994,Struchtrup2006}. For rarefied gases and micro or nano flows, the Boltzmann-Bhatnagar–Gross–Krook (BGK) equation is often used as a simplification of the full Boltzmann equation when developing mathematical models. However, it is still difficult to solve even the one-dimensional Boltzmann-BGK equation due to the additional dimension of the phase space, which consists of time, physical space and velocity space \cite{Cercignani1972}. On the other hand, simpler standard continuum dynamics models like the Euler or Navier-Stokes equations do not give accurate solutions and even miss so called rarefaction effects in the regime of large Knudsen number \cite{Grad1949} so that extended models are necessary, \cite{Bobylev2006,Torrilhon2016}.

One approach is the direct discretization of velocity space, called Discrete Velocity Method (DVM) \cite{Baranger2012,dechriste2014,Filbet2013}. However, a large number of discrete velocities is necessary to achieve a sufficient accuracy. These models thus often lack computational efficiency, even though they are easy to parallelize. Locally adaptive versions try to circumvent this by shifting the velocity grid in every step \cite{Brull2014}. A method that is also based on the use of discrete velocity points is the Unified Gas Kinetic Scheme (UGKS) \cite{Xu2010}.

A straightforward extension of the Navier-Stokes equations leads to the so-called Burnett equations, which have been proven to be unstable \cite{Struchtrup2006}. Another method to derive models for simulations of rarefied gases is via the Nonlinear Coupled Constitutive Relations (NCCR) \cite{Zhongzheng2019}.

A different approach is to employ moment models that use a special ansatz for the distribution function, including only a small number of variables that capture the distribution function well for small to medium deviations from equilibrium \cite{Struchtrup2006,Torrilhon2016}. This is based on the seminal work by Grad \cite{Grad1949}, who used a weighted Hermite polynomial ansatz for the distribution function. The moment models were extended with great success in recent years by \cite{Cai2010,Torrilhon2015}. The resulting system of equations is called moment equations, for example, the Hyperbolic Moment Equations (HME) or Quadrature-Based Moment Equations (QBME) method \cite{Cai2013,Fan2016,Koellermeier2014}, and is based on a polynomial expansion or the maximum entropy method \cite{Levermore1996,McDonald} using an exponential ansatz. For recent progress and efficient implementations of the maximum entropy method, we refer to \cite{boehmer2020,Sadr2020,schaerer2015}. All existing models typically balance between computational efficiency, accuracy, and analytical properties like hyperbolicity, conservation properties or existence of an entropy. We refer to the review article \cite{Torrilhon2016} for further details.

The explicit terms of the moment equations are largely decided by the choice of the ansatz for the distribution function. It can be a simple, piecewise constant distribution function as in the DVM case, a sum of globally supported basis functions as in the case of HME, or an even more non-linear exponential ansatz like the maximum entropy method, where no closed form of the equations exists in the general case.

In this paper we present a new model called Spline Moment Equations (SME), which is based on the choice of weighted spline functions as basis for the reconstruction of the distribution function. Even though this is to some extend a consistent extension of the DVM method, up to our knowledge, there is no work in the literature about moment models based on piecewise, higher degree polynomials like splines.

Several works using piecewise polynomials for the Boltzmann equation are limited to the space-time domain, see e.g. \cite{Kitzler2017,KITZLER2015}. There exists some work in the context of a direct discontinuous Galerkin discretization of the velocity space for computation of the collision operator \cite{Zhang2018,Alekseenko2012,Alekseenko2014}, with a focus on the collision term of the full Boltzmann equation. However, the simultaneous discretization in the full phase space does not lead to a closed moment model in the sense of the above methods. Furthermore, the discretization using a fixed grid in velocity space does not allow for an efficient, adaptive approximation in velocity space. We will cover both points in this paper.

Standard splines, i.e. piecewise polynomials with bounded support and certain continuity properties, are not a-priori suitable for the reconstruction of typical distribution functions in kinetic theory. This is because conservation of mass, momentum, and energy, which are integral values of the distribution function, need to be conserved. We thus use a specific linear combination of four weighted B-splines that ensures compatibility with the respective conservation laws. We call those \emph{weighted fundamental constrained splines}. We study the approximation properties of these basis functions in comparison to an unweighted or unconstrained spline basis to identify these weighted constrained splines as a suitable basis.

We combine the spline approach with a shifting and scaling of the distribution function to allow for an efficient and accurate approximation of the Boltzmann equation, first used in \cite{Kauf2011} and applied to globally defined Hermite basis functions in \cite{Koellermeier2014b}. This leads to a transformed adaptive grid in the velocity space and greatly decreases the number of variables needed. The SME are obtained by insertion of the spline basis ansatz into the Boltzmann equation and by Galerkin projecting onto respective spline test functions.

A slightly different approach was described in \cite{Kauf2011} and briefly tested for some spline functions, without an in depth study of the approximation properties or numerics used. In \cite{Koellermeier2014,Koellermeier2014b} the approach by Kauf was used for weighted Hermite basis functions to derive hyperbolic moment equations. Without any modification, the Hermite-based moment equations will only be hyperbolic in a bounded domain of the variable space. This might lead to numerical instabilities and a breakdown of the simulation away from equilibrium. Many models have been developed to solve this problem and to arrive at globally hyperbolic models \cite{Fan2016,Koellermeier2017b}, see also the new diagram based approach to study hyperbolicity \cite{Koellermeier2020b}. We therefore study hyperbolicity and show that the SME has a similar behavior as Grad's equations. In the appendix we introduce a simple, linearized version of the equations that retains hyperbolicity at the expense of some accuracy.

The simulation results of the one-dimensional Boltzmann-BGK equation for the shock tube test case yield the first systematic investigation of a moment model based on piecewise polynomial spline functions. Further tests for the symmetric two-beam problem and the stationary shock structure test case show that the new SME model outperforms several existing hyperbolic moment models, despite using only a small number of variables and a straightforward basis definition.

The rest of this paper is organized as follows: In Section \ref{sec:2} the one-dimensional Boltzmann-BGK equation is briefly described. Section \ref{sec:3} outlines the construction of approximation spaces via spline basis functions and examples. Different spline basis functions are investigated in detail with respect to their ability to approximate standard distribution functions in Section \ref{sec:4}. The derivation of moment equations from the spline basis functions and its stability analysis is described in Section \ref{sec:5}, before three different numerical tests highlight the good approximation quality of the new SME model in Section \ref{sec:6}. The paper ends with a short conclusion.

\section{Boltzmann transport equation}
\label{sec:2}
According to \cite{bookRGD} the Boltzmann transport equation (BTE) is used to model fluid motions where the molecules' mean free path $\lambda$ is large in comparison to a reference length scale $L$ (for example, the diameter of a tube or the curvature radius of a space shuttle). The mean free path is the average distance that a molecule travels before it collides with another molecule. The dimensionless flow parameter, the Knudsen number, is denoted as $\textrm{Kn} = \frac{\lambda}{L}$, for details see e.g. \cite{Koellermeier2017b,Struchtrup2006}.
Typical applications for a large Knudsen number are rarefied gases and vacuum technology, as stated in \cite{bookRGD}. For small Knudsen numbers the BTE blends into the standard macroscopic conservation laws, i.e. the Navier-Stokes or Euler equations \cite{bookRGD}.

For conciseness, we only use the one-dimensional spatial case here and leave the extension to the multi-dimensional case for further work.

The BTE describes the dynamics of the distribution function $f$ which is a probability density defined at every position in space $x \in \reals$, time $t \in \reals_+$ and velocity space $c \in \reals$ \cite{bookRGD,Koellermeier2017b}.
\begin{equation}
 \dt f(x,t,c) + c \dx f(x,t,c) = S(f)
\end{equation}
The left-hand side describes the transport of particles, whereas the right-hand side contains the collision operator. Throughout the paper, we will use the BGK collision term $S(f)$ \cite{paperPGK}. In this model the distribution function $f$ relaxes towards the equilibrium function $f_{M}$ as a consequence of the collisions following
\begin{equation}
 S(f)=-\frac{1}{\tau}(f-f_{M}).
\end{equation}
Here $\tau$ is the relaxation time. Depending on the gas properties, it may be large, i.e. few collisions and slow relaxation towards $f_M$, or small, i.e. many collisions and fast relaxation of particle velocities towards equilibrium. The equilibrium limit $f_M$ is the Maxwellian distribution function, where the particle velocities are distributed according to a bell curve with its position and scaling determined by the macroscopic quantities density $\rho \in \reals$, temperature $\theta \in \reals$ and velocity $v \in \reals$.
\begin{equation}
\label{e2:Maxwellian}
 f_{M}=\frac{\rho(x,t)}{\sqrt{2\pi \theta(x,t)}} \exp{(-\xi^2 /2)}, \quad \xi = \frac{c-v(t,x)}{\sqrt{\theta(t,x)}}
\end{equation}
Density $\rho$, velocity $v$ and temperature $\theta$ are linked to the distribution function $f$ via integration in velocity space.
\begin{gather}
\label{makro}
 \begin{split}
 \rho(t,x) &= \inftyint f(t,x,c)dc \\
 \rho(t,x)v(t,x) &= \inftyint c f(t,x,c)dc \\
 \rho(t,x)\theta(t,x) &= \inftyint \vert c-v\vert^2f(t,x,c)dc.
\end{split}
\end{gather}
The Boltzmann-BGK equation thus forms a non-linear integro-differential equation. It is hard to solve the equation numerically because discretization is necessary in the microscopic velocity $c \in \reals$, the time $t \in \reals_+$ and the position in space $x \in \reals$.

\section{Spline basis functions}
\label{sec:3}
Many approximations exist to derive accurate discretizations of the BTE in velocity space. However, up to our knowledge there is no work investigating the natural use of spline ansatz functions, although splines are used in many applications in engineering with great success. In this paper, we will therefore use the flexible spline ansatz for the distribution function. We parameterize the distribution function by means of ansatz functions multiplied by respective coefficients. In this section we will introduce the ansatz functions that compose the basis space. The two classes of possible ansatz functions presented here are B-splines and fundamental constrained splines (FCS).

\subsection{B-splines}
\label{bsplineexamples}
B-splines were first presented as \emph{fundamental spline functions} to serve as a basis for the known spline functions \cite{paperSC}. The class of spline functions is very popular and for example used in polynomial interpolation. Its members are piecewise functions $S: \reals \mapsto \reals$ consisting of polynomials of equal degree.
The pieces are defined by a grid $G_{\xi}=[\xi_0, \hdots, \xi_{n-1}]$ and the polynomial degree is called the spline's \emph{order} $k \in \naturals$. A spline of order $k$ is at least $k-1$ times continuously differentiable by construction.

Because of its simplicity we will here use the recursive definition for B-splines \cite{bookPGS}. The initial definition is not recursive and can be found in \cite{paperSC}.

\begin{definition}[B-spline]
    Let $G_{\xi}: \hdots < \xi_{-2} < \xi_{-1} < \xi_{0} < \xi_{1} <  \xi_{2} < \hdots  ,\: \: \xi_i \in \reals$ be a grid.

    The $j$th B-spline of order $k=1$ is defined by
    \begin{equation}
    \label{BSplineOder1Def}
    B_{j1}(\xi) :=
    \begin{cases}
     \frac{\xi-\xi_j}{\xi_{j+1}-\xi_{j}} & \xi \in [\xi_j,\xi_{j+1}] \\
     \frac{-\xi+\xi_{j+2}}{\xi_{j+2}-\xi{j+1}} & \xi \in [\xi_{j+1},\xi_{j+2}] \\
     0 & \textrm{else}.
    \end{cases}.
    \end{equation}
    B-splines of higher order are obtained by the recursion formula
    \begin{equation}
    \label{SplineDef}
     B_{jk}(\xi)=\omega_{jk}^{(1)}B_{j,k-1}(\xi)+\omega_{jk}^{(2)} B_{j+1,k-1}(\xi),
    \end{equation}
    with
    \begin{equation}
     \omega_{jk}^{(1)}(\xi):=\frac{\xi-\xi_j}{\xi_{j+k}-\xi_j} \quad \textrm{and} \quad \omega_{jk}^{(2)}(\xi):=\frac{-\xi+\xi_{j+k+1}}{\xi_{j+k+1}-\xi_{j+1}}.
    \end{equation}
\end{definition}

\begin{remark}
    The recursion formula implies that a B-spline of order $k$ is a composition of two neighboring B-splines of order $k-1$. The B-splines' support therefore increases by one grid cell per order, that is
    $\supp B_{jk}(\xi)=[\xi_j,\xi_{j+k+1}]$.
    A proof can be found in \cite{bookPGS} p. 91. On an equidistant grid of width $\Delta \xi$ this simplifies to $\supp B_{jk}(\xi)=[\xi_j,\xi_{j}+(k+1)\Delta \xi]$.
\end{remark}

\begin{figure}[H]
 \centering
 \includegraphics[scale=0.3]{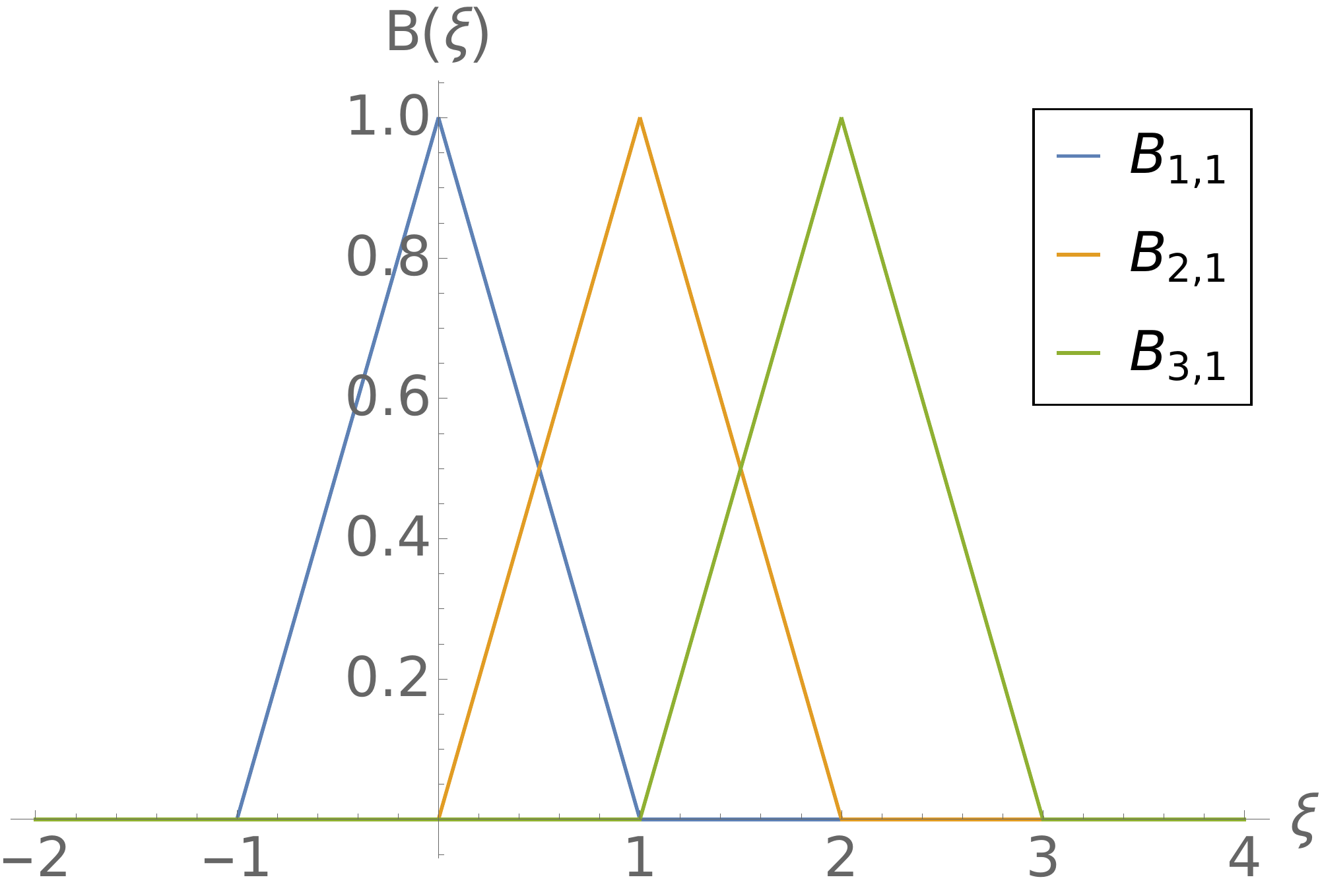}
 \caption{Linear splines on $G^{(1)}_{\xi}$, $\Delta \xi=1$.}
 \label{border1}
\end{figure}

We exemplify B-splines of first and second order. For simplicity we choose the grid $G_{\xi}$ in such a way that the first B-spline $B_{1k}$ is centered around the origin. Order $k=1$ then results in the grid
\begin{equation}
 G^{(1)}_{\xi} : \hdots < \xi_0 = -2\;\Delta \xi < \xi_1 = -\Delta \xi < \xi_2 = 0 < \xi_3 = \Delta \xi <  \hdots \: .
\end{equation}
The first B-spline of first order on $G^{(1)}_{\xi}$ is
\begin{equation}
B_{1,1}(\xi)=
\begin{cases}
 1-\frac{\xi }{\Delta \xi } & 0\leq \xi \leq \Delta \xi  \\
 \frac{\xi }{\Delta \xi }+1 & -\Delta \xi \leq \xi \leq 0 \\
 0 & |\xi|>\Delta \xi.
\end{cases}
\end{equation}
Figure \ref{border1} contains this B-spline $B_{1,1}$ together with its neighbors $B_{2,1}$ and $B_{3,1}$ on the grid $G^{(1)}_{\xi}$ for $\Delta \xi=1$. These first order B-splines are often denoted as \emph{hat functions}.

Splines of second order $k=2$ are piecewise quadratic polynomials. In comparison with \ref{border1} a different grid $G^{(2)}_{\xi}$ is used

\begin{equation}
G^{(2)}_{\xi} : \hdots < \xi_0 = -\frac{5}{2}\Delta \xi\; < \xi_1 = -\frac{3}{2}\Delta \xi < \xi_2 = -\frac{1}{2}\Delta \xi < \xi_2 = \frac{1}{2} \Delta \xi < \hdots \: .
\end{equation}

On this grid the first B-spline of second order is

\begin{equation}
B_{2,1}(\xi)=
\begin{cases}
 \frac{1}{8} \left(\frac{2 \xi }{\Delta \xi }-3\right)^2 & \frac{\Delta \xi }{2}\leq \xi \leq \frac{3 \Delta \xi }{2} \\
 \frac{3}{4}-\frac{\xi ^2}{\Delta \xi ^2} & -\frac{\Delta \xi }{2}\leq \xi \leq \frac{\Delta \xi }{2} \\
 \frac{1}{8} \left(\frac{2 \xi }{\Delta \xi }+3\right)^2 & -\frac{3 \Delta \xi }{2}\leq \xi \leq -\frac{\Delta \xi }{2}\\
 0 & |\xi|>\frac{3}{2} \Delta \xi.
\end{cases}
\end{equation}

\begin{figure}[H]
 \centering
 \includegraphics[scale=0.3]{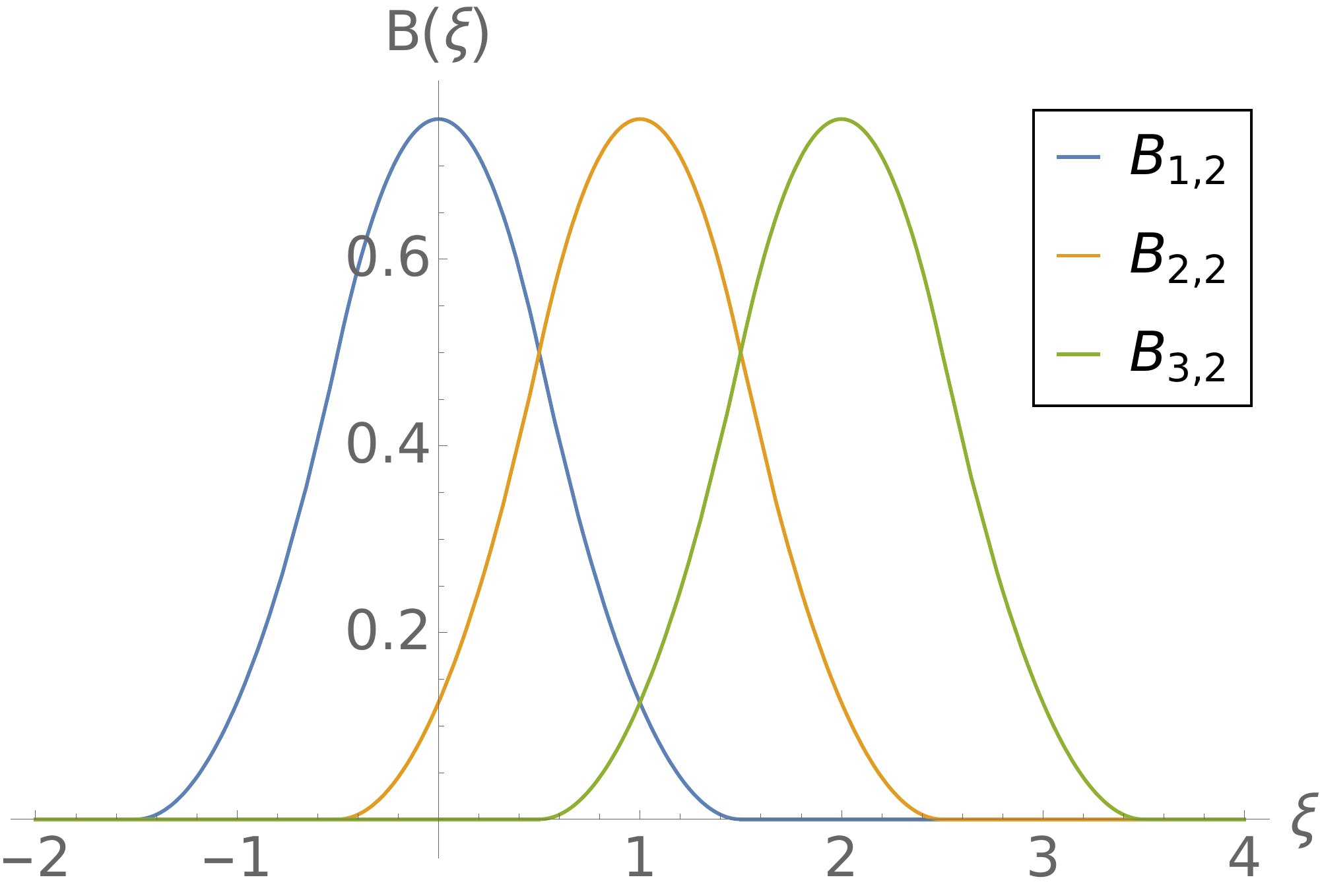}
 \caption{Quadratic splines on $G^{(2)}_{\xi}$, $\Delta \xi=1$.}
 \label{border2}
\end{figure}

This B-spline is depicted in Figure \ref{border2} along with its neighbors.

While B-splines are legitimate basis for piecewise polynomials, they do not conserve mass, momentum, and energy of the distribution function when used in an expansion in velocity space.

\subsection{Fundamental constrained splines}
Constrained splines are a subset of the spline functions. They will be used to discretize the velocity space for the Boltzmann equation because they preserve the distribution function's mass, momentum, and energy, which are related to the integral values $\rho$, $v$ and $\theta$, see Equation \eqref{makro}. In the derivations in Section \ref{sec:5} the distribution function $f(t,x,c)$ is scaled and shifted in the velocity variable to a new distribution function $f(t,x,\xi)$ so that the scaled distribution function has $\rho=1$, $v=0$ and $\theta=1$. Any constrained spline approximation function $S^{constr}$ then needs to fulfill the constraints
\begin{equation}
    \int_{-\infty}^{\infty} \maboxinormal S^{constr}(\xi) \begin{pmatrix}1\\\xi\\\xi^2\end{pmatrix} d\xi=\begin{pmatrix}0\\0\\0\end{pmatrix}.
    \label{compatibility}
\end{equation}
Any basis function for these constrained splines is called \emph{fundamental constrained spline} (FCS) $F^{constr}$. Their construction from an existing B-spline basis is described in \cite{Kauf2011} and will be presented below.

\begin{definition}[Fundamental constrained splines]
    Given a grid with four consecutive B-splines $B_{jk}, B_{j+1,k}, B_{j+2,k}, B_{j+3,k}$, a \emph{fundamental constrained spline} (FCS) $F_i^{constr}$ is assembled by linear combination

    \begin{equation}
     F_j^{constr}=a_0 B_{jk}+a_1 B_{j+1,k} + a_2 B_{j+2,k} + a_3 B_{j+3,k} \textrm{{\;,}}
    \end{equation}

    where the coefficients $\alpha_0, \ldots, \alpha_3 \in \reals$ are determined by

    \begin{gather}
    \begin{split}
    a_0=\frac{\tilde{a_0}}{\sqrt{\tilde{a_0}^2+\tilde{a_1}^2+\tilde{a_2}^2+1}}\;,\quad
    a_1=\frac{\tilde{a_1}}{\sqrt{\tilde{a_0}^2+\tilde{a_1}^2+\tilde{a_2}^2+1}}\;,\\
    a_2=\frac{\tilde{a_2}}{\sqrt{\tilde{a_0}^2+\tilde{a_1}^2+\tilde{a_2}^2+1}}\;,\quad
    a_3=\frac{\tilde{a_3}}{\sqrt{\tilde{a_0}^2+\tilde{a_1}^2+\tilde{a_2}^2+1}}\;.
    \end{split}
    \end{gather}

    Setting $\tilde{a_3}=1$, the other variables $\tilde{a_0}$, $\tilde{a_1}$ and $\tilde{a_2}$ are found by solving the linear system

    \begin{equation}
    \label{fcslgs}
    \sum_{j=0}^2 \tilde{a_j} \int_{-\infty}^{\infty} \begin{pmatrix}1\\\xi\\\xi^2 \end{pmatrix} e^{(-\xi^2/2)} B_{\alpha+j,k}(\xi)d\xi=-\int_{-\infty}^{\infty} \begin{pmatrix}1\\\xi\\\xi^2 \end{pmatrix} e^{(-\xi^2/2)} B_{\alpha+3,k}(\xi)d\xi.
    \end{equation}
\end{definition}

\begin{remark}
    A consequence of the necessity of four B-splines for the construction of one single fundamental constrained spline is that the dimension of the function space decreases by three when changing from a B-spline basis to a basis of fundamental constrained splines.
\end{remark}

\begin{remark}
    In the 1D setting discussed here, the three constraints \eqref{compatibility} lead to the reduction of the dimension of the function space by three. However, in a multi-dimensional setting, e.g. with a multi-dimensional spline basis obtained from the tensor product of one-dimensional splines, there will be more constraints and thus more splines that build a fundamental constrained spline. In a general, $d-$dimensional setting, there are $d+2$ constraints.
\end{remark}

\begin{figure}[h]
 \centering
 \begin{subfigure}[t]{0.47\textwidth}
 \includegraphics[width=\textwidth]{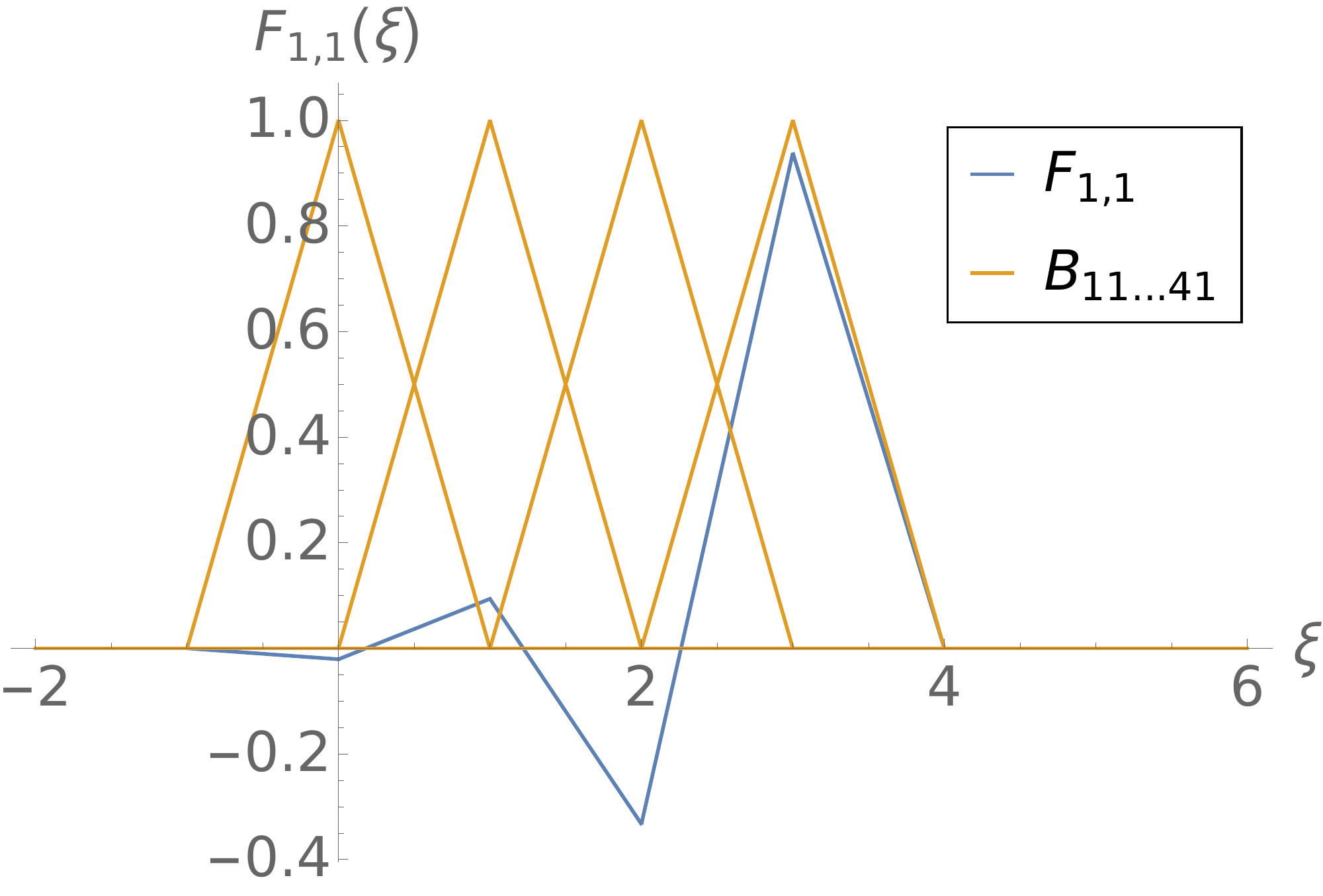}
 \caption{Linear element $F_{1,1}$.}
 \label{forder1}
 \end{subfigure}
 \hfill
 \begin{subfigure}[t]{0.47\textwidth}
 \centering
 \includegraphics[width=\textwidth]{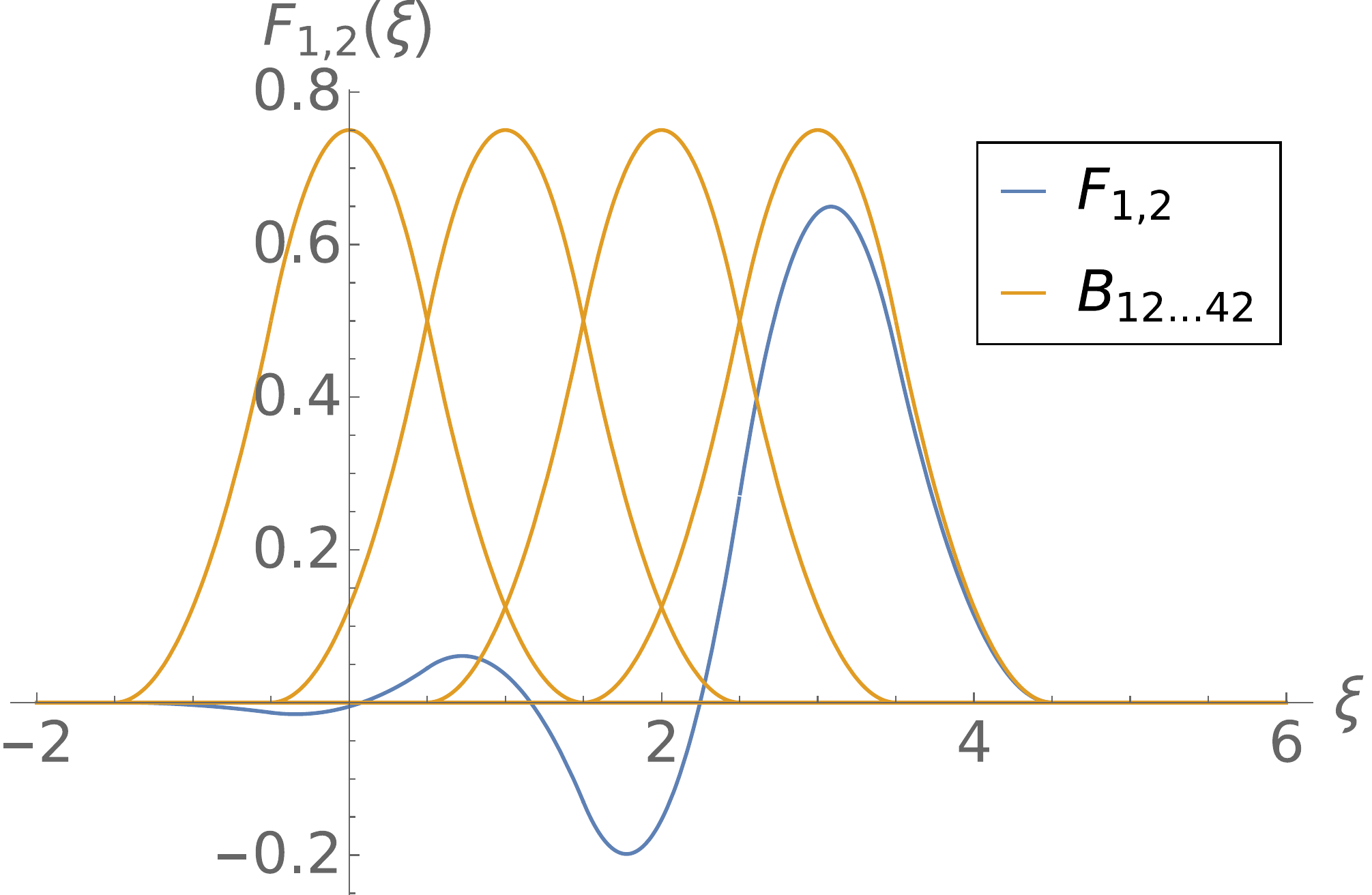}
 \caption{Quadratic element $F_{1,2}$.}
 \label{forder2}
 \end{subfigure}
\hfill
 \begin{subfigure}[t]{0.47\textwidth}
 \centering
 \includegraphics[width=\textwidth]{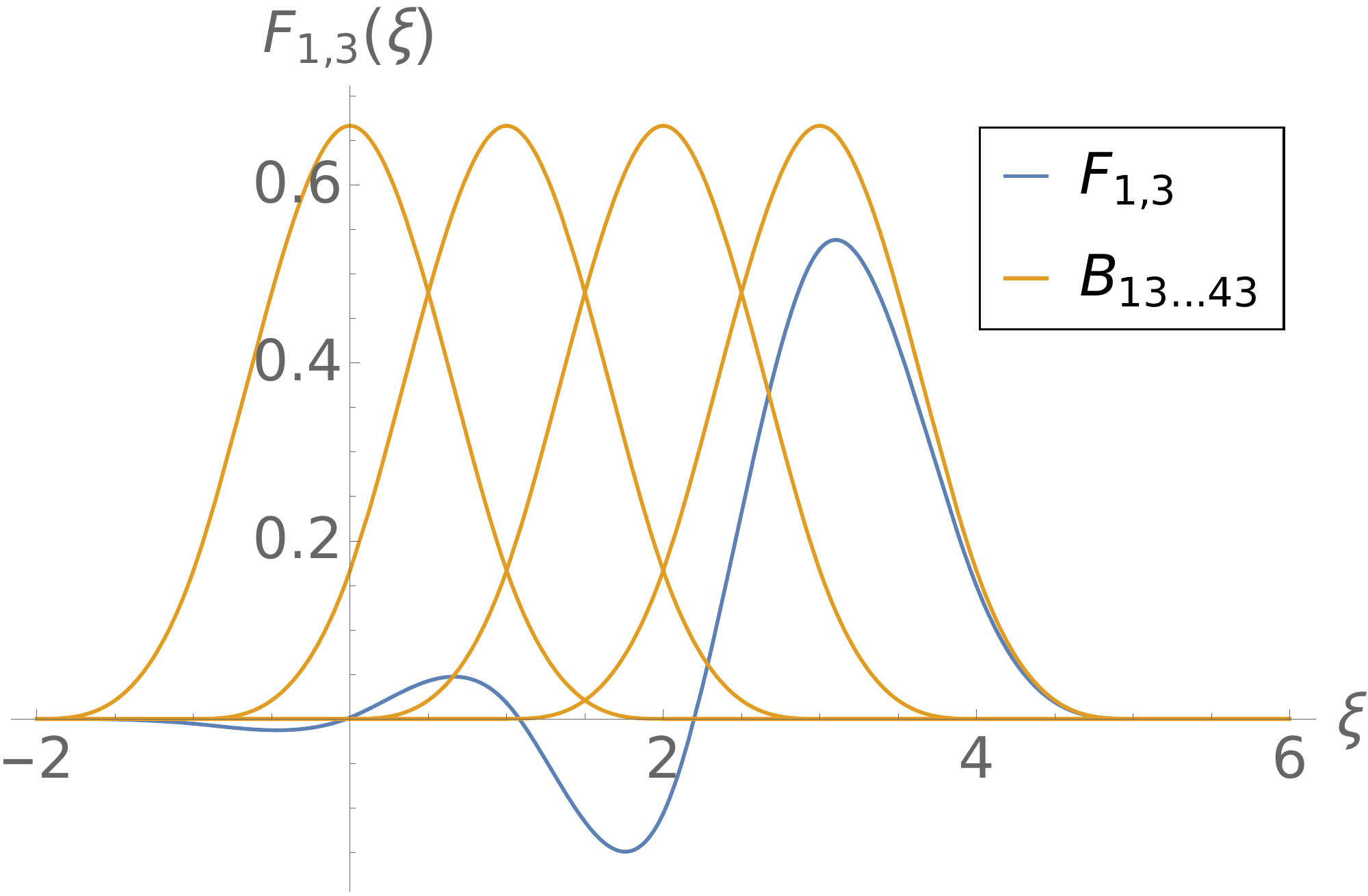}
 \caption{Cubic element $F_{1,3}$.}
 \label{forder3}
 \end{subfigure}
 \hfill
  \begin{subfigure}[t]{0.47\textwidth}
 \centering
 \includegraphics[width=\textwidth]{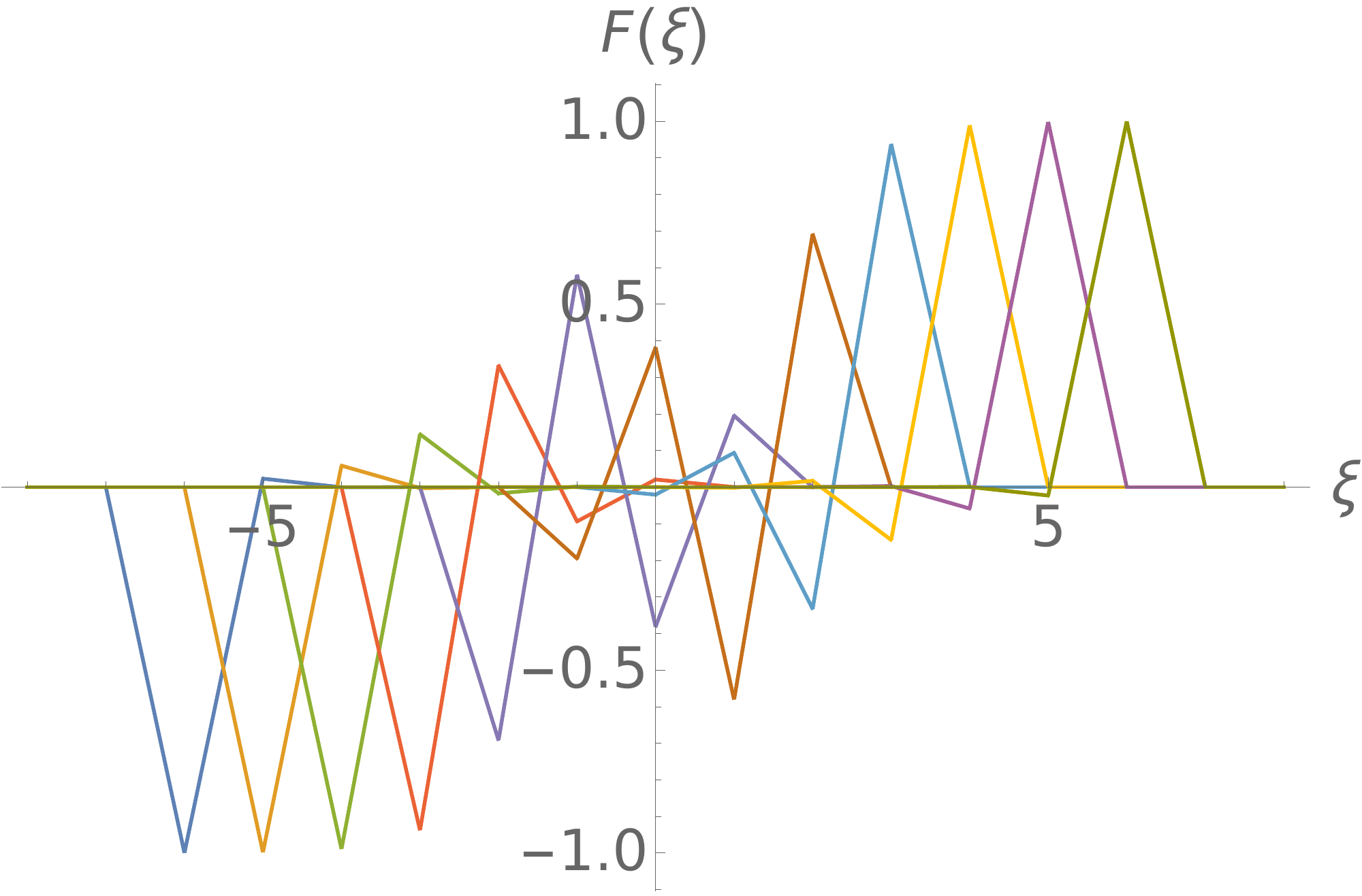}
 \caption{Complete linear basis, $13$ elements.}
 \label{forder4}
 \end{subfigure}
 \caption{Fundamental contrained spline basis for different order $k$, $\Delta \xi=1$. (a) $k=1$, (b) $k=2$, (c) $k=3$, (d) basis for $k=1$.}
  \label{forders}
\end{figure}

Figure \ref{forder1} shows a fundamental constrained spline of first order $F_{1,1}$. The underlying grid is $G^{(1)}_{\xi}$ from \Eqn \ref{bsplineexamples} with $\Delta \xi=1$. It can be seen that $F_{1,1}$ is a linear combination of the B-splines $B_{1,1}$, $B_{2,1}$, $B_{3,1}$ and $B_{4,1}$. The equivalent holds true for spline orders $k=2$ and $k=3$ as depicted in Figures \ref{forder2} and \ref{forder3}.

The full basis of fundamental constrained splines for this example with $13$ basis elements is pictured in Figure \ref{forder4}.

\begin{remark}
  In contrast to the B-splines in Figure \ref{border1}, the FCS have a larger support and change their shape depending on their position on the $\xi$-axis. The FCS will later be multiplied with a local equilibrium weight function, which is the scaled Maxwellian \ref{e2:Maxwellian}, i.e. $\maboxinormalflat$. Values in the region of large $\xi$ will thus be scaled down accordingly.
\end{remark}
\section{Spline approximations of distribution functions}
\label{sec:4}

Using B-splines and FCS, we investigate three possible ways for an expansion of the distribution function via splines called unweighted splines, weighted splines and weighted FCS. As test cases we will use selected bimodal functions, which often occur as distribution functions in kinetic equations.

Firstly, the distribution function will later be scaled and shifted in the velocity variable to a new velocity space $\xi$. The scaled distribution function then has $\rho=1$, $v=0$ and $\theta=1$, see Section \ref{sec:5}. We thus only consider examples of that form here.

Secondly, the ansatz should not be a mere linear combination of splines, due to the fact that the equilibrium function is a Maxwellian (or a Gaussian in the transformed space). Instead, we assume the distribution function to be a sum of a Gaussian distribution $\maboxinormalflat$ plus a linear combination of basis functions. The linear combination will thus ultimately describe the distribution function's \emph{deviation} from the Gaussian distribution, while a Maxwellian can be represented exactly by setting all basis coefficients to zero. This reduces the approximation error to zero in equilibrium.

Thirdly, to describe deviations from Maxwellians easily, also the basis functions can be weighted with that Maxwellian. To that extend, we include a weighting of the ansatz functions with the Gaussian distribution $\maboxinormalflat$ in the transformed space. This leads to suppressed oscillations towards the boundary of the discretization grid $\xi \in [\xi_{min},\xi_{max}]$. Another advantage is to have more balanced basis coefficients in terms of their absolute values because deviations from the equilibrium distribution are larger close the origin where $\xi$ is small.

We test weighted and unweighted B-splines as well as the weighted FCS as ansatz function (unweighted FCS do not make sense as the compatibility constraints already make use of the weight function).
Finally, the coefficients are found by means of a Galerkin method. The expansion is multiplied with test functions followed by integrating over velocity space and then solving the resulting system of equations.

We investigate the approximation properties using two bimodal functions
\begin{equation}
    f_{b1}(\xi)=0.527399 e^{-5.55556 (0.566569 \xi -0.3)^2}+0.169521 e^{-3.125 (0.566569 \xi +0.7)^2},
\end{equation}
\begin{equation}
    f_{b2}(\xi)=0.274485 e^{-3.125 (1.10085 \xi -0.525)^2}+0.274485 e^{-0.347222 (1.10085 \xi +0.175)^2}
\end{equation}
and the bimodal function used by Kauf in \cite{Kauf2011}
\begin{equation}
    f_{PK}(\xi)=0.16241 e^{-0.45116 (\xi -0.8)^2}+0.812051 e^{-6.34444 (\xi +0.6)^2}.
\end{equation}
The example functions were chosen such that they represent different shapes of bimodal functions frequently occurring in numerical simulations, for example, in the shock structure test case in Section \ref{sec:Shock_structure}.

\begin{figure}[H]
 \centering
 \includegraphics[scale=0.4]{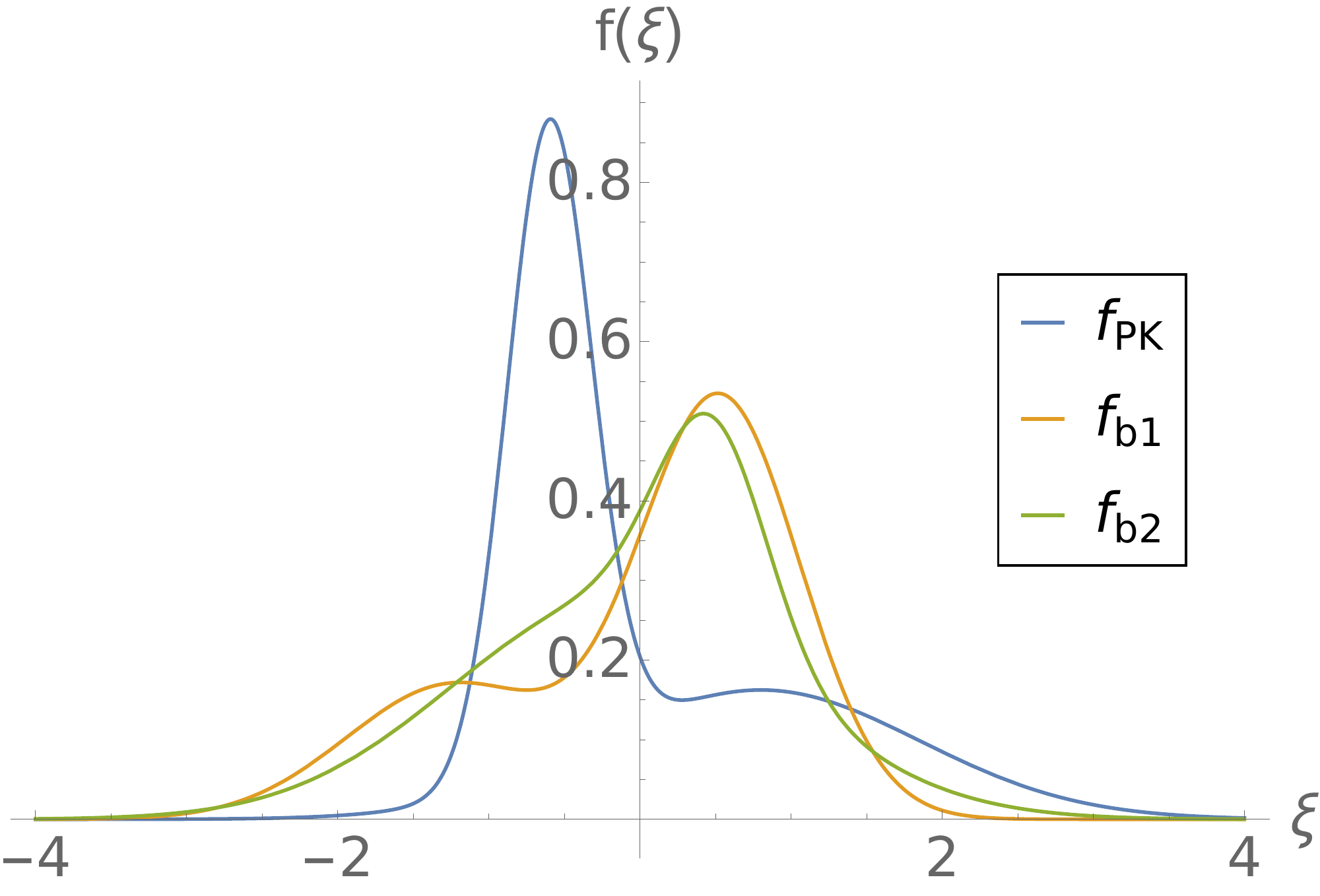}
 \caption{{Bimodal test functions} $f_{PK}$, $f_{b1}$, $f_{b2}$.}
 \label{bimods}
\end{figure}

Note that all distribution functions are chosen to have $\rho=1$, $v=0$ and $\theta=1$, to comply with the constraints in Equation \eqref{compatibility} in the transformed space, see Equation \eqref{trans} and Section \ref{sec:5}. All bimodal functions are plotted in Figure \ref{bimods}.

\subsection{Unweighted B-splines}
\label{vorgehen}
Unweighted B-splines expand the deviation of the distribution function from the equilibrium Gaussian in a standard spline series using the following definition.
\begin{definition}[unweighted B-spline expansion]
    Let $n \in \naturals$ and $k \in \naturals$ be fixed numbers. Let the $\xi$-range $[\xi_{min},\xi_{max}]$ also be fixed. Then we define the equidistant grid $G_{\xi}$ as
    \begin{equation}
        G_{\xi}: \xi_1 = \xi_{min}-\frac{\Delta \xi(k+1)}{2} < \xi_2 = \xi_1+\Delta \xi < \hdots < \xi_{n+k+1}=\xi_{max}+\frac{\Delta \xi(k+1)}{2}
    \end{equation}
    for the grid spacing $\Delta \xi=\frac{\xi_{max}-\xi_{min}}{n-1}$.
    The \emph{unweighted B-spline} expansion of a distribution function $f: \reals \mapsto \reals$ with a number of $n$ spline functions is then given by
    \begin{equation}
        \label{ansatz}
        \hat{f}(\xi)=\maboxinormal+\sum_{i=1}^{n} \alpha_i \phi_i(\xi).
    \end{equation}
    for B-splines $\phi_i = B_{ki}$ of order $k$ on $G_{\xi}$.
\end{definition}

\begin{figure}[H]
 \centering
 \includegraphics[scale=0.55]{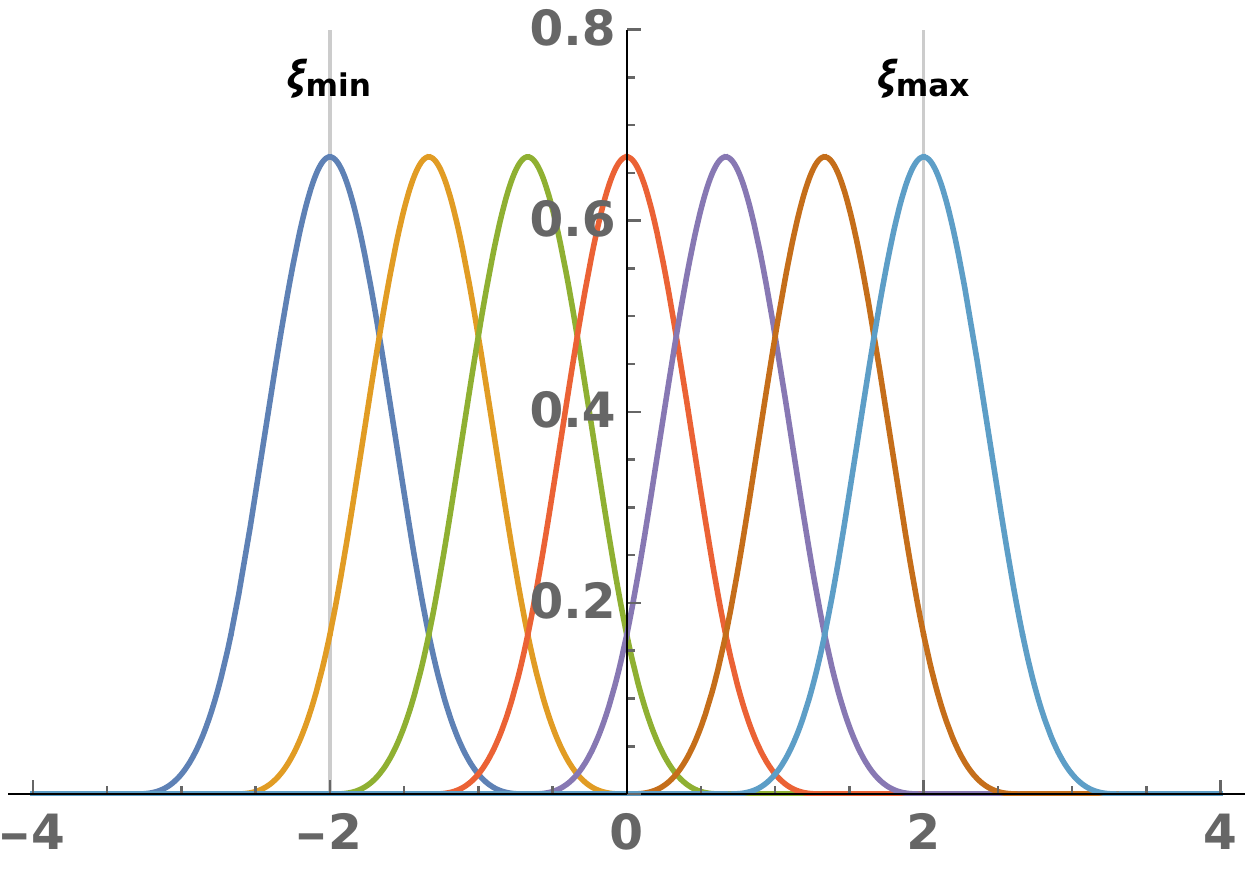}
 \caption{Example velocity grid with $7$ third-order splines. $[\xi_{min},\xi_{max}]=[-2,2]$.}
 \label{bereich}
\end{figure}

Constructing the B-spline according to Section \ref{sec:3}, we obtain exactly $n$ B-splines for an unweighted B-spline basis, where the outer B-splines' midpoints coincide with $\xi_{min}$ and $\xi_{max}$ (see Figure \ref{bereich}).

The basis coefficients are calculated by Galerkin projection of \ref{ansatz} onto test functions $\psi_j(x)=\phi_j$, $j \in 1,..,n$. This leads to the linear system of equations
\begin{equation}
    \label{lgs}
    A \alpha = u,
\end{equation}
with solution $\alpha \in \mathbb{R}^{n}$, while the entries $u \in \mathbb{R}^{n}$ and $A \in \mathbb{R}^{n \times n}$ are given by
\begin{equation}
    u_j =  \inftyint f(x) \psi_j dx, \quad A_{ij}=\inftyint \phi_i(x) \psi_j(x) dx. \label{ansatzend}
\end{equation}

\begin{figure}[H]
\centering
\begin{subfigure}[b]{0.47\textwidth}
         \centering
         \includegraphics[width=\textwidth]{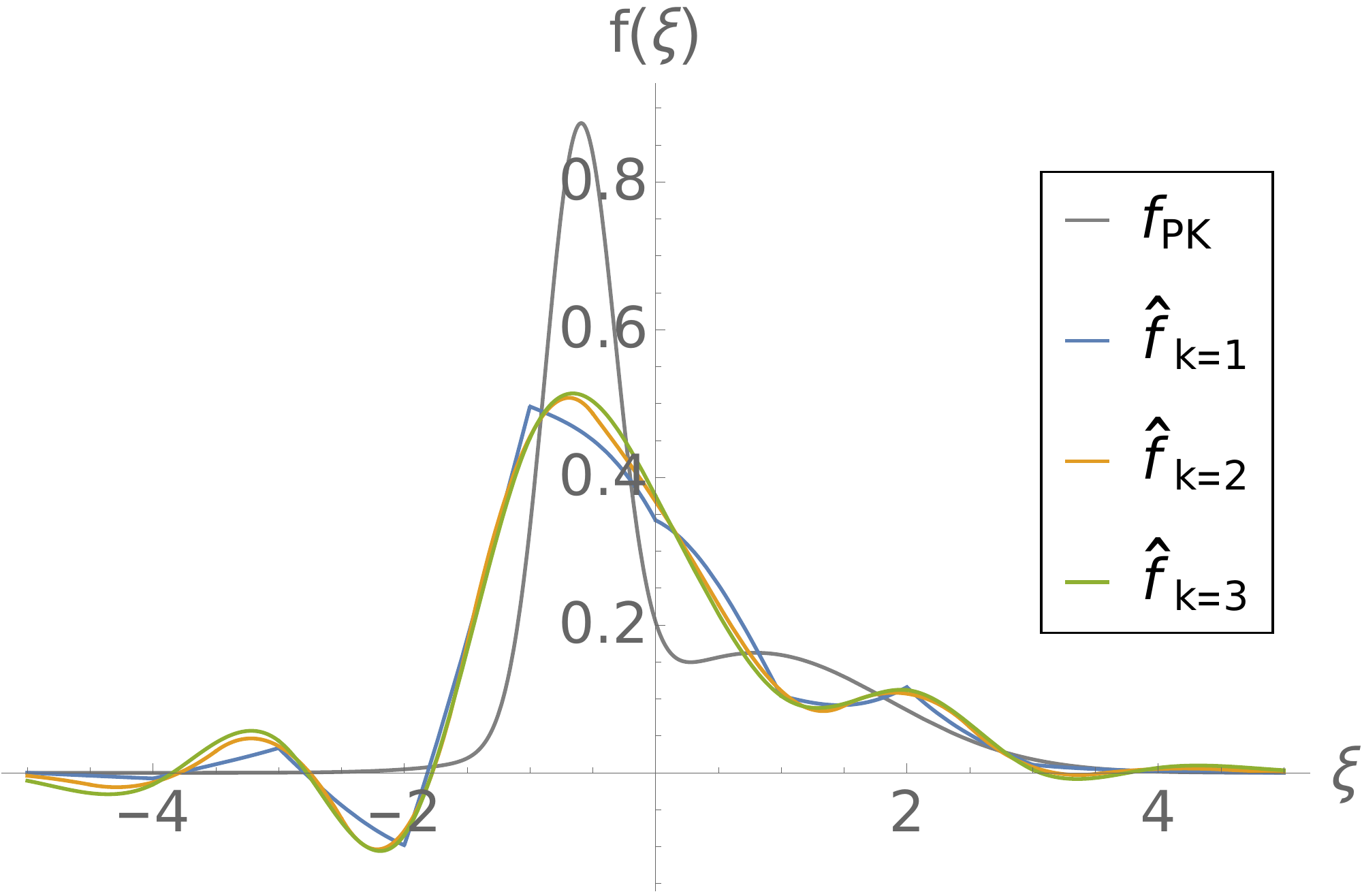}
         \caption{$n=9$.}
         \label{unweightedb/pk1}
\end{subfigure}
\hfill
\begin{subfigure}[b]{0.47\textwidth}
         \centering
         \includegraphics[width=\textwidth]{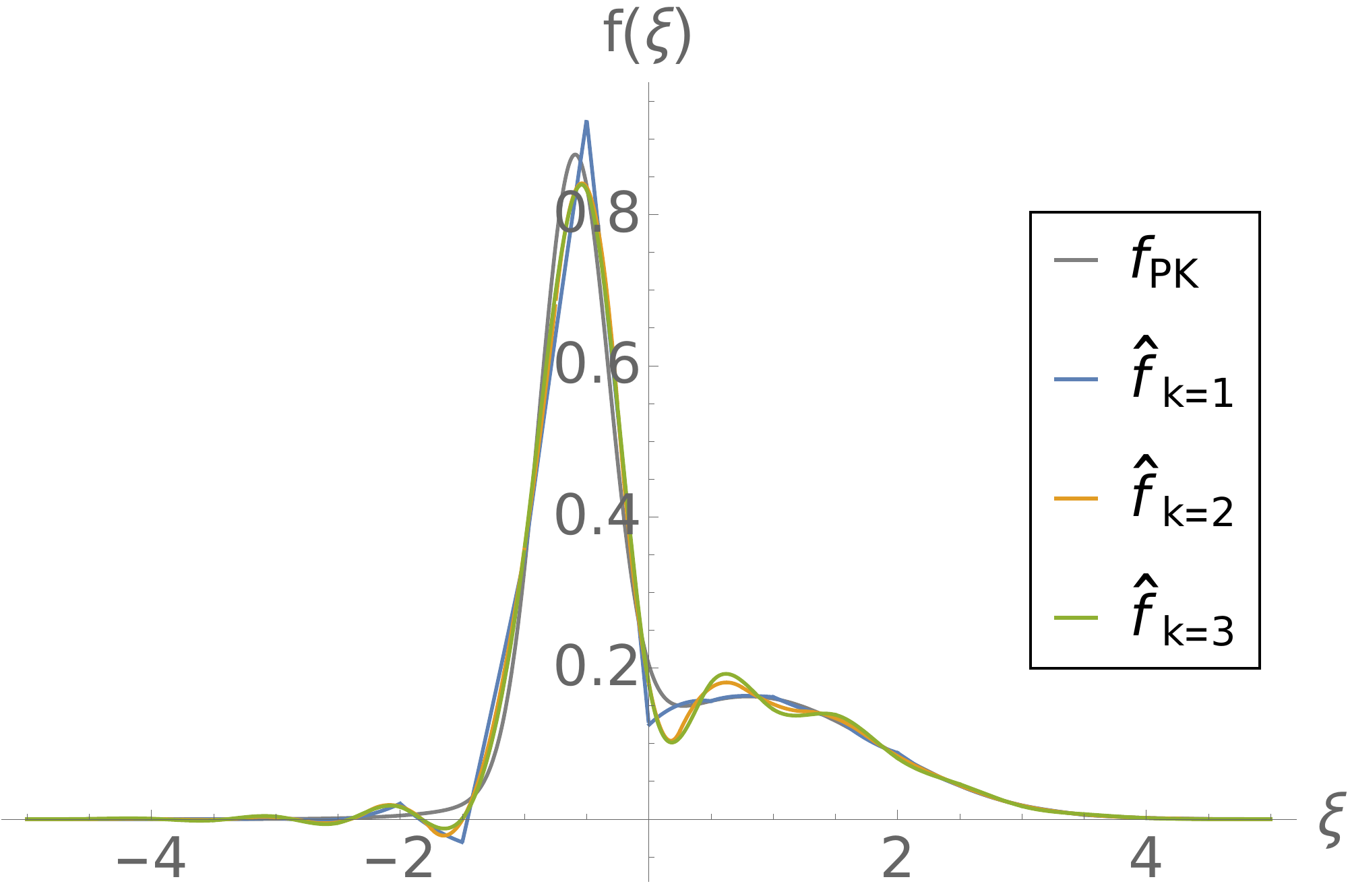}
         \caption{$n=17$.}
         \label{unweightedb/pk2}
     \end{subfigure}
\hfill
\begin{subfigure}[b]{0.47\textwidth}
         \centering
         \includegraphics[width=\textwidth]{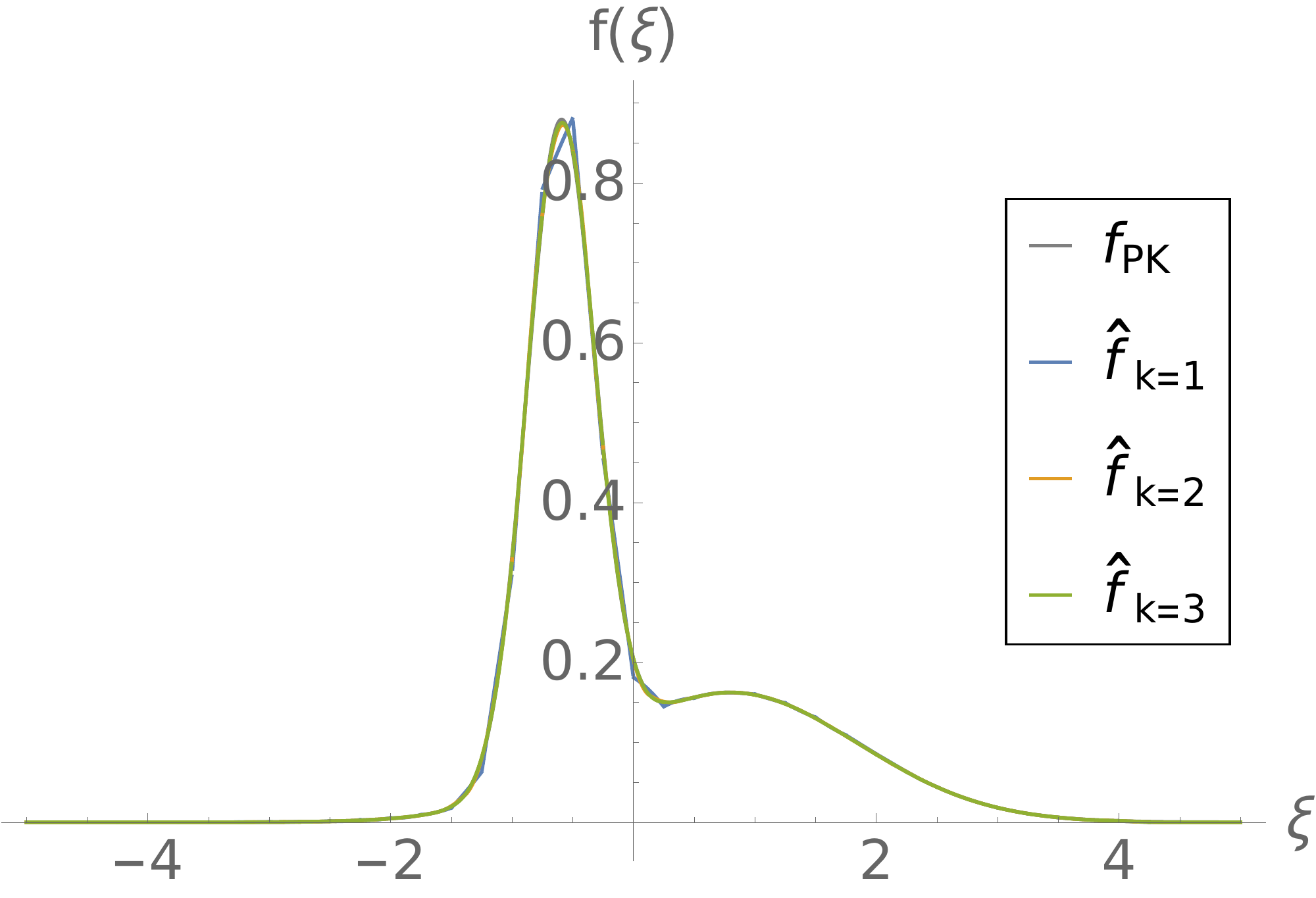}
         \caption{$n=33$.}
         \label{unweightedb/pk3}
     \end{subfigure}
\hfill
\begin{subfigure}[b]{0.47\textwidth}
         \centering
         \includegraphics[width=\textwidth]{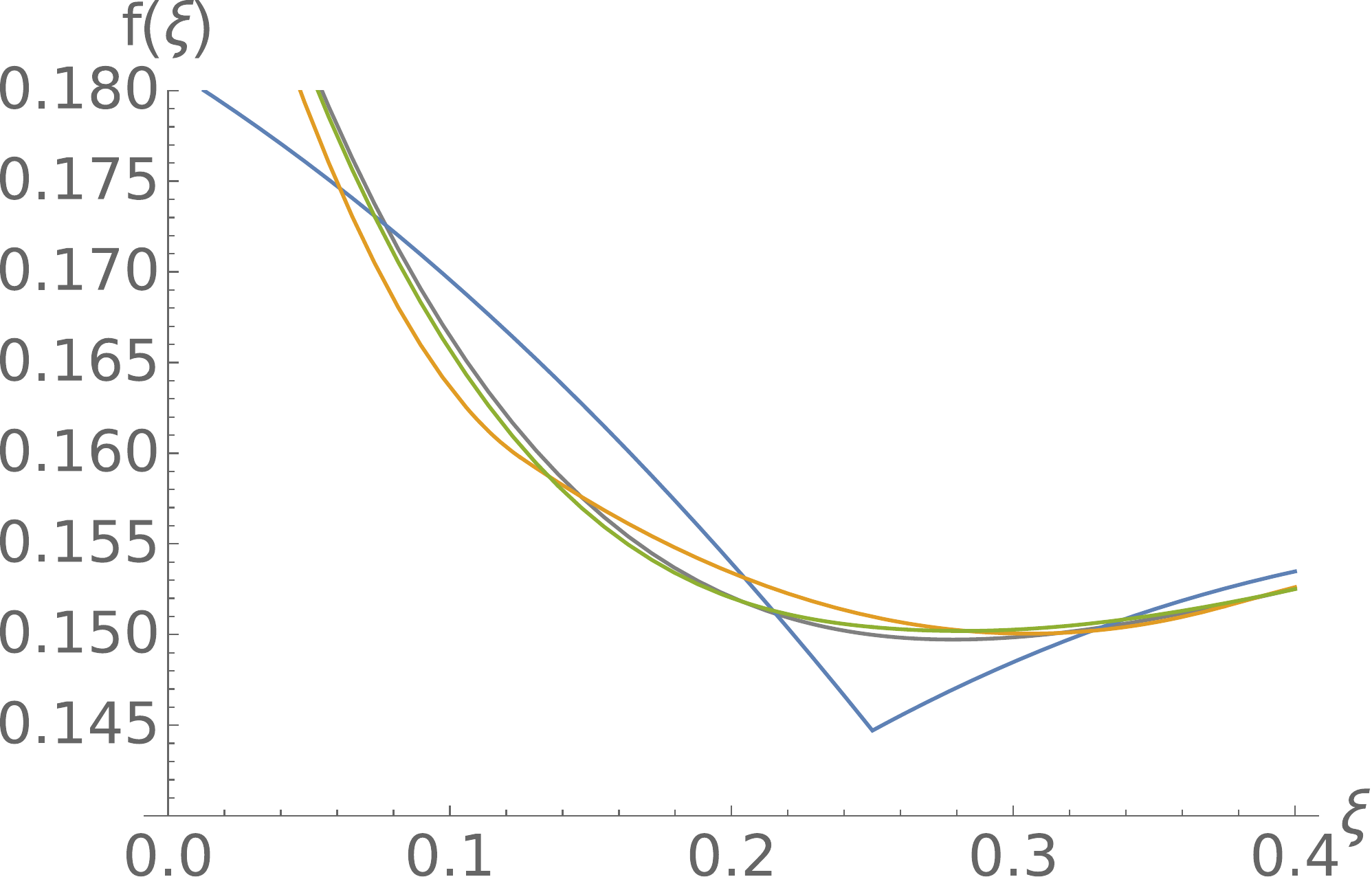}
         \caption{Zoom into (c). $n=33$.}
         \label{unweightedb/pk4}
\end{subfigure}
\caption{Approximation of $f_{PK}$ with unweighted B-splines of different orders $k$, $[\xi_{min},\xi_{max}]=[-4,4]$. Number of splines $n$ ranging from (a) $n=9$, (b) $n=17$, (c)-(d) $n=33$.}
\label{approxUnweightedBPK}
\end{figure}

The unweighted B-spline ansatz was tested in a numerical study for the three bimodal distribution functions $f_{PK}$, $f_{b1}$ and $f_{b2}$ for different numbers of splines $n$ and various spline orders $k$ for constant $\xi$-range $\xi \in [\xi_{min},\xi_{max}]=[-4,4]$ with results shown in Figure \ref{approxUnweightedBPK}.

The approximation error decreases with the number of splines. We can also observe that B-splines of higher order are superior in approximation, although this has less effect than the number of splines.
Figure \ref{unweightedb/pk4}, which is an enlarged view of Figure \ref{unweightedb/pk3}, shows that splines of higher order fit the curves of the original function better and are thus more precise.

\begin{figure}[ht]
\centering
\begin{subfigure}[b]{0.48\textwidth}
         \centering
         \includegraphics[width=\textwidth]{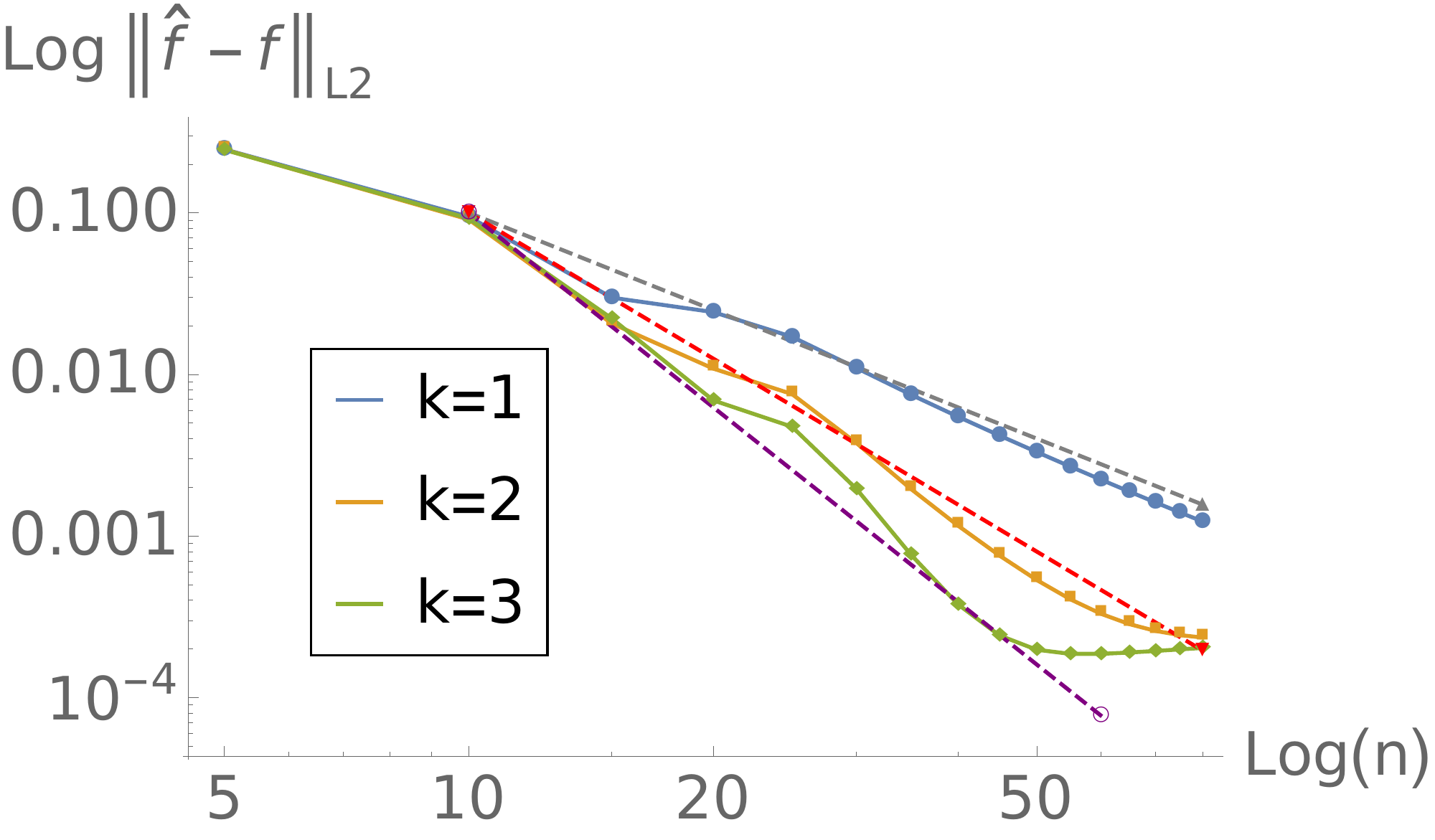}
         \caption{$[\xi_{min},\xi_{max}]=[-4,4]$.}
         \label{unweightedb/plot}
\end{subfigure}
\hfill
\begin{subfigure}[b]{0.48\textwidth}
         \centering
         \includegraphics[width=\textwidth]{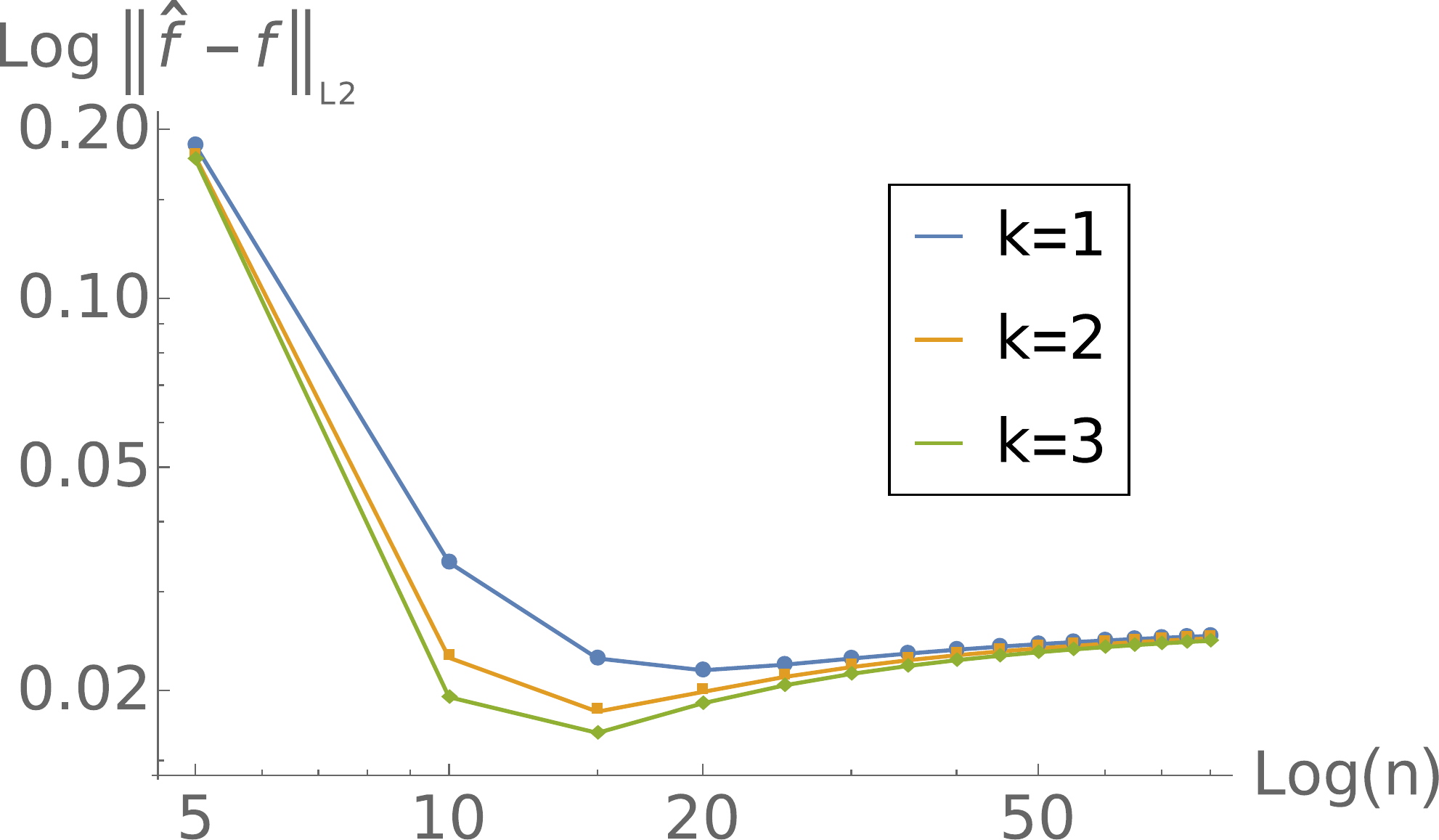}
         \caption{$[\xi_{min},\xi_{max}]=[-2,2]$.}
         \label{unweightedb/plot_smol}
\end{subfigure}
\caption{Approximation error depending on $n$ for unweighted B-splines and two different $\xi$-ranges. (a) $[\xi_{min},\xi_{max}]=[-4,4]$, (b) $[\xi_{min},\xi_{max}]=[-2,2]$.}
 \label{unweightedb}
\end{figure}

The relation between the number of splines and the error $(\Delta f):=\|f-\hat f\|_{L2}$ for the three orders $k=1,2,3$ is displayed in Figure \ref{unweightedb/plot}. The error curves describe the arithmetically averaged relative error of all three distribution functions. The plot also shows lines with convergence rates $1,2,3$, corresponding to the (empirical) rate of convergence following the approximations with the respective orders $k=1,2,3$. Note that from $50$ splines onwards there is no further improvement for the splines of third order. This can be attributed to the fact that the splines have finite support, whereas the bimodal function's support is the whole $\xi$-axis. A small error of $10^{-4}$ thus remains. In order to clearly demonstrate this effect we consider the error on an even smaller $\xi$-range $[\xi_{min},\xi_{max}]=[-2,2]$ in Figure \ref{unweightedb/plot_smol}. Interestingly, a reduced $\xi$-range leads to better approximation results for a small number of spline functions. This is because for an unchanged number of splines, the B-splines now lie much denser ($\Delta \xi$ has decreased) and approximate the function better on that small domain around the origin. However, the limit is reached very soon and an error of $0.025$ remains. From $20$ splines onwards the approximation on the range $[-4,4]$ is clearly superior and reaches lower error values.
We will investigate this effect further when comparing with other ansatzes.


\subsection{{Weighted B-splines}}
Weighted B-splines are multiplied with the Gaussian distribution and allow for a smoother distribution function with more weight around the center.
\begin{definition}[weighted B-spline expansion]
 Let $k \in \naturals$ be the order and $n \in \naturals$ be a fixed number of splines. Let $G_{\xi}$ be a standard grid using a $\xi$-range $[\xi_{min},\xi_{max}]$ and equidistant grid spacing $\Delta \xi=\frac{(\xi_{max}-\xi_{min})}{(n-1)}$.

    The \emph{weighted B-spline} expansion of a distribution function $f: \reals \mapsto \reals$ is then given by
    \begin{equation}
        \hat{f}(\xi)=\maboxinormal \left( 1+ \sum_{i=1}^{n} \alpha_i \phi_i(\xi)\right).
    \end{equation}
    for B-splines $\phi_i = B_{ki}$ of order $k$.
\end{definition}

The approximation of the three bimodal functions $f_{PK}$, $f_{b1}$ and $f_{b2}$ can then be computed using Galerkin projection with B-spline test functions $\psi_j=B_{jk}$. The coefficient vector $\alpha$ is determined according to the method presented in Subsection \ref{vorgehen}.

\begin{figure}[H]
\centering
\begin{subfigure}[b]{0.47\textwidth}
         \centering
         \includegraphics[width=\textwidth]{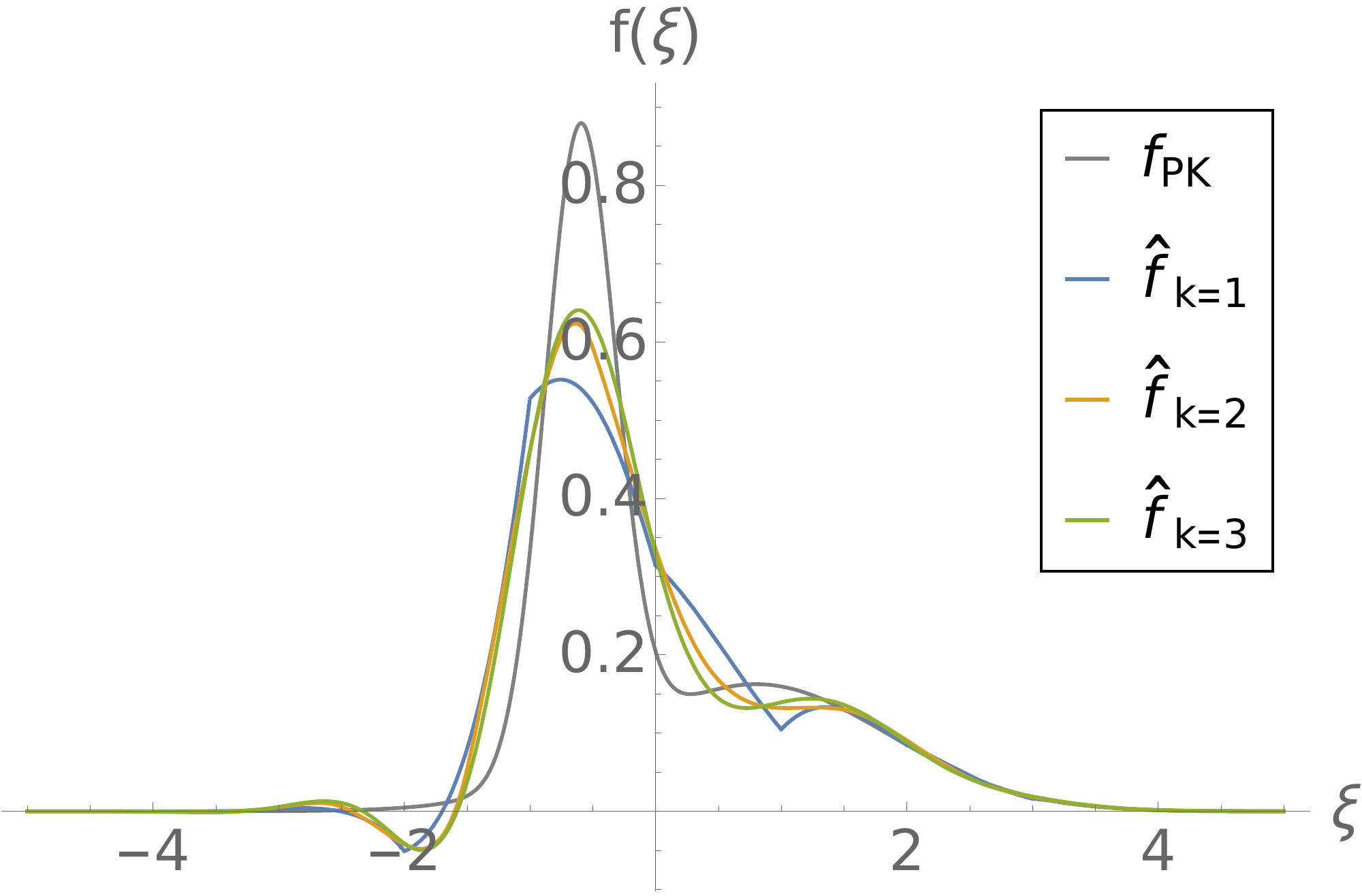}
         \caption{$n=9$.}
         \label{weightedb/pk1}
\end{subfigure}
\hfill
\begin{subfigure}[b]{0.47\textwidth}
         \centering
         \includegraphics[width=\textwidth]{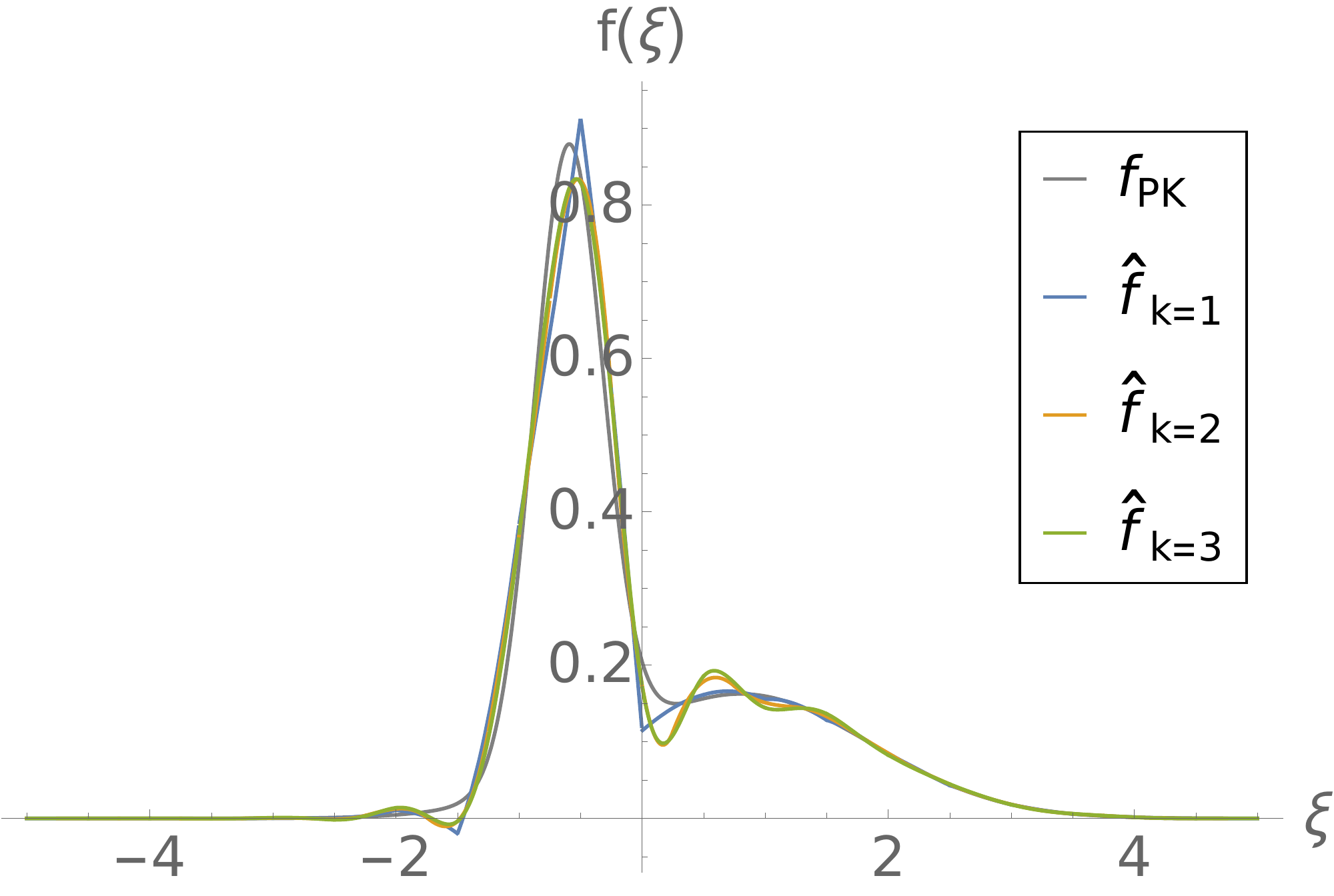}
         \caption{$n=17$.}
         \label{ightedb/pk2}
     \end{subfigure}
\caption{Approximation of $f_{PK}$ with weighted B-splines of different orders $k$, $[\xi_{min},\xi_{max}]=[-4,4]$. Number of splines $n$ ranging from (a) $n=9$, (b) $n=17$.}
\label{approxWeighted}
\end{figure}

\begin{figure}[H]
\centering
\begin{subfigure}[b]{0.46\textwidth}
         \centering
         \includegraphics[width=\textwidth]{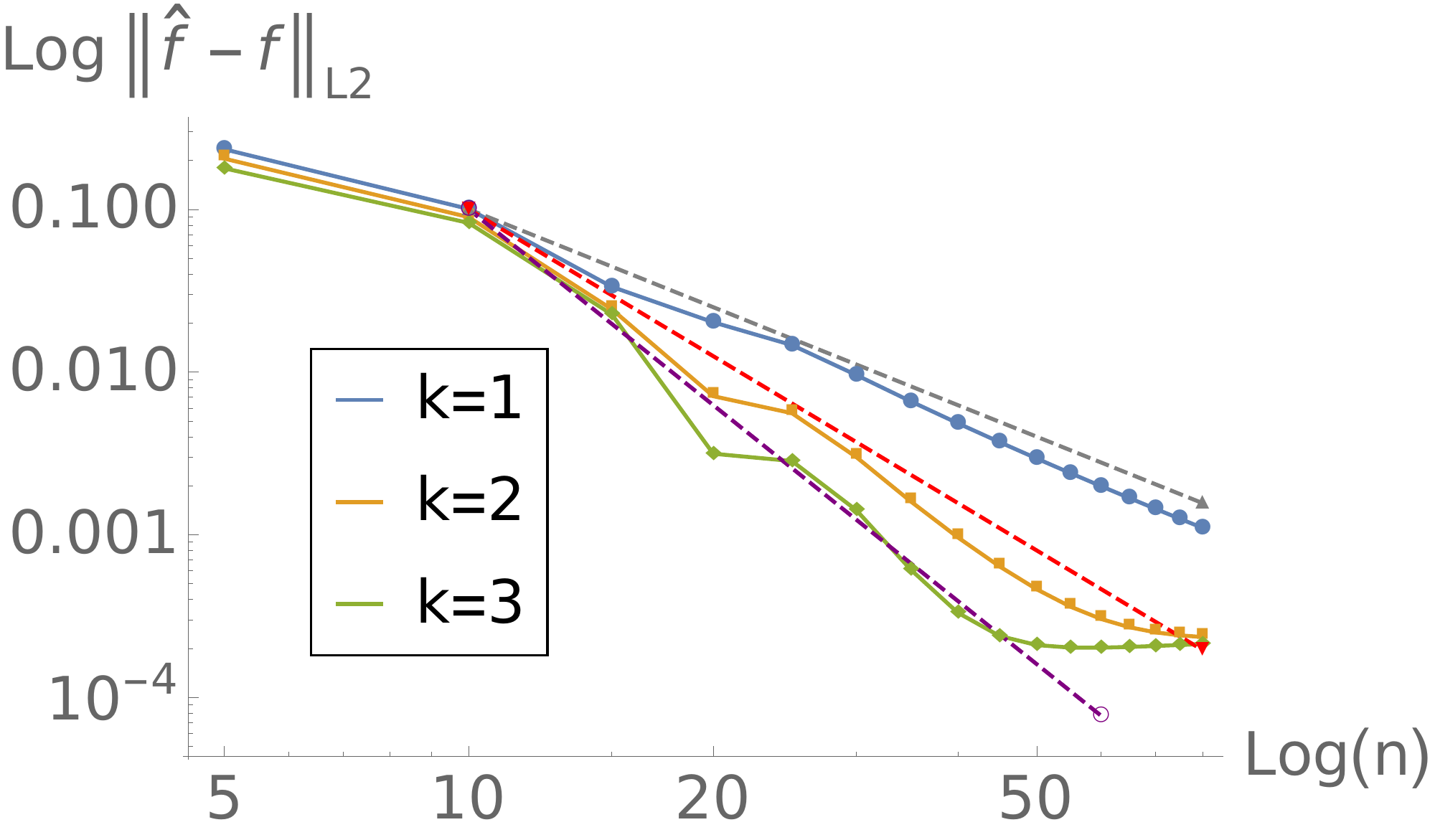}
         \caption{$[\xi_{min},\xi_{max}]=[-4,4]$.}
         \label{weightedb/plot}
\end{subfigure}
\hfill
\begin{subfigure}[b]{0.46\textwidth}
         \centering
         \includegraphics[width=\textwidth]{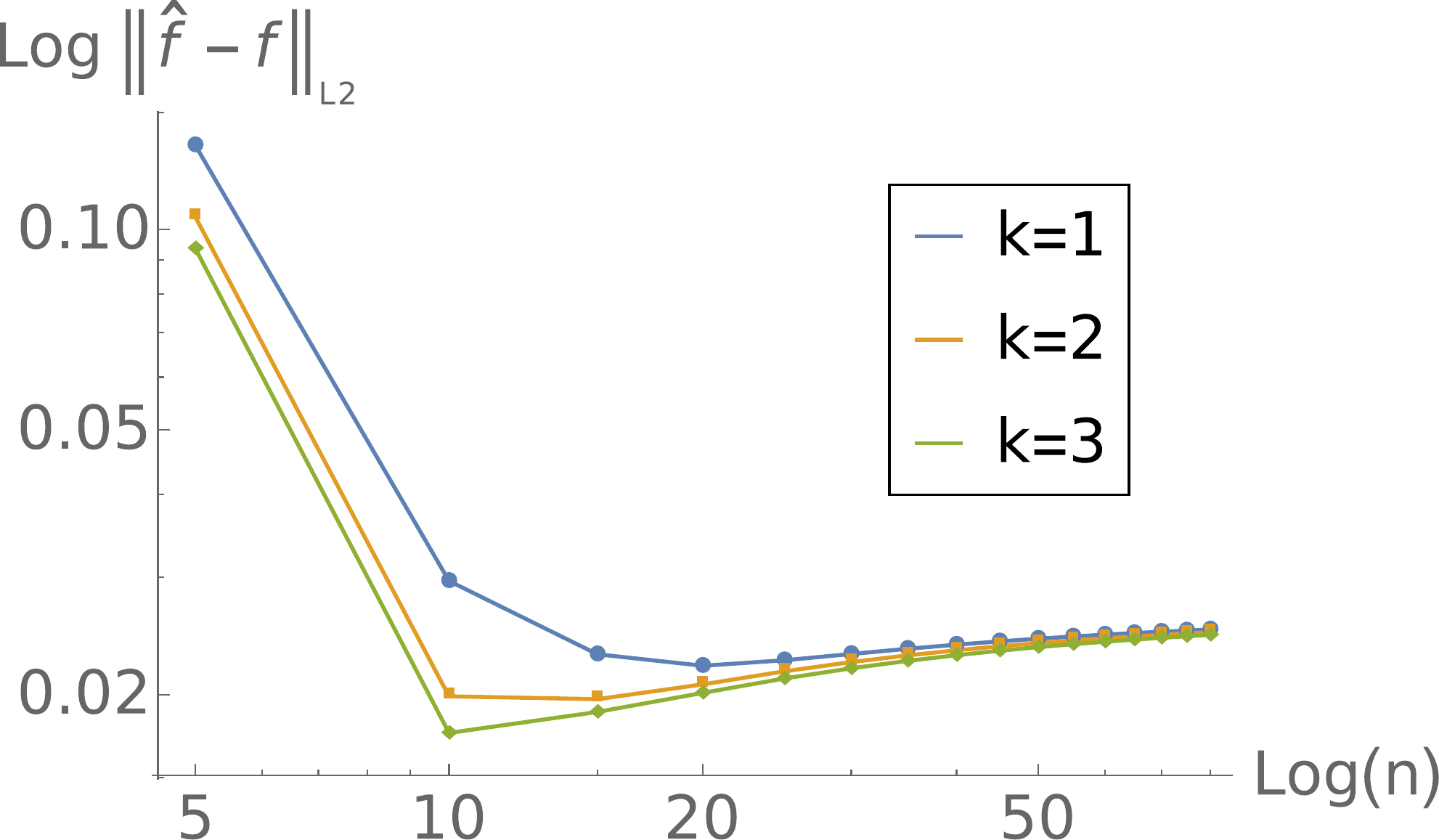}
         \caption{$[\xi_{min},\xi_{max}]=[-2,2]$.}
         \label{weightedb/plot_smol}
\end{subfigure}
       \caption{Approximation error depending on $n$ for weighted B-splines and two different $\xi$-ranges. (a) $[\xi_{min},\xi_{max}]=[-4,4]$, (b) $[\xi_{min},\xi_{max}]=[-2,2]$.}
 \label{weightedb}
\end{figure}

Figure \ref{approxWeighted} shows that there are fewer oscillations close to the boundary of the grid in comparison to the unweighted ansatz (compare for example Figure \ref{weightedb/pk1} and Figure \ref{unweightedb/pk1}). Apart from that the approximation looks equally accurate. The error reduction and the behavior for a smaller grid size in Figure \ref{weightedb} are the same as for the unweighted ansatz.

\subsection{{Weighted fundamental constrained splines}}
Weighted FCS are used to conserve mass, momentum, and energy of the distribution function by using the previously developed FCS basis function plus Gaussian weight.
\begin{definition}[weighted FCS expansion]
    Let $n \in \naturals$ be a fixed number. Let $k \in \naturals$ be the fixed order and $G_{\xi}$ be a standard grid using a $\xi$-range $[\xi_{min},\xi_{max}]$ and equidistant grid spacing $\Delta \xi=\frac{(\xi_{max}-\xi_{min})}{(n-1)}$.

    The \emph{weighted fundamental constrained splines} expansion of a distribution function $f: \reals \mapsto \reals$ is then given by
    \begin{equation}
        \hat{f}(\xi)=\maboxinormal \left( 1+ \sum_{i=1}^{n-3} \alpha_i \phi_i(\xi)\right)
    \end{equation}
    with the $n-3$ fundamental constrained splines $\phi_i$, that are created from the $n$ B-splines of order $k$ on the Grid $G_\xi$.
\end{definition}

\begin{remark}
    By construction of the fundamental constrained splines, the following compatibility conditions are exactly fulfilled
    \begin{equation}
    \label{compat}
     \inftyint \hat{f}(\xi)dx=1,\quad
     \inftyint \xi \hat{f}(\xi)dx=0, \quad \textrm{and} \quad
     \inftyint \xi^2\hat{f}(\xi)=1,
    \end{equation}
    which equally hold for the bimodal functions $f_{PK}$, $f_{b1}$ and $f_{b2}$.
\end{remark}

\begin{figure}[htb!]
\centering
\begin{subfigure}[b]{0.45\textwidth}
         \centering
         \includegraphics[width=\textwidth]{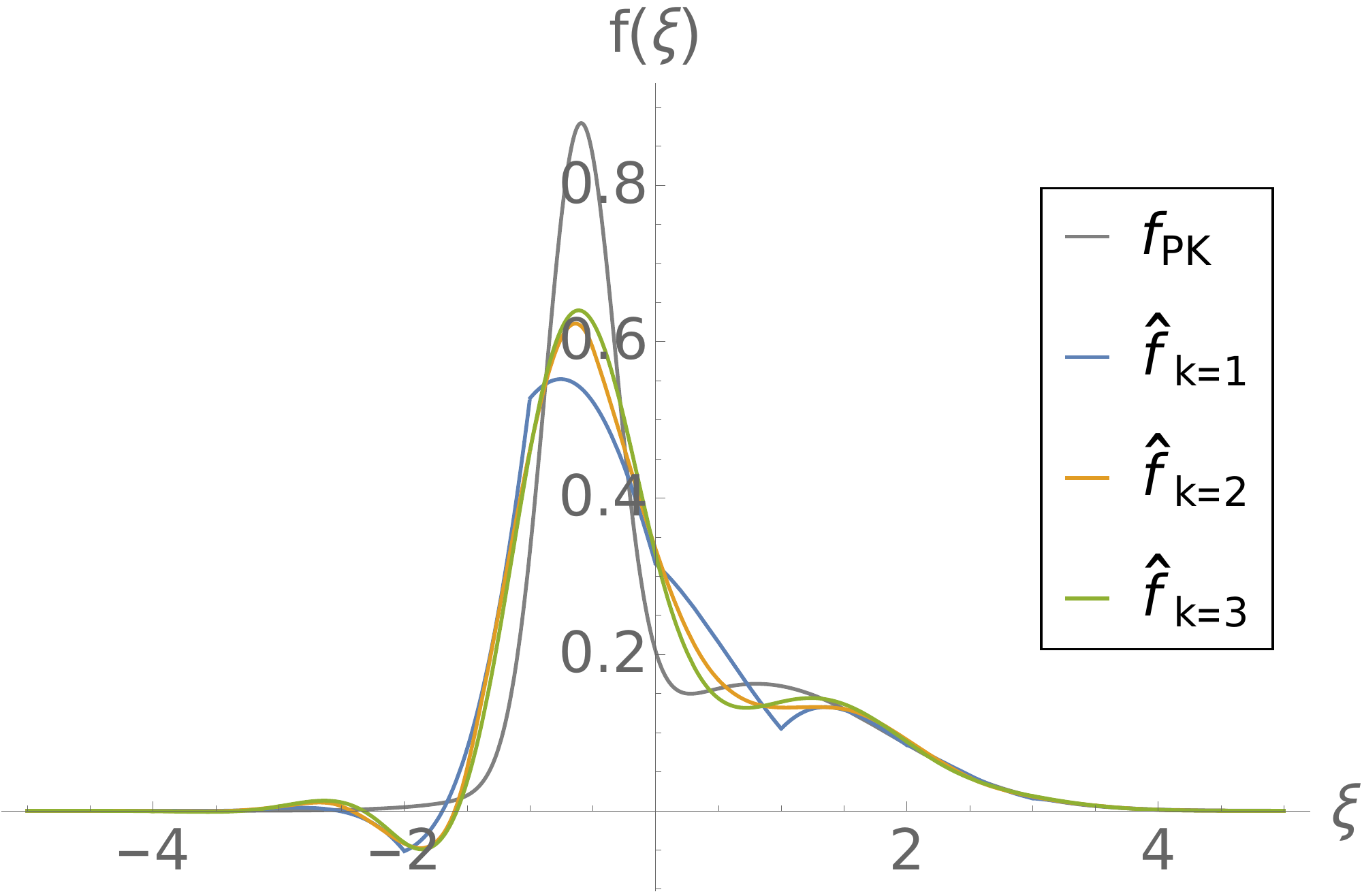}
         \caption{$n=9$.}
         \label{weightedfcs/pk1}
\end{subfigure}
\hfill
\begin{subfigure}[b]{0.45\textwidth}
         \centering
         \includegraphics[width=\textwidth]{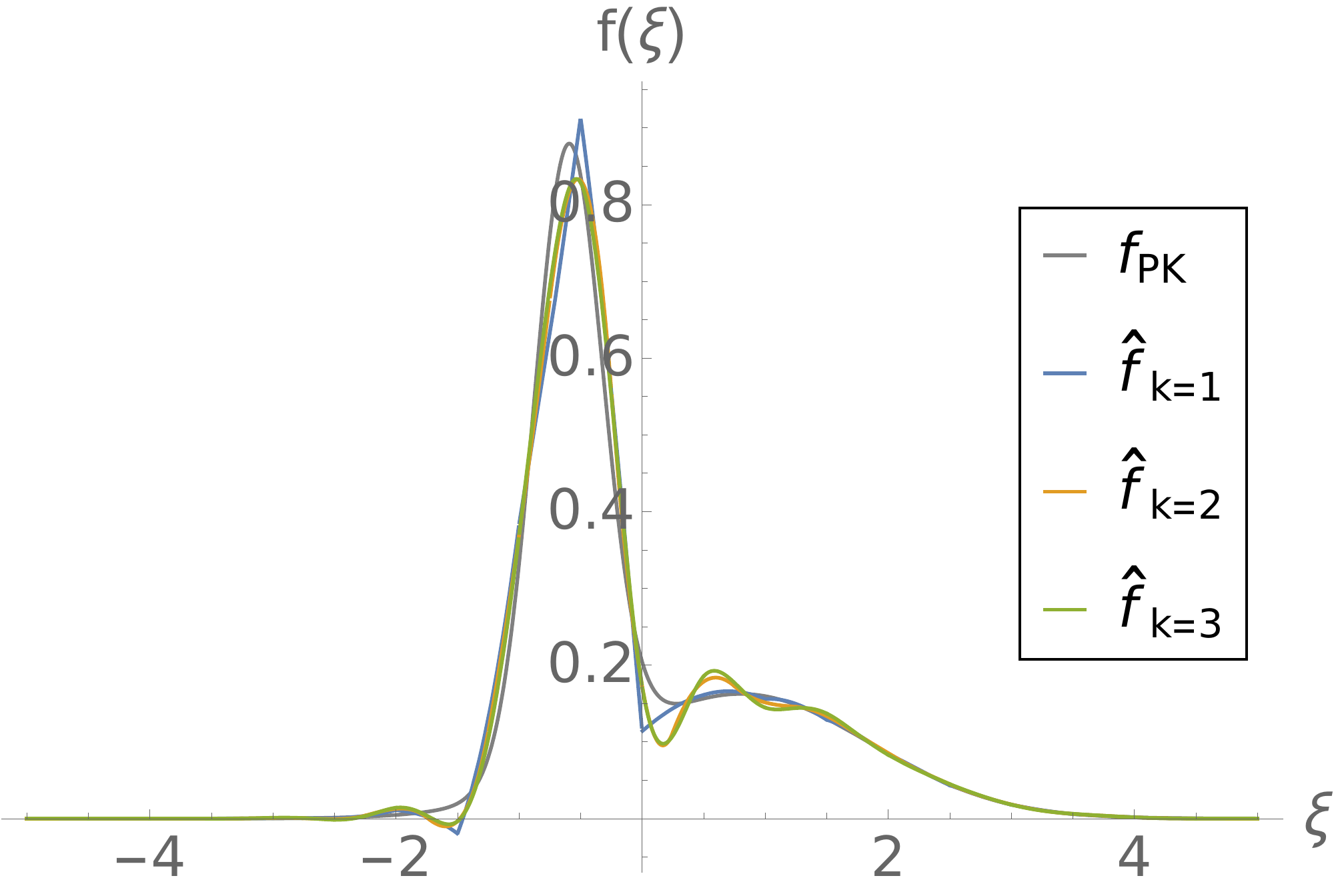}
         \caption{$n=17$.}
         \label{weightedfcs/pk2}
     \end{subfigure}
\hfill
\caption{Approximation of $f_{PK}$ with weighted FCS of different orders $k$, $[\xi_{min},\xi_{max}]=[-4,4]$. Number of splines $n$ ranging from (a) $n=9$, (b) $n=17$.}
\label{approxWeightedFCS}
\end{figure}

\begin{figure}[H]
\centering
\begin{subfigure}[b]{0.47\textwidth}
         \centering
         \includegraphics[width=\textwidth]{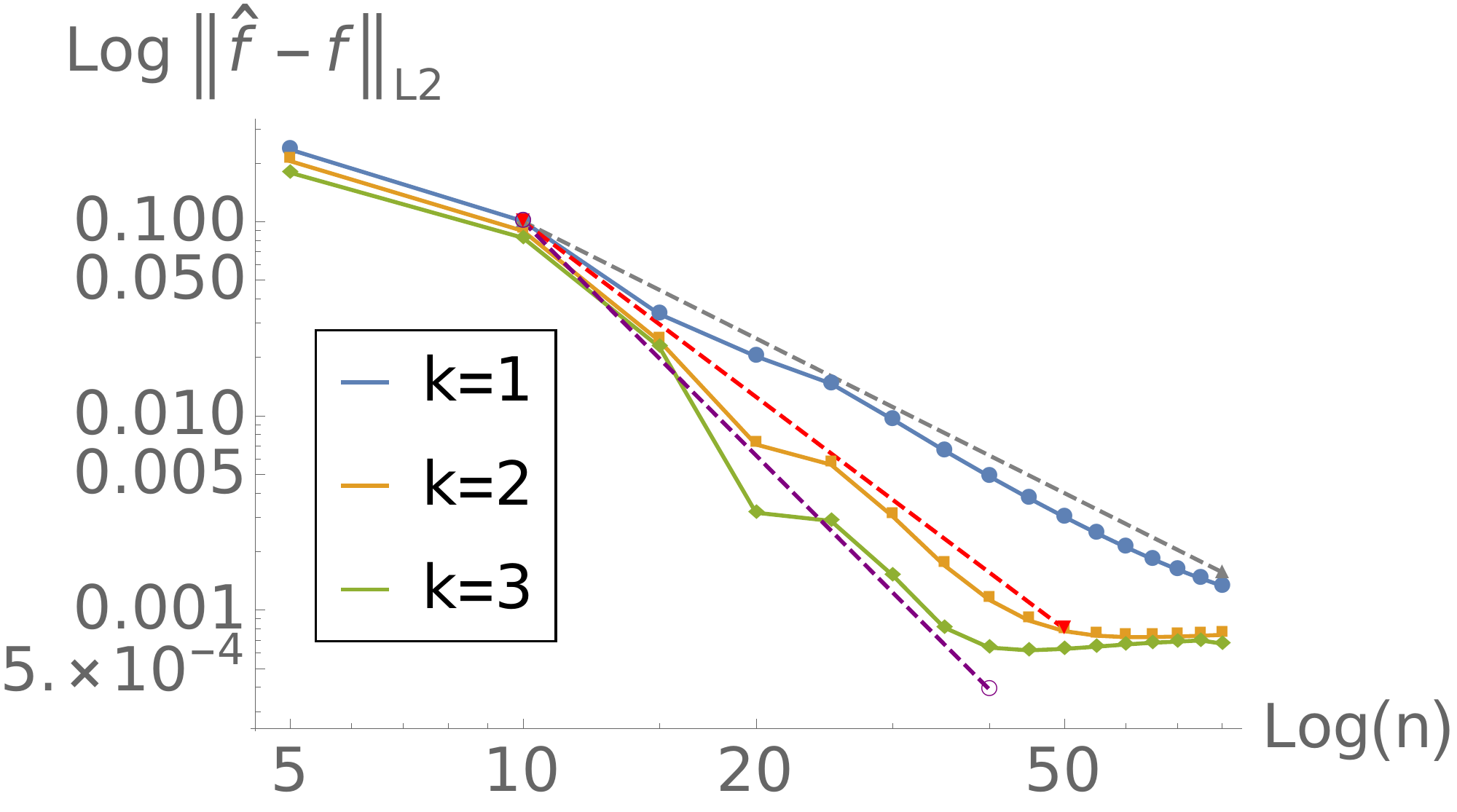}
         \caption{$[\xi_{min},\xi_{max}]=[-4,4]$.}
         \label{weightedfcs/plot}
\end{subfigure}
\hfill\begin{subfigure}[b]{0.47\textwidth}
         \centering
         \includegraphics[width=\textwidth]{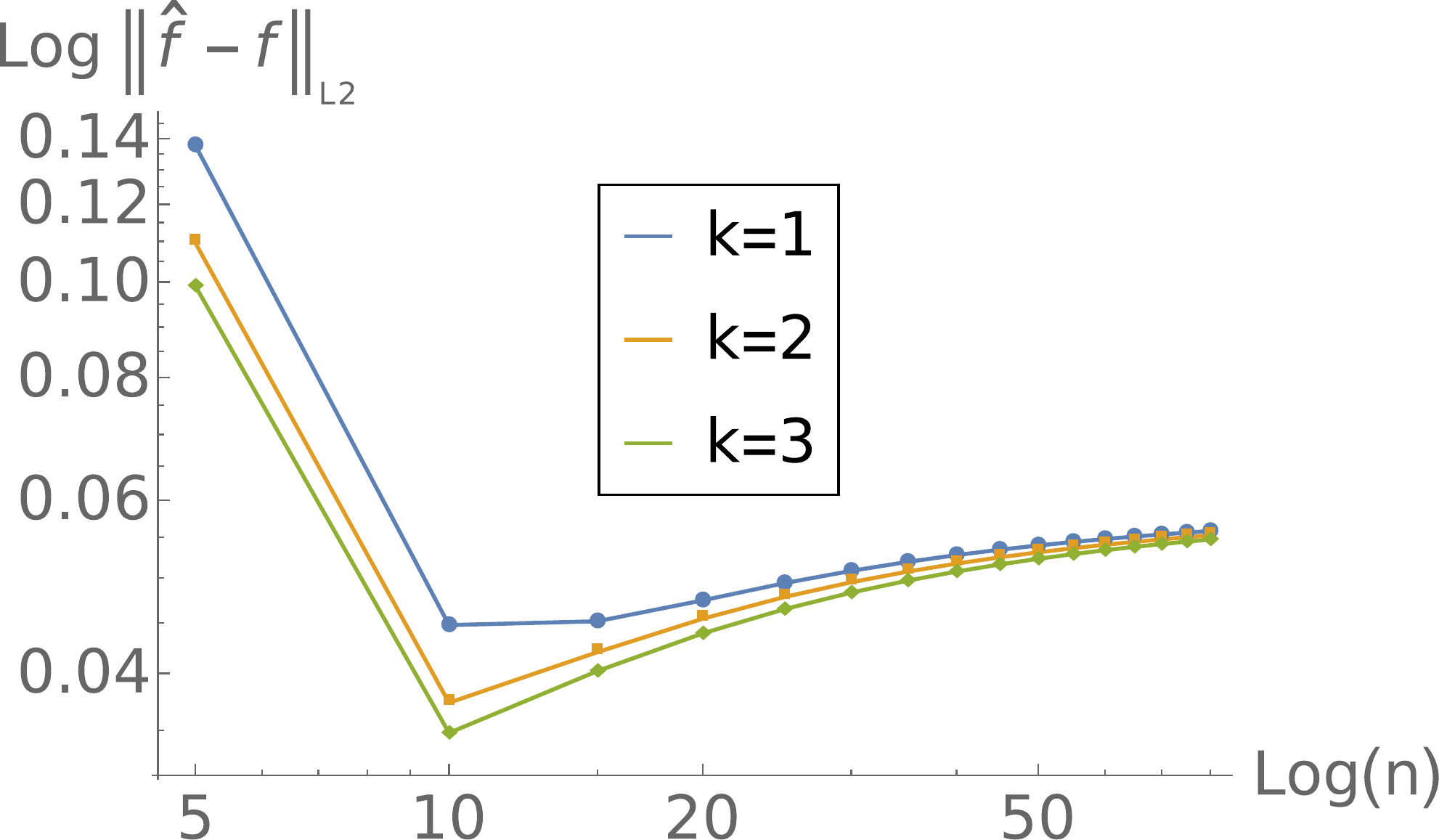}
         \caption{$[\xi_{min},\xi_{max}]=[-2,2]$.}
         \label{weightedfcs/plot_smol}
\end{subfigure}
       \caption{Approximation error depending on $n$ for weighted FCS and two different $\xi$-ranges. (a) $[\xi_{min},\xi_{max}]=[-4,4]$, (b) $[\xi_{min},\xi_{max}]=[-2,2]$.}
 \label{weightedfcs}
\end{figure}

 Figure \ref{approxWeightedFCS} again shows fewer oscillations around the boundaries of the grid due to the weighting. A similar decrease of the approximation error with increasing number of splines is seen. Note that the number of splines does not correspond to $n$ but to $n-3$ by construction of the fundamental constrained splines. There is some residual error due to the finite $\xi$-range as shown in Figure \ref{weightedfcs/plot}. We note that when extending the limits to $[\xi_{min}=-6,\xi_{max}=6]$, the error saturation happens at much lower error values, at around $10^{-6}$ in contrast to the saturation at $10^{-4}$ that we observe on the range $[-4,4]$. This result is not pictured for conciseness.

\subsection{Approximation property}
The approximation properties of B-splines for standard continuous functions are well-known and can be found in the literature, for example, in \cite{bookPGS}, Chapter XII. Therein, the following theorem for the approximation error can be found, which reads for an equidistant grid:

\begin{theorem}[Spline approximation error]
    \label{th:Jackson Type}
    Let $\$_{k,t}$ denote the space of all B-splines of order $k$ on the equidistant grid $G: \, a=t_0 < t_1 < t_2 < \hdots < t_n = b$ with mesh size $|t|$ and $\| \cdot \|_{max}:=\max_{a\leq x \leq b}|g(x)|$.
    Then for each order $k$ there exists a constant $C_{k}$ so that for all smooth functions $g \in C^{(k+1)}[a\hdots b]$
    $$dist(g,\$_{k,t}) \coloneqq min\{\| g - s \|_{max} : s \in \$_{k,t}\} \leq C_k |t|^{k+1} \|g^{(k+1)}\|_{max}\;.$$
\end{theorem}

Theorem \ref{th:Jackson Type} states that the approximation error decreases at least with the $k$th power of the mesh size, corresponding to our observations regarding the convergence speed, where we observed the respective higher rates of convergence for splines of higher order. With the additional restriction that the function $g$ from \ref{th:Jackson Type} fulfills the compatibility conditions in Equation \ref{compat}, we can formulate the same theorem for the weighted FCS.

\subsection{{Spline approximation convergence}}
\begin{figure}[H]
\centering
\begin{subfigure}[b]{0.49\textwidth}
         \centering
         \includegraphics[width=\textwidth]{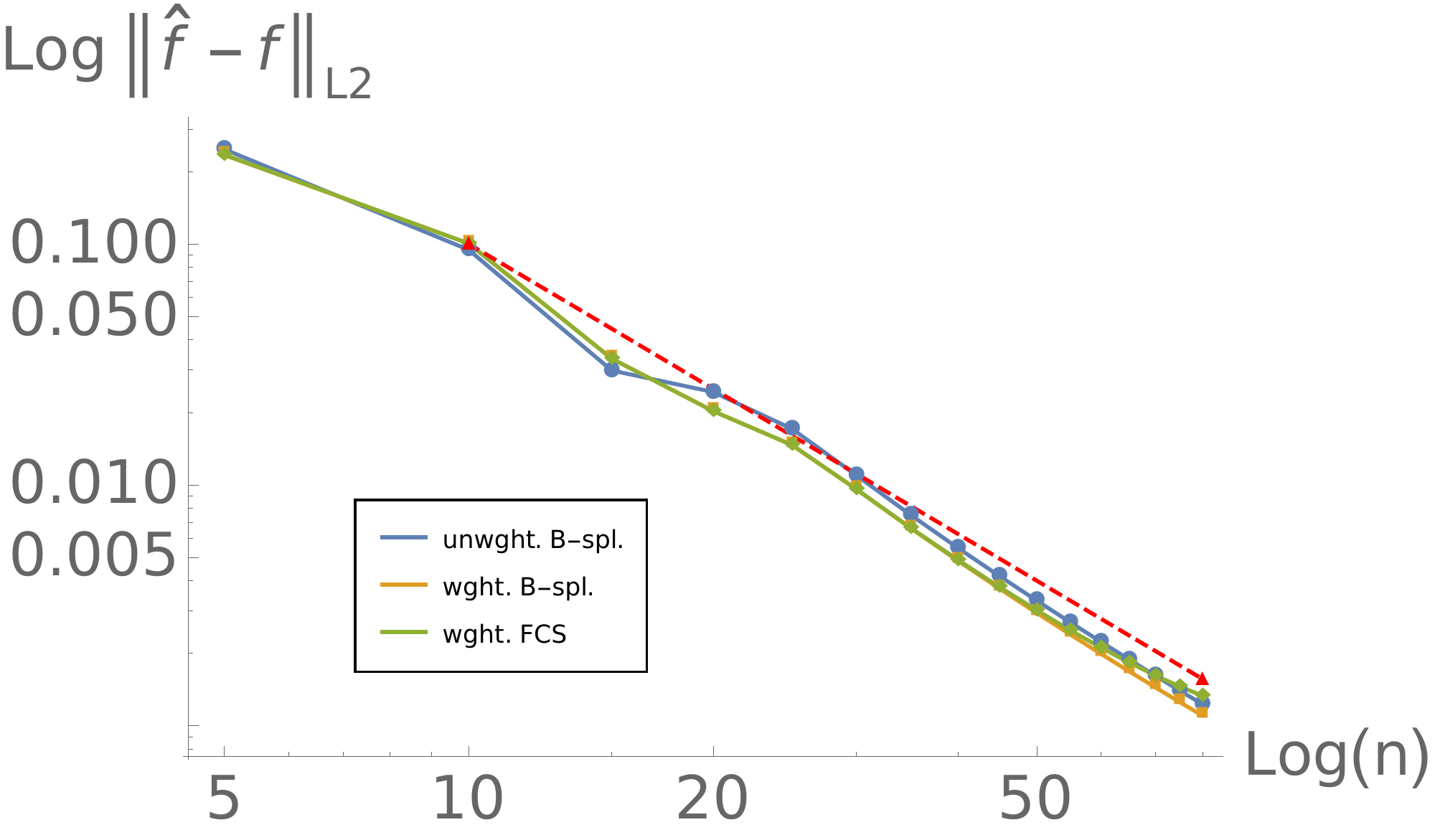}
         \caption{$k=1$.}
         \label{cmpLin}
\end{subfigure}
\hfill
\hfill
\begin{subfigure}[b]{0.49\textwidth}
         \centering
         \includegraphics[width=\textwidth]{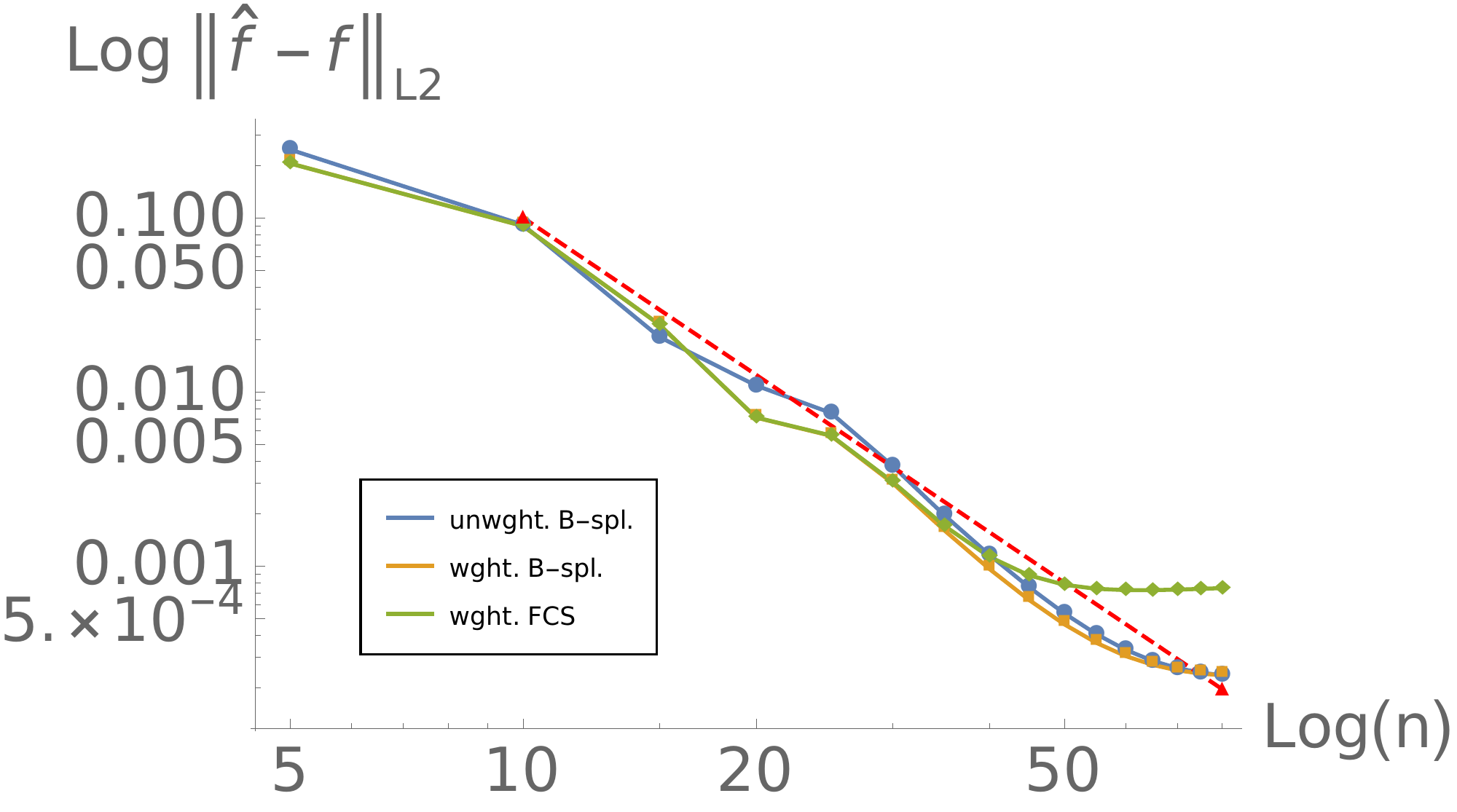}
         \caption{$k=2$.}
         \label{cmpQuad}
\end{subfigure}
\hfill
\hfill
\begin{subfigure}[b]{0.49\textwidth}
         \centering
         \includegraphics[width=\textwidth]{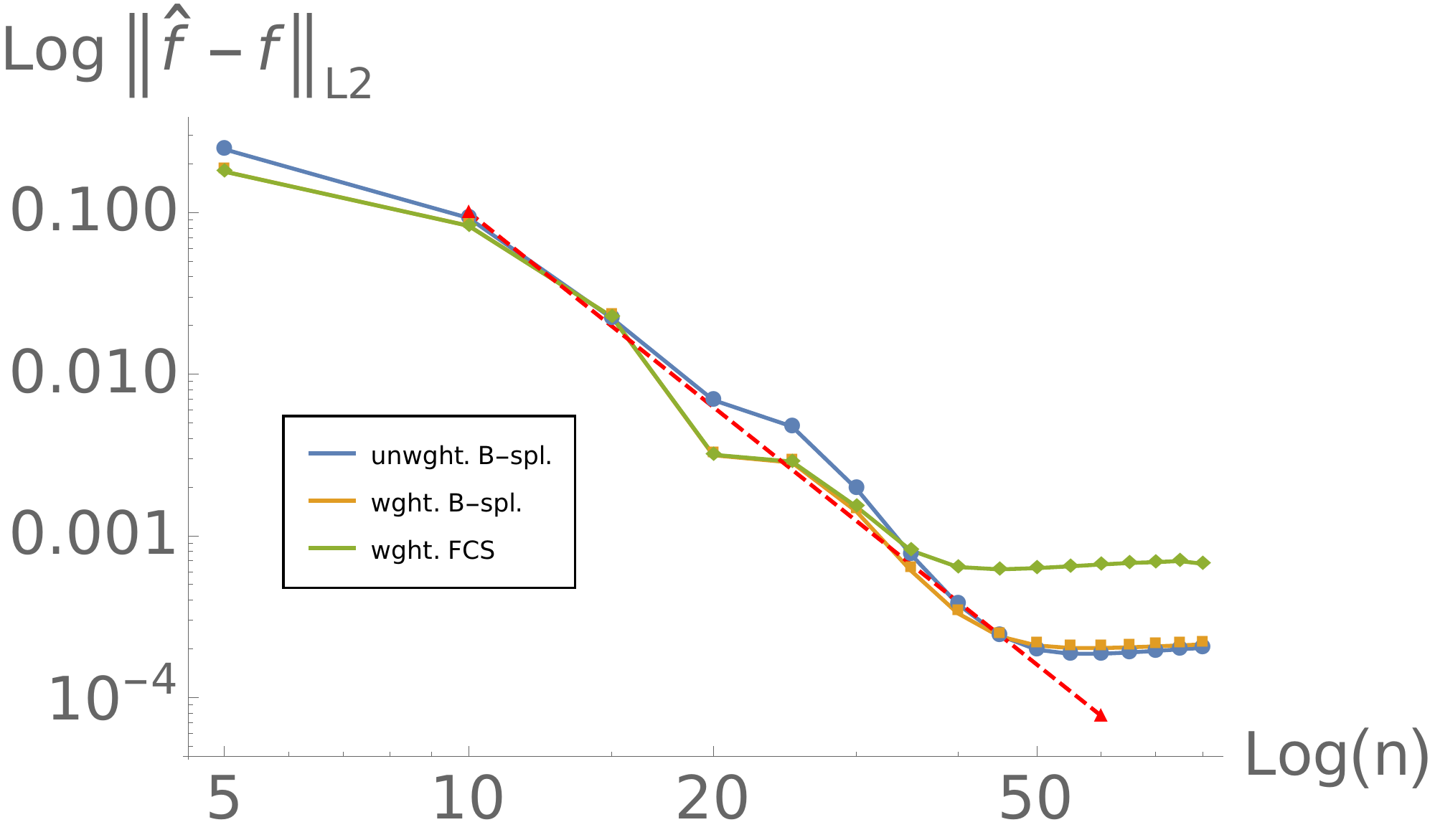}
         \caption{$k=3$.}
         \label{cmpCubic}
\end{subfigure}
\hfill
\caption{Approximation error depending on the number of underlying B-splines $n$ for different orders $k$ and basis functions, $[\xi_{min},\xi_{max}]=[-4,4]$. (a) $k=1$, (b) $k=2$, (c) $k=3$.}
 \label{cmpConvergence}
\end{figure}

In Figure \ref{cmpConvergence} the respective error evolution of each ansatz is compared using the standard support $[\xi_{min},\xi_{max}]=[-4,4]$. Note that we compare by the number of underlying B-splines $n$ and not by the number of splines which equals $n$ for the weighted and unweighted B-splines and $n-3$ for the FCS. For the linear splines in Figure \ref{cmpLin}, weighted and unweighted B-splines as well as weighted FCS show the same convergence behavior. For all models we observed that the error reaches a certain plateau. The plateaus for the linear case in Figure \ref{cmpLin} are not visible because they are reached beyond $n=50$. For a smaller support, the plateaus are reached earlier, resulting in a larger remaining error (not shown here for conciseness). We also observed that the remaining error for the FCS is slightly higher than the remaining error for the other models. The plateaus are in all cases a combination of the limited support of the spline basis and of rounding errors that occurred during the computation of the coefficients, particularly while solving the linear system \eqref{lgs} using an ill-conditioned matrix. Additional errors for the FCS can occur during the construction of the fundamental constrained splines, i.e. solving the linear system in Equation \ref{fcslgs}. However, in practice we do not expect to use more than $15-20$ splines to balance accuracy and efficiency. In this region, no problems with stability occur and the error remains very small in all tested cases.

In summary, we observe that the approximation of bimodal distribution functions via weighted FCS yields very high accuracy while conserving mass, momentum, and energy of the distribution function. We will thus also use the weighted FCS basis for the derivation of a set of PDEs for the dynamic evolution of the coefficients in the following section.

\begin{remark}
    In this paper, we focus on the one-dimensional Boltzmann-BGK equation. There exist several possibilities to extend the spline basis to a multi-dimensional setting. In principle, a tensorized ansatz with or without certain adaptivity for different directions as outlined in \cite{Koellermeier2018,Torrilhon2015} can be used. According to the tests in this section, a bounded support for the spline basis yields sufficient accuracy in a 1D setup, so that the full multi-dimensional basis might not need to include too many basis functions. This is especially true as we use the spline basis for the transformed Boltzmann-BGK equation introduced in Section \ref{sec:5}. We note that this is a crucial difference to existing models based on DG discretizations in non-adaptive velocity spaces, see e.g. \cite{Zhang2018,Alekseenko2012,Alekseenko2014}. In our case, the resulting moment model will be highly non-linear, but efficient and accurate at the same time.
\end{remark}

\subsection{{Approximation of a discontinuous distribution function}}
\begin{figure}[H]
	\centering
	\includegraphics[scale=0.4]{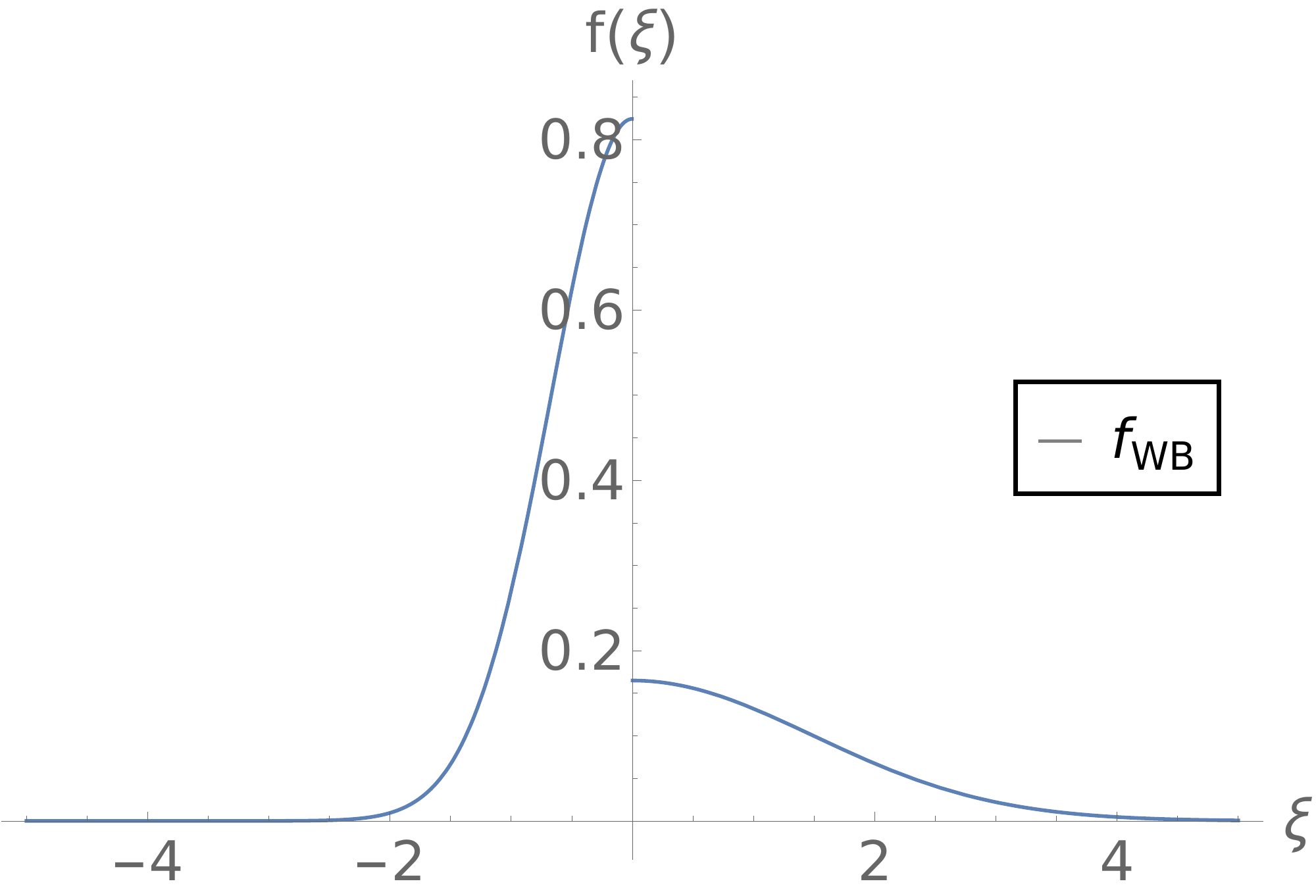}
	\caption{Discontinuous distribution function $f_{WB}$ from a wall boundary condition.}
	\label{fig:wallBoundary}
\end{figure}

The splines of order $k$ are $k$ times continuously differentiable and are thus most often used for the approximation of continuous functions. However, typical distribution functions during the simulation of the Boltzmann-BGK equation do include discontinuities. Especially in near-wall flows, where the wall boundary condition is given by some external temperature. A typical example for such a distribution function is the wall-boundary distribution function
\begin{equation*}
f_{WB}(\xi)=
\begin{cases}
0.824423 e^{-1.11803 \xi ^2} & \xi < 0, \\
0.164885 e^{-0.223607 \xi ^2} & \xi \geq 0, \\
\end{cases}
\end{equation*}
where the numerical values are chosen such that the distribution function has moments $\rho=1$, $v=0$ and $\theta=1$, corresponding to the previous examples in section \ref{sec:4}. The wall-boundary distribution function is commonly encountered in boundary value problems, where the wall has a different temperature than the fluid, for example in the numerical tests in \cite{Fan2019}. It is depicted in Figure \ref{fig:wallBoundary}.

\begin{figure}[htb!]
	\centering
	\begin{subfigure}[b]{0.45\textwidth}
		\centering
		\includegraphics[width=\textwidth]{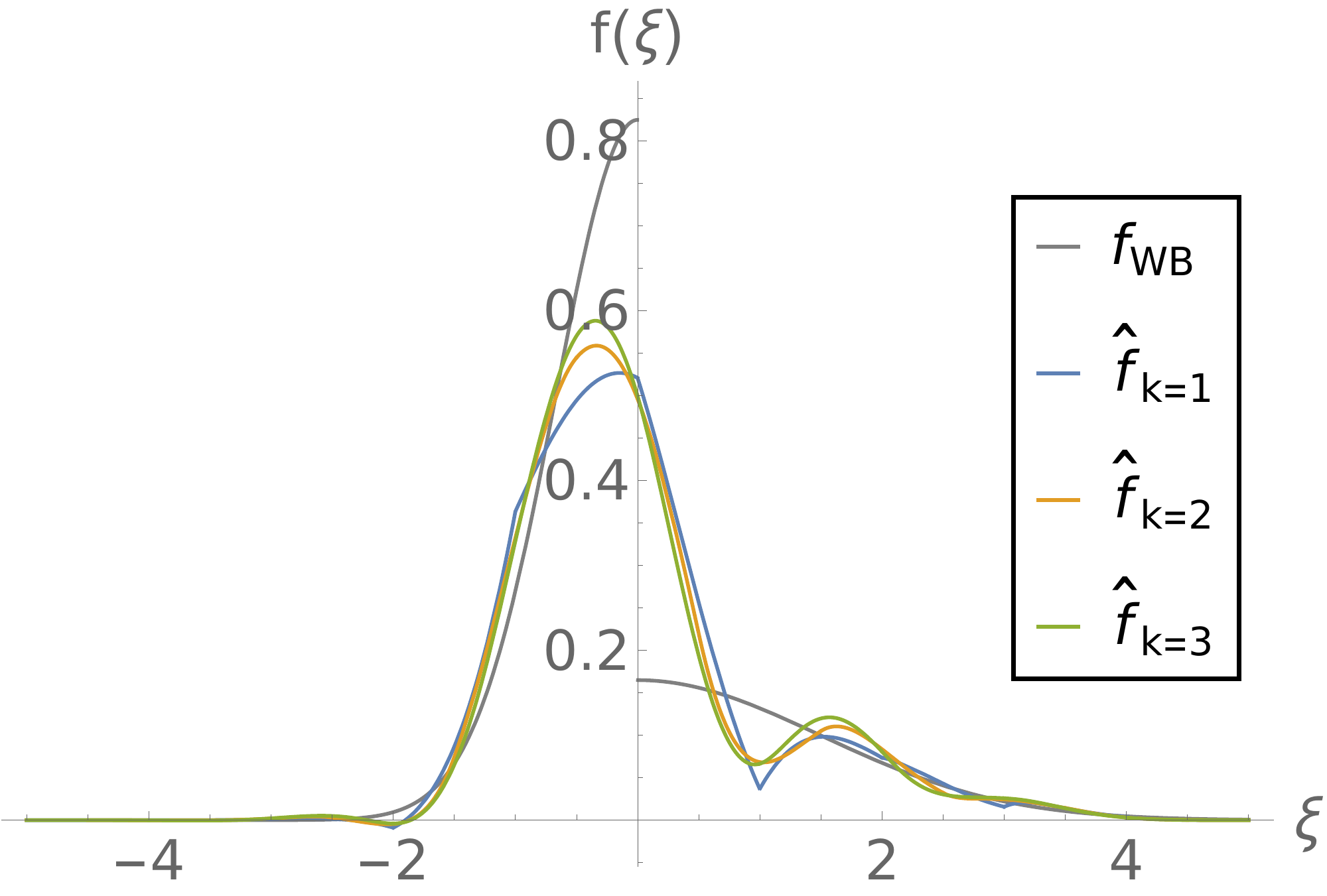}
		\caption{$n=9$.}
		\label{fig:approxWB1}
	\end{subfigure}
	\hfill
	\begin{subfigure}[b]{0.45\textwidth}
		\centering
		\includegraphics[width=\textwidth]{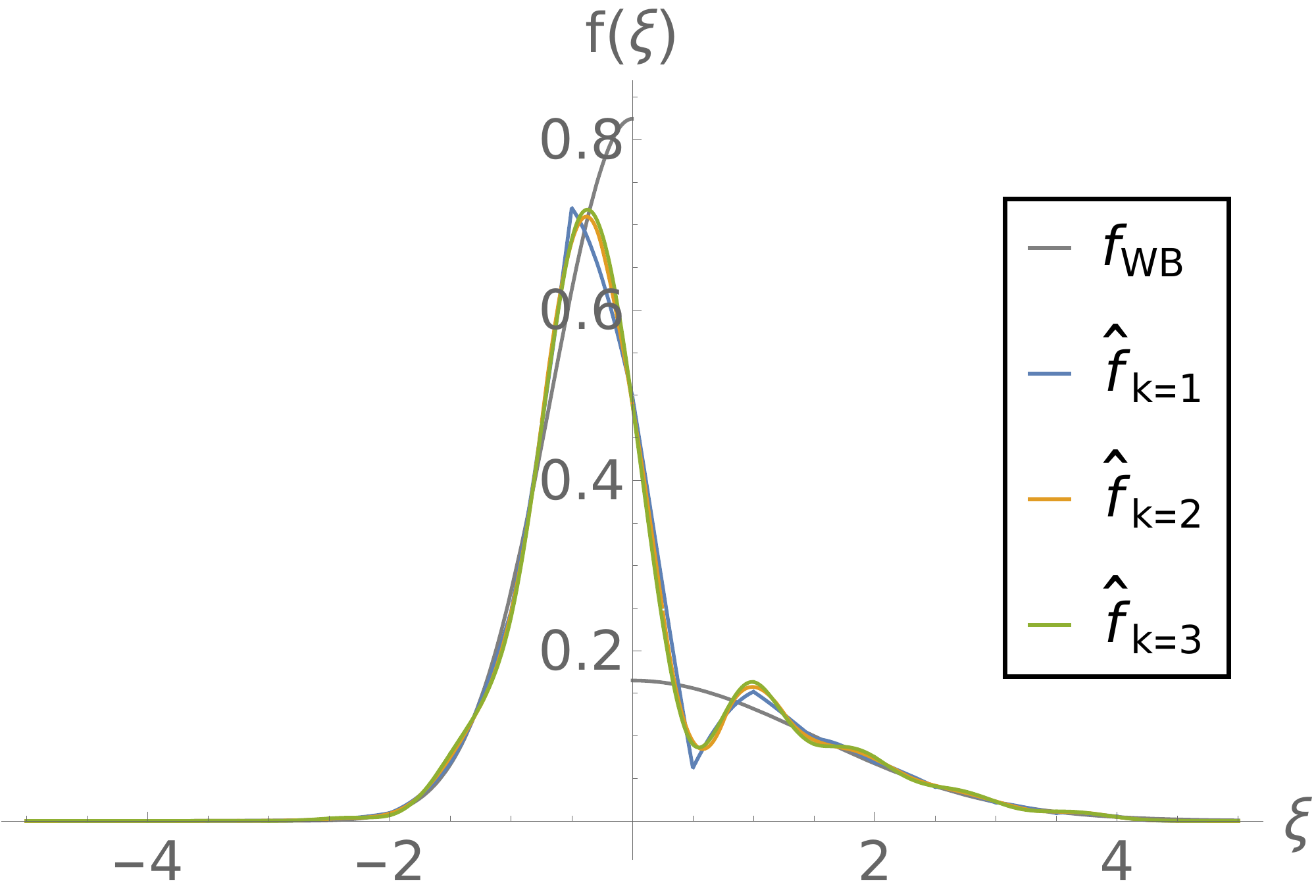}
		\caption{$n=17$.}
		\label{fig:approxWB2}
	\end{subfigure}
	\hfill
	\begin{subfigure}[b]{0.47\textwidth}
		\centering
		\includegraphics[width=\textwidth]{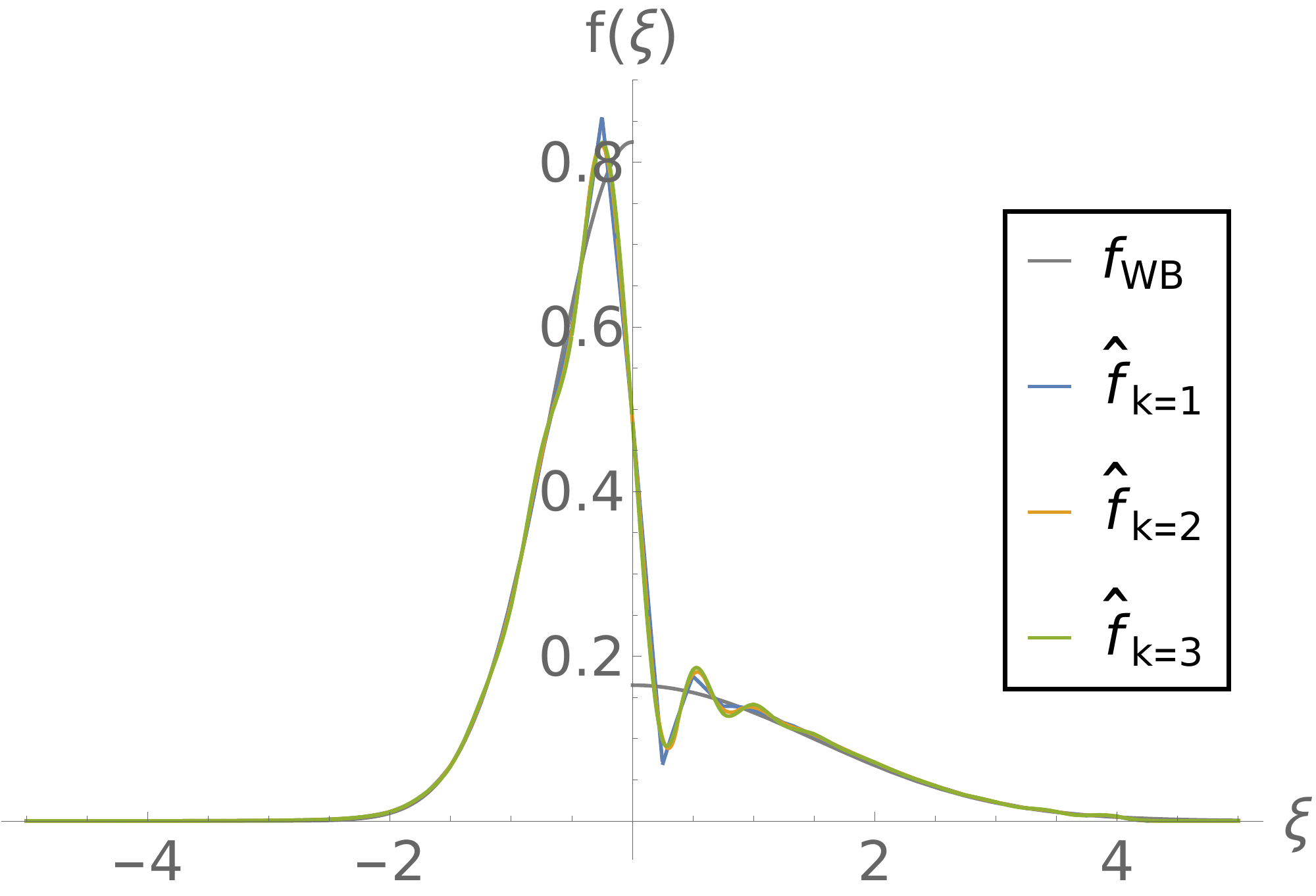}
		\caption{$n=33$}
		\label{fig:approxWB3}
	\end{subfigure}
	\hfill\begin{subfigure}[b]{0.47\textwidth}
		\centering
		\includegraphics[width=\textwidth]{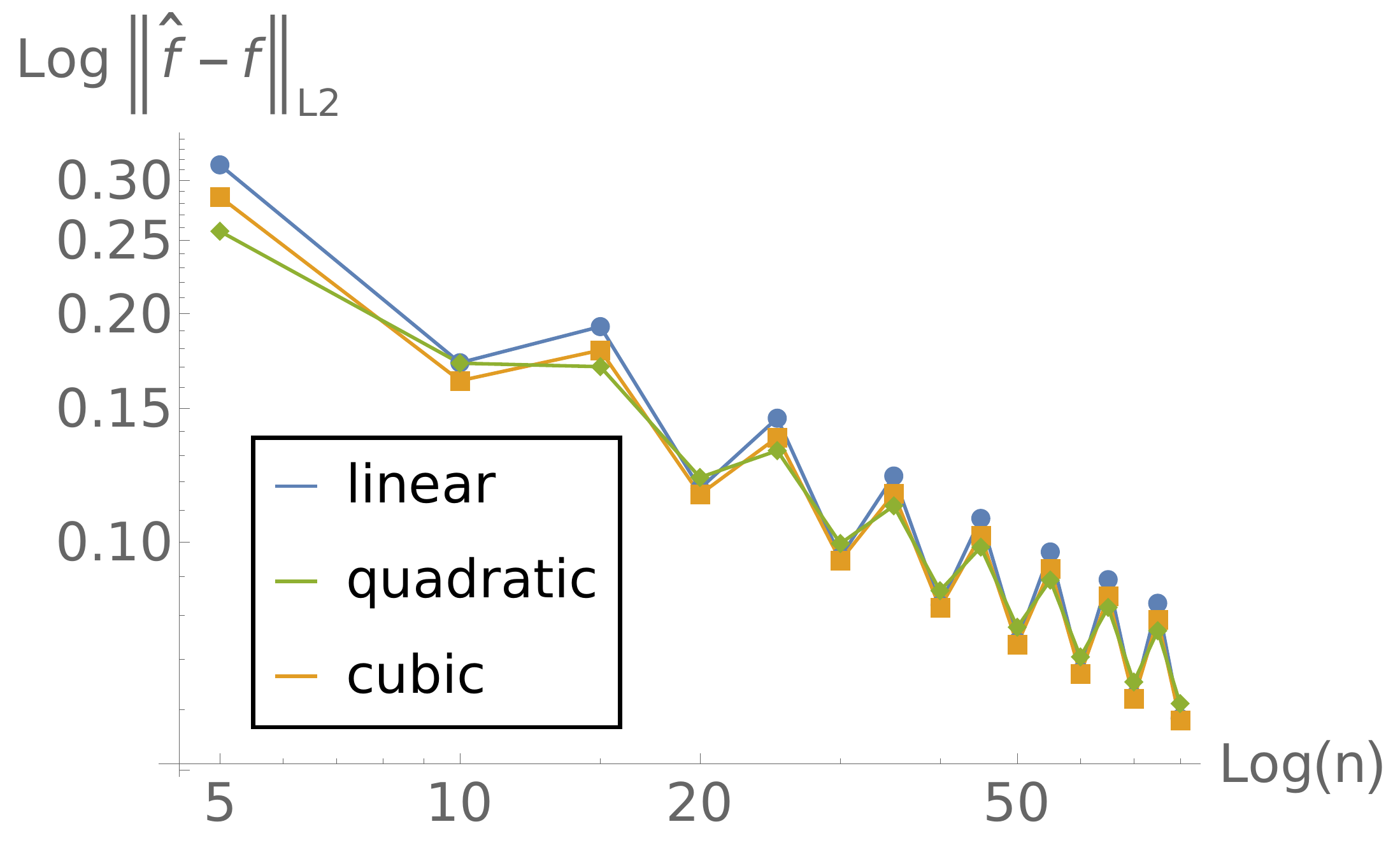}
		\caption{Logarithmic error for $n=5, \ldots, 80$.}
		\label{fig:errorWB}
	\end{subfigure}
	\caption{(a)-(c) approximation of discontinuous $f_{WB}$ for different $n=9,17,33$ and (d) approximation error depending on $n$. $[\xi_{min},\xi_{max}]=[-4,4]$.}
	\label{approxWB}
\end{figure}

We tested the approximation quality of a basis with fundamental constrained splines for this distribution function. As expected for this more challenging example, the error pictured in Figures \ref{fig:approxWB1} to \ref{fig:approxWB3} for linear ($k=1$) weighted FCS is slightly larger than for the previous distribution functions of this section. However, the function can still be approximated well, although there remain some oscillations around the discontinuity, which is a common and known issue when approximating discontinuous functions with continuous polynomials (Gibbs effect). The L2-error evolution for linear, quadratic and cubic FCS can be seen in Figure \ref{fig:errorWB}. It can be seen that the error is almost the same for the linear, quadratic, and cubic splines. Most importantly, an increase in $n$ leads to a more accurate approximation of the distribution function. The error does not decrease monotonically, but shows an even-odd decoupling of the number of splines. This is also a known behavior for the approximation of discontinuous functions by polynomials. In the context of moment methods, a similar effect has been observed in \cite{Torrilhon2015}. An additional effect for the spline basis is the fact that for even number of underlying B-splines $n$, the spline grid places grid points away from $\xi=0$, which is a better way for modelling the discontinuity than in the odd case, where one grid point is placed at $\xi=0$, compare Figure \ref{bereich}. However, in both even and odd cases, the error is clearly decreasing with $n$ indicating a convergence and sufficient accuracy of the spline approximation. We conclude that the weighted FCS ansatz is able to approximate even discontinuous distribution functions with sufficient accuracy.

\section{Spline Moment Equations for transformed equation}
\label{sec:5}
{%
We have described and analyzed how typical distribution functions can be approximated using a weighted FCS ansatz. In this section a set of PDEs for the evolution of the basis coefficients including the macroscopic values $\rho, v, \theta$ will be derived. The Boltzmann-BGK equation is transformed analogously to \cite{Koellermeier2014}. Then we use the aforementioned weighted FCS ansatz, perform a Galerkin projection and rewrite the resulting system in concise form.

Beginning with the one-dimensional Boltzmann-BGK equation
\begin{equation}
 \label{BoltzBKG}
 \dt f(t,x,c) + c\: \dx f(t,x,c) = S(f)  = \frac{1}{\tau}\left(\frac{\rho}{\sqrt{2 \pi \theta}} \exp\left(-\frac{(c-v)^2}{2\theta}\right)-f\right),
\end{equation}
we allow for an efficient discretization using spline functions by performing a non-linear transformation of the velocity space used in \cite{Kauf2011,Koellermeier2014}
\begin{equation}
\label{trans}
 \xi(x,t,c):=\frac{c-v(x,t)}{\sqrt{\theta(x,t)}}
\end{equation}
to obtain the following transformed equation denoting $f=f(t,x,\xi)$
\begin{multline}
\label{TransformedB}
 D_t f+\seta \xi \dx f + \partial_{\xi} f \left(-\frac{1}{\seta} (D_t v + \seta \xi \dx v) - \frac{1}{2\theta} \xi (D_t \theta + \seta \xi \dx \theta)\right) \\ =S(f,\rho,v,\theta) = \frac{1}{\tau}\left(\frac{\rho}{\sqrt{2 \pi \theta}} \exp(-\xi^2/2)-f\right),
\end{multline}
where $D_t:=\dt + v\:\dx$ denotes the convective time derivative.

Equation \eqref{TransformedB} indeed looks more complicated than the standard Boltzmann-BGK equation. However, we will see below that a concise system of equations can be derived from it in a straightforward way using the spline basis. The transformation \eqref{trans} makes it possible to use a very coarse grid in the velocity space. This grid typically only consists of a few spline basis functions. It is important to note that both $v(t,x)$ and $\theta(t,x)$ depend on the solution as moments of the distribution function and also depend on time and space. Due to the non-linearity of the transformation the ansatz in velocity space can be interpreted as using a moving grid in the velocity space, which is shifted by the velocity $v(t,x)$ and scaled by the variance $\sqrt{\theta(x,t)}$. This leads to a gain in efficiency and allows to achieve high accuracy with very few moments in comparison with standard schemes that simply expand around a global equilibrium function in the velocity space.

Following the literature, a scaled distribution function $\tilde f$ is used
\begin{equation}
\label{DefTildeF}
 \tilde f := \frac{\seta}{\rho}f,
\end{equation}
which leads to the following relations for moments of the scaled distribution function $\tilde f$
\begin{equation}
\label{one-o-one}
 \inftyint \tilde f(t,x,\xi)d\xi =1, \quad
 \inftyint \xi f(t,x,\xi)d\xi =0, \quad
 \inftyint \xi^2 f(t,x,\xi)d\xi =1.
\end{equation}
and Equation \ref{TransformedB} then results in the transformed Boltzmann-BGK equation

\begin{equation}
 \begin{split}
\label{FinalB}
\left(\frac{1}{\rho} D_t \rho - \frac{d}{2\theta} D_t \theta \right) \tilde f + D_t \tilde f \\+\: \seta \xi \left( \left( \frac{1}{\rho} \dx \rho - \frac{d}{2\theta}\dx \theta \right) \tilde f + \dx \tilde f \right)  \\+\: \partial_{\xi} \tilde f \left( - \frac{1}{\seta}(D_t v + \seta \xi \dx v ) - \frac{1}{2 \theta} \xi (D_t \theta + \seta \xi \dx \theta ) \right) \\=\: \frac{1}{\tau}\left(\frac{\rho}{\sqrt{2 \pi \theta}} \exp(-\xi^2/2)-\frac{\rho}{\seta}\tilde f\right).
 \end{split}
\end{equation}

According to the previous sections, we now choose a weighted FCS ansatz for $\tilde f$ so that Equations \ref{one-o-one} are fulfilled. For the expansion we need to fix a number $n \in \naturals$, $n \geq 4$, which will be the number of B-splines that underlie the FCS.
\begin{equation}
 \tilde f = \frac{1}{\sqrt{2\pi}} \exp(-\xi^2/2) \left( 1+ \sum_{i=1}^{n-3} \kappa_{i} \phi_i \right).
\end{equation}
The full vector of unknowns is given by $\left(\rho, v, \theta, \kappa_1, \dots, \kappa_{n-3}\right) \in \mathbb{R}^{n}$. Note that the number of unknowns -- and thus the number of equations in the upcoming equation system -- corresponds exactly to the number of underlying B-splines $n$ from the FCS construction.

For the Galerkin projection, the test functions include the first three monomials to reproduce the Euler equations, which are the conservation laws of mass, momentum, and energy
\begin{equation}
    \psi \in \{ 1, \xi, \xi^2 \} \cup (\hat \psi_j)_{j \in \{1,\ldots,n-3\}}
\end{equation}
and $(n-3)$ remaining FCS $\hat \psi_j$.

The resulting PDE system then reads
\begin{equation}
 \begin{split}
\label{angesetztB}
 \left( \frac{1}{\rho}D_t \rho - \frac{1}{2 \theta} D_t \theta \right) (M \kappa+V) + M D_t \kappa \\
 + \seta \left( \left( \frac{1}{\rho} \dx \rho - \frac{1}{2\theta} \dx \theta \right) (M^{\xi} \kappa + V^{\xi}) + M^{\xi} \dx \kappa \right) \\
 - \frac{1}{\seta}\left( D_t v (M^{\dxi} \kappa+V^{\dxi}) + \seta \dx v (M^{\xi \dxi} \kappa+V^{\xi \dxi}) \right)\\
 - \frac{1}{2 \theta}\left(D_t \theta (M^{\xi \dxi} \kappa+V^{\xi \dxi}) + \seta \dx \theta (M^{\xi \xi \dxi} \kappa+V^{\xi \xi \dxi}) \right)= \frac{1}{\tau} M \kappa
\end{split}
\end{equation}
{with} $n \times (n-3)$-{matrices} $M, M^{\xi}, M^{\dxi}, M^{\xi \dxi}$ and $M^{\xi \xi \dxi}$ as well as vectors $V$, $V^{\xi}$, $V^{\dxi}$, $V^{\xi \dxi}$ and $V^{\xi \xi \dxi}$ of length $n$.
{Denoting the scalar product} $\langle f,g \rangle:=\inftyint f(\xi)g(\xi) \; d\xi$ {and} $w=\frac{1}{\sqrt{2\pi}} \exp(-\xi^2/2)$ { the matrices and vectors can be written as}:
\begin{center}
$M_{ij}=\langle \psi_i,w \phi_j \rangle$, \quad
$V_{i}=\langle \psi_i,w \rangle$,\\
$M^{\xi}_{ij}=\langle \psi_i,\xi w\phi_j \rangle$, \quad
$V^{\xi}=\langle \psi_i,\xi w \rangle$,\\
$M^{\dxi}_{ij}=\langle \psi_i,\dxi(w\phi_j) \rangle$, \quad
$V^{\dxi}_{ij}=\langle \psi_i,\dxi w \rangle$,\\
$M^{\xi \dxi}_{ij}=\langle \psi_i,\xi \dxi (w\phi_j) \rangle$, \quad
$V^{\xi \dxi}_{ij}=\langle \psi_i,\xi \dxi w \rangle$, \\
$M^{\xi \xi \dxi}_{ij}=\langle \psi_i,\xi \xi \dxi (w\phi_j) \rangle$, \quad
{and} \quad
$V^{\xi \xi \dxi}_{ij}=\langle \psi_i,\xi \xi \dxi w \rangle$.
\end{center}
{%
In compact notation Equation \ref{angesetztB} can be written as
}
\begin{equation}
\label{ABsystem}
 B \: \dt  \begin{pmatrix}\rho\\v\\\theta\\\kappa_1\\\vdots\\\kappa_{n-3}\end{pmatrix} + A \: \dx \begin{pmatrix}\rho\\v\\\theta\\\kappa_1\\\vdots\\\kappa_{n-3}\end{pmatrix} = \frac{1}{\tau} M \begin{pmatrix}\kappa_1\\\vdots\\\kappa_{n-3}\end{pmatrix},
\end{equation}
{%
where matrices $A$ and $B$ are defined as
}
\begin{equation}
 A=
 \begin{pmatrix}
    A_1 \: \vline \: A_2 \: \vline \: A_3 \: \vline \: v M + \seta M^{\xi}
   \end{pmatrix}
   \in \reals^{n\times n}
\end{equation}
{with}
\begin {gather}
 A_1=\frac{v}{\rho}\left(M \kappa +V\right)+\frac{\seta}{\rho} \left(M^{\xi} \kappa+ V^{\xi}\right), \\
 A_2=-\frac{v}{\seta}\left(M^{\dxi} \kappa + V^{\dxi}\right) - \left(M^{\xi \dxi} \kappa +V^{\xi \dxi}\right), \\
 A_3=-\frac{1}{2 \theta}\left(v \left(\left(M+M^{\xi \dxi}\right) \kappa +V+V^{\xi \dxi}\right)+ \seta \left(\left(M^{\xi}+M^{\xi \xi \dxi}\right) \kappa +V^{\xi}+V^{\xi \xi \dxi}\right)\right)
\end {gather}
{and}
\begin{equation}
 B=
 \begin{pmatrix}
  B_1 \: \vline \: B_2 \: \vline \: B_3 \: \vline \: M
 \end{pmatrix}
 \in \reals^{n \times n}
\end{equation}
{with}
\begin{gather}
 B_1=\frac{1}{\rho}\left(M \kappa + V\right), \\B_2=-\frac{1}{\seta}\left(M^{\dxi} \kappa +V^{\dxi}\right),\\
 B_3=-\frac{1}{2\theta}\left(\left(M+M^{\xi \dxi}\right) \kappa +V+V^{\xi \dxi}\right).
\end{gather}

{Because of} $B^{-1} \left(B_1 \: \vline \: B_2 \: \vline \: B_3 \: \vline \: M\right)=B^{-1}B=I_n$, {we know that}
$B^{-1} M =
\begin{pmatrix}
0 \\ I_{n-3}
\end{pmatrix}.
$

After multiplication with $B^{-1}$ Equation \ref{ABsystem} yields the final system of PDEs
\begin{equation}
    \label{AsysDGL}
    \dt  \begin{pmatrix}\rho\\v\\\theta\\\kappa_1\\\vdots\\\kappa_{n-3}\end{pmatrix} + B^{-1}A \: \dx \begin{pmatrix}\rho\\v\\\theta\\\kappa_1\\\vdots\\\kappa_{n-3}\end{pmatrix} = \frac{1}{\tau} \begin{pmatrix}0\\0\\0\\\kappa_1\\\vdots\\\kappa_{n-3}\end{pmatrix}.
\end{equation}

The specific entries of the system matrix $\asys(\rho,v,\kappa)=B^{-1} A$ can be easily computed after choosing the number of unknowns $n$ and a corresponding grid for the FCS. The system \eqref{AsysDGL} is then called Spline Moment Equations (SME).

\begin{remark}
    Due to choosing the first three test functions as $1, \xi, \xi^2$, the first three rows of Equation \ref{AsysDGL} are precisely the Euler equations
    \begin{equation}
        D_t \rho + \rho \dx v = 0, \quad D_t v + \frac{1}{\rho} \dx p =0, \quad D_t \theta + \frac{1}{\rho} \dx q + \frac{2p}{\rho} \dx v = 0
    \end{equation}
    with pressure tensor $p$ and heat flux $q$
    \begin{equation}\label{e:Q}
        p = \rho \theta \inftyint \xi^2 \tilde f d\xi, \quad q = \rho \theta^{(3/2)} \inftyint \xi^3 \tilde f d\xi,
    \end{equation}
    which can be verified by inserting the test functions ${1, \xi, \xi^2}$ into Equation \ref{angesetztB}.
\end{remark}

\begin{remark}
    Due to the transformation \ref{trans} the velocity space discretization becomes non-uniform as it depends on $v(x,t)$ and $\theta(x,t)$. This means that the spline functions are adaptively positioned around the mean velocity $v(x,t)$ and the range of their support scales with the temperature $\theta(x,t)$. The range of the support is then effectively $[v-\xi_{min}\sqrt{\theta},v+\xi_{max}\sqrt{\theta}]$. Different regions of the flow field with different $v(x,t)$ and $\theta(x,t)$ can thus be accurately approximated using one single basis with relatively few unknowns. A similar procedure for discrete velocity methods was used with $\xi_{min/max} = \pm 4$ in \cite{Brull2014}. As the same range led to good results for the approximation properties, in Section \ref{sec:4} we will make use of the same range for most of the next section.
\end{remark}

\subsection{Eigenvalues and hyperbolicity}
We briefly consider the properties of the system matrix. Only if all eigenvalues of the system matrix are real-valued, the PDE system is globally hyperbolic and constitutes a meaningful physical model for the simulation of rarefied gases. Previously developed models using Hermite polynomials are not globally hyperbolic, see \cite{Koellermeier2014}, but there are some known procedures that yield global hyperbolicity, see \cite{Fan2016,Koellermeier2017b}.

We compute the eigenvalues for different numbers $(n-3)$ of fundamental constrained splines, for order $k=1$ and $\xi$-range $[\xi_{min}, \xi_{max}]=[-4,4]$.

Similarly as in \cite{Cai2013} and \cite{Koellermeier2014} all eigenvalues $\zeta_i$ of $\asys$ are shifted by the mean velocity $v$ and scaled with the square root of the temperature $\theta$, i.e.
\begin{equation}
\label{zeta}
 \zeta_i = v \pm \seta \beta_i \kappas.
\end{equation}

The eigenvalue $\zeta_i$ is thus real-valued if and only if $\beta_i \kappas$ is real. $\beta_i \kappas$ is the zero of the modified characteristic polynomial
\begin{equation}
 P(\lambda) = \det (\asys - (\lambda \seta +v)\:E) = 0.
\end{equation}

In Figure \ref{evhigh} we show the hyperbolicity area of the SME using different number of equations $n=5,7,9,11$, depending on the parameters $\kappa_i$ that correspond to the \emph{center} basis functions with an index close to $\frac{(n-3)}{2}$ while setting the outer kappas to zero. Hyperbolicity is only obtained in a small domain around the origin, which represents equilibrium. The hyperbolicity domain shrinks with increasing $n$. Furthermore, the hyperbolicity plots show striking similarity with the hyperbolicity plots for the Hermite-based ansatz in \cite{Cai2013}.

\begin{figure}[H]
\centering
 \begin{subfigure}[b]{0.47\textwidth}
 \centering
 \includegraphics[width=\textwidth]{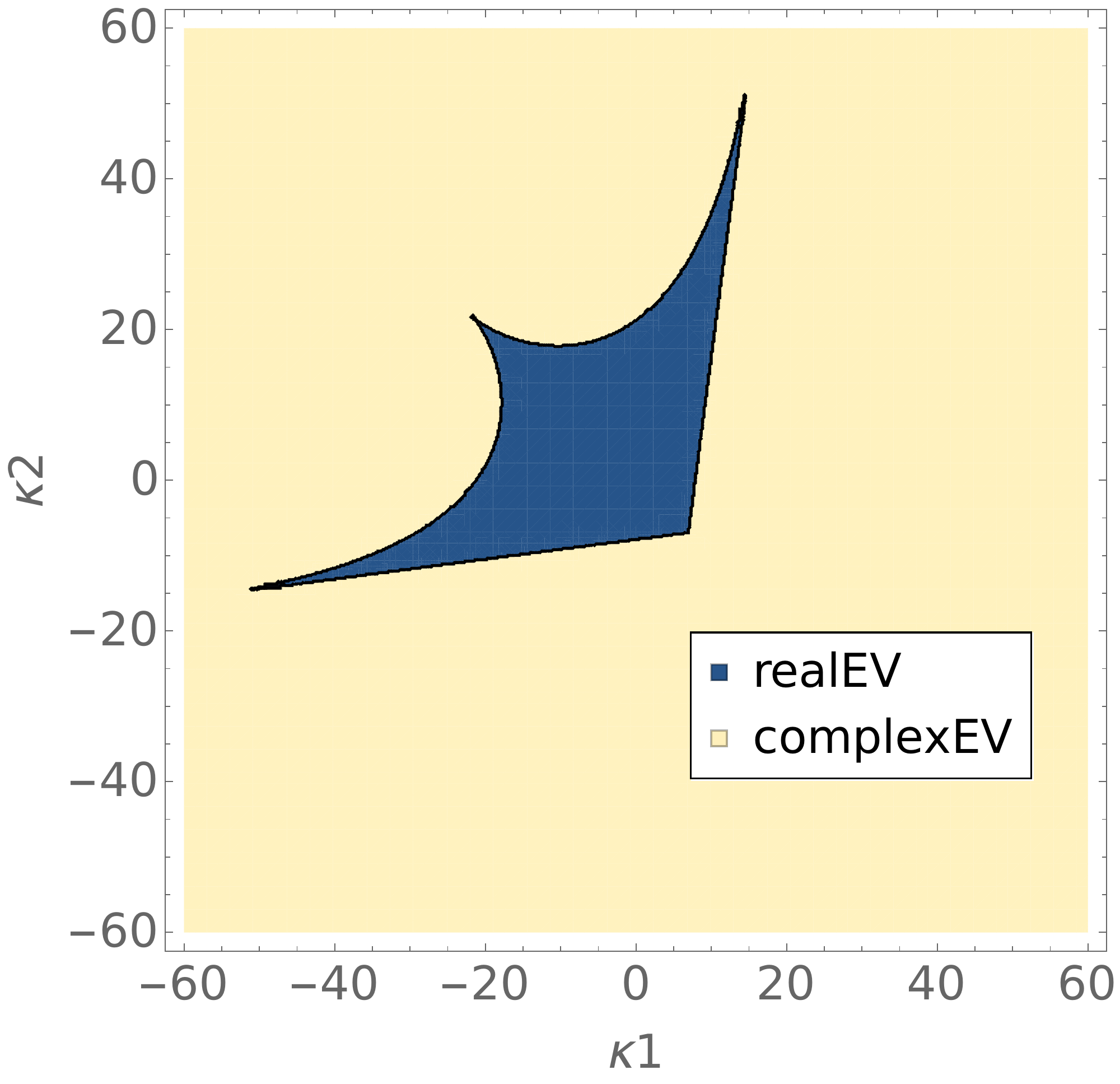}
 \caption{$n=5$.}
 \label{evunf1}
\end{subfigure}
\hfill
\begin{subfigure}[b]{0.47\textwidth}
 \centering
 \includegraphics[width=\textwidth]{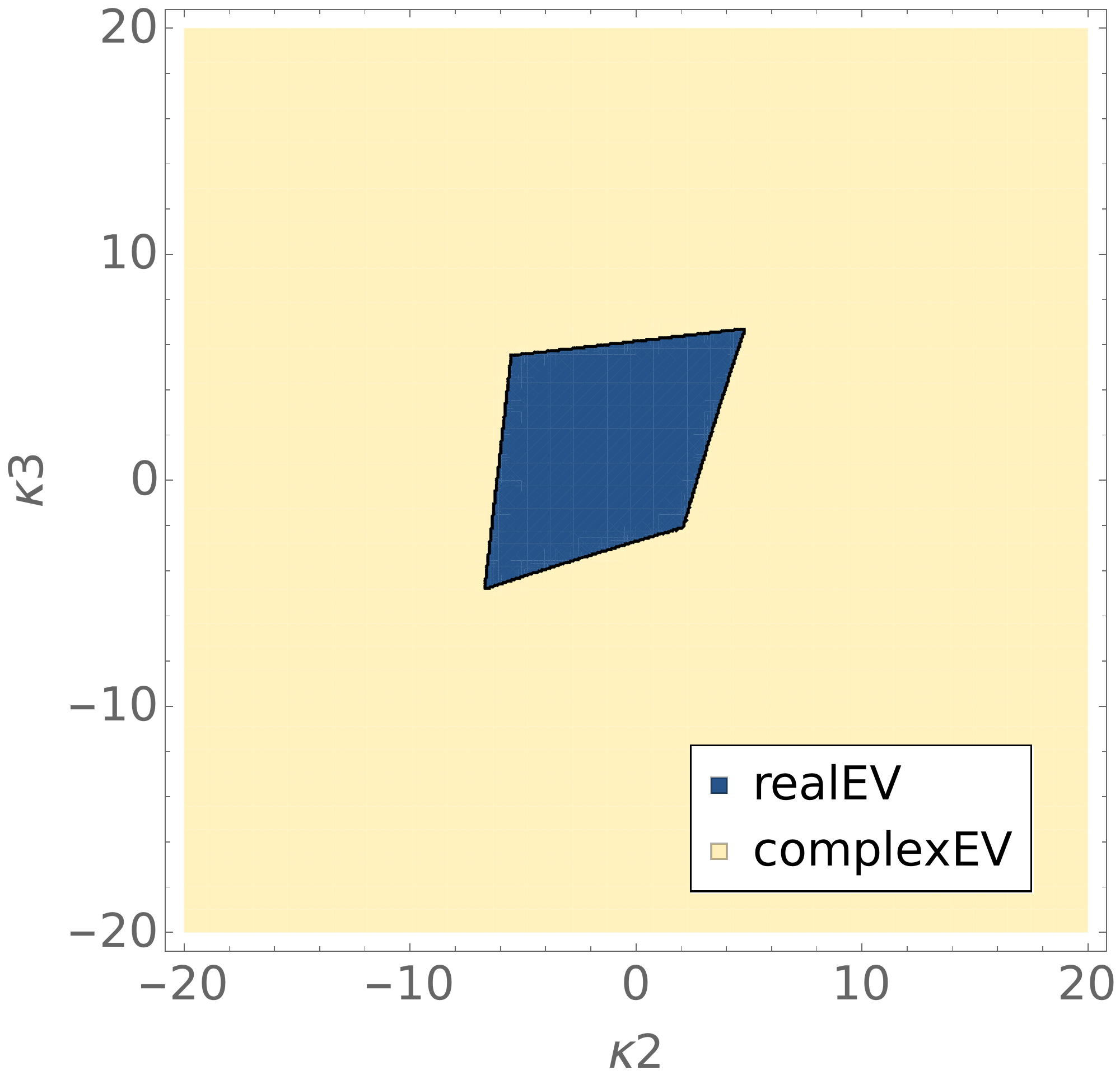}
 \caption{$n=7$.}
 \label{evunf2}
\end{subfigure}
\hfill
\begin{subfigure}[b]{0.47\textwidth}
 \centering
 \includegraphics[width=\textwidth]{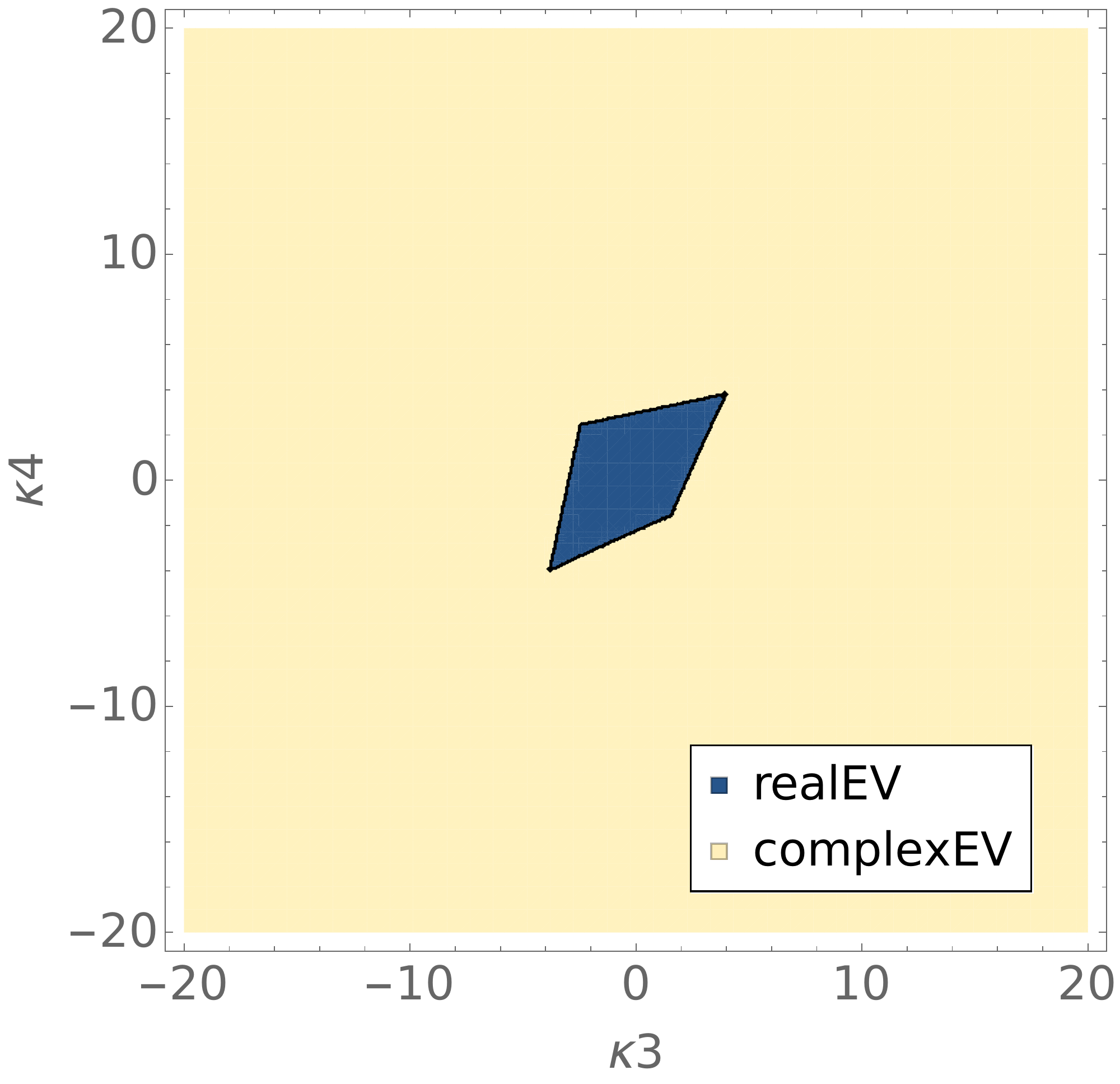}
 \caption{$n=9$.}
 \label{evunf3}
 \end{subfigure}
 \begin{subfigure}[b]{0.47\textwidth}
 \centering
 \includegraphics[width=\textwidth]{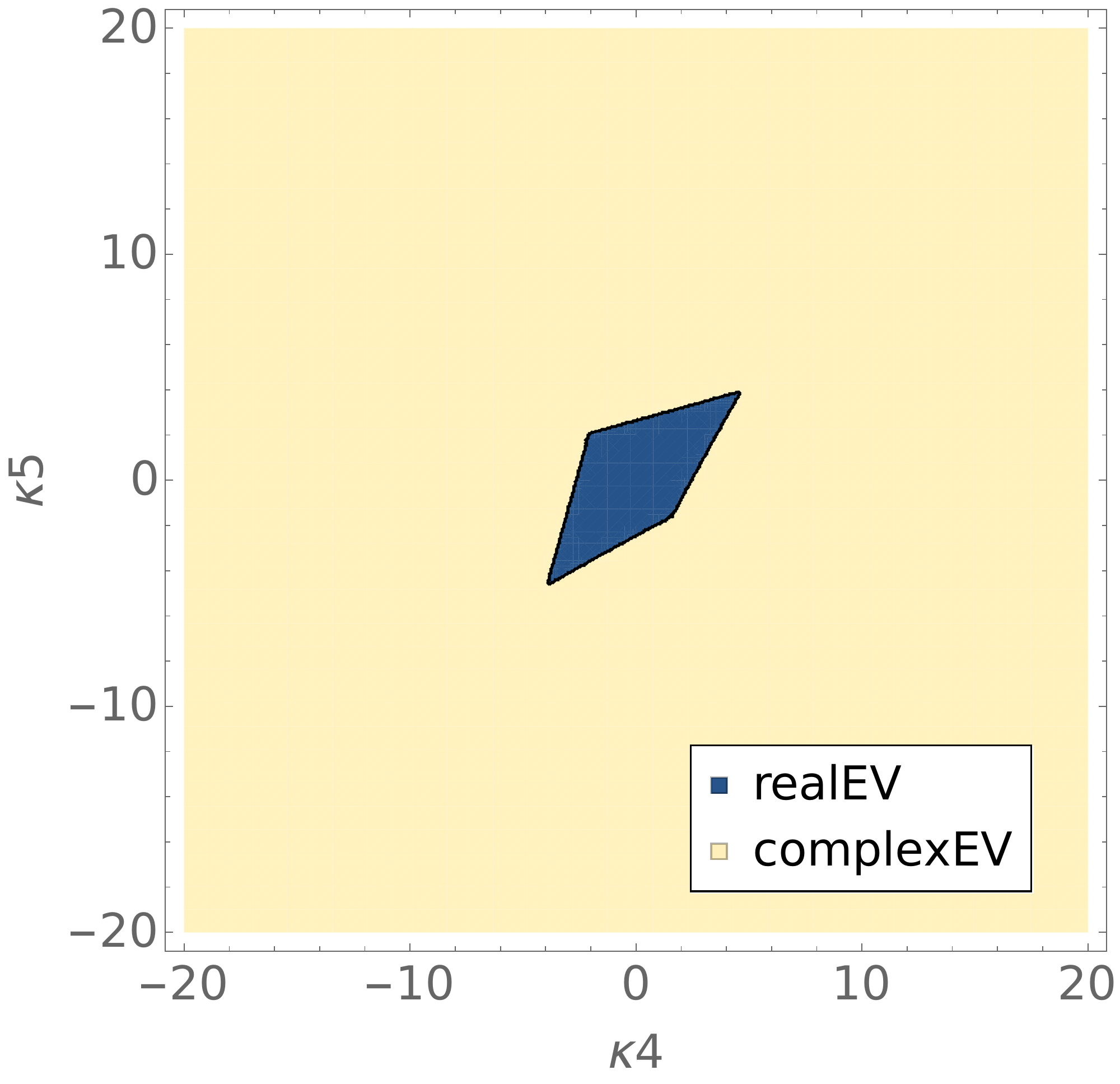}
 \caption{$n=11$.}
 \label{evunf4}
\end{subfigure}
\caption{{(a)-(d) Hyperbolicity domain for varying $n=5,7,9,11$ depending on center $\kappa_{i/i+1}$.}}
\label{evhigh}
\end{figure}
{%
Similar results are obtained for different $n \in \{1,\ldots, 12\}$, support ranges $[\xi_{min}, \xi_{max}]$, and orders $k$, for which we have always observed bounded hyperbolicity domains due to the occurrence of imaginary eigenvalues for larger values of the coefficients.

The SME thus lack global hyperbolicity, similar to Grad's system \cite{Grad1949}.

\begin{remark}
    Despite the lack of hyperbolicity, the model is expected to give satisfactory results within the hyperbolic domain, just like Grad's model. We note that this paper focusses on the spline approximation itself and further research is necessary for the hyperbolicity of the model. For details on hyperbolic regularisation of different models, we refer to \cite{Fan2016,Koellermeier2014b} and the recently developed diagram notation for the derivation of hyperbolic model in \cite{Koellermeier2020b}. As we do not focus on a hyperbolic regularization in this paper, we only consider a very basic possibility to achieve a hyperbolic PDE system in this work. It is based on a linearization around the equilibrium state, which is similar to the regularization in \cite{Cai2013} in a different variable setting and \cite{Koellermeier2020c} for another application in shallow flows. In Appendix \ref{app} we show the application of this linearized SME model yielding global hyperbolicity. Unfortunately, the simplifications of the linearization are too severe and the model does no longer converge to the solution of the Boltzmann-BGK equation. It is thus necessary to try one of the more advanced methods mentioned in the references above.
\end{remark}
\section{Simulation results}
\label{sec:6}

Three different test cases are used for the simulation of the new SME model. We will, however, focus on standard benchmark initial value problems and stationary problems, without the influence of boundary values. The definition of boundary values for moment models is a separate topic. An outline on how to derive boundary conditions can already be found in \cite{Grad1949} and applications can be found e.g. in \cite{Cai2018}. The derivation of boundary conditions for the SME model can in principle be done in the same way but is left for future work. In order to assess the model accuracy and not have the accuracy spoilt by numerical errors resulting from the spatial discretization, we choose a rather fine spatial discretization similar to the tests in \cite{Cai2013, Koellermeier2017b, Koellermeier2017, schaerer2015}. However, it is clear that the model can also be used for simulations with coarser spatial grids, which will only affect the spatial discretization error and not the model accuracy itself.

\subsection{Shock tube test case}
We employ a standard shock tube test case for the first simulations of the new SME model \eqref{AsysDGL}. This test case is widely used with the same settings as in \cite{Cai2013,Koellermeier2017b}.
\begin{figure}[H]
 \centering
 \includegraphics[scale=.3]{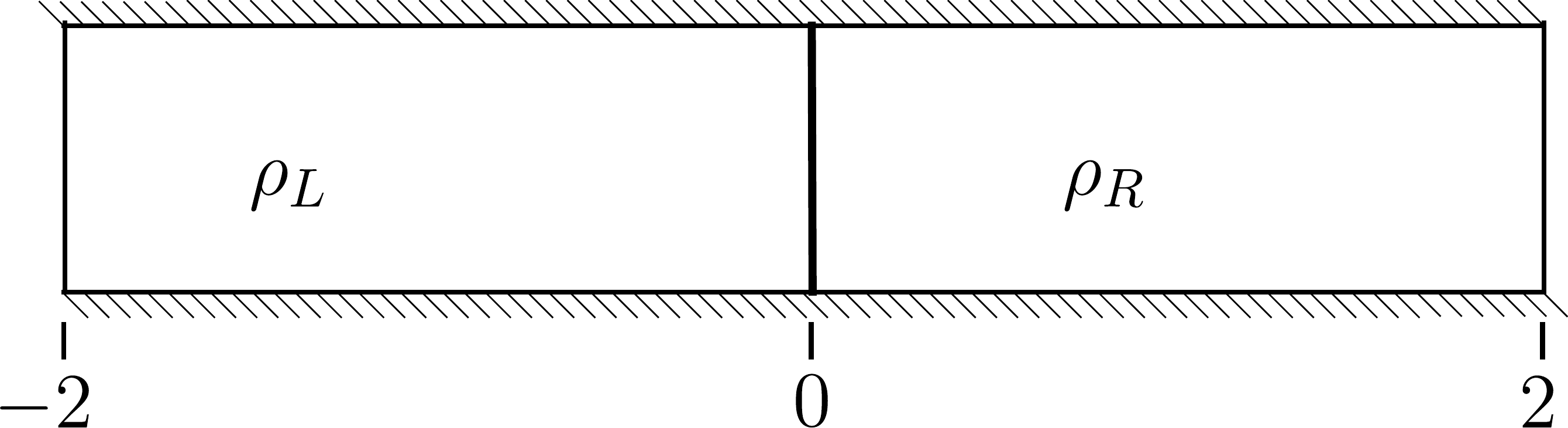}
 \caption{Set-up shock tube test.}
 \label{shockTubeSketch}
\end{figure}

The initial density ratio is given by $\rho_L=7$ and $\rho_R=1$, while the gas is at a rest with uniform temperature of $\theta=1$. The pressure is computed using $p = \rho \theta$. The BGK collision operator uses a non-linear relaxation time $\tau = \textrm{Kn}/ \rho$ to model collisions.

The numerical simulation was conducted with the explicit first order path-consistent non-conservative finite volume solver from \cite{Koellermeier2017b} using $4000$ points for time step size $\Delta t = 0.0001$, which corresponds to a CFL number of approximately $0.5$. The results are shown at $t_\textrm{END}=0.3$. We consider a small Knudsen number $\textrm{Kn}=0.05$ and a larger Knudsen number $\textrm{Kn}=0.5$ for more rarefied conditions. We will start from a base test setting of $n=4$ and $k=1$, $[\xi_{min},\xi_{max}]=[-2,2]$ and subsequently vary the different parameters. As a reference for the exact solution we use results from a high-resolution discrete velocity method also used in \cite{Koellermeier2017b}.

\subsubsection{Effect of the spline order $k$}
In Figure \ref{diforders} a test case with varying spline order is shown and only minimal differences can be observed. Only in the enlarged view slightly different subshock positions can be observed, whereas the main features remain unchanged.
\begin{figure}[ht]
 \centering
 \begin{subfigure}[t]{0.49\textwidth}
 \centering
 \includegraphics[width=\textwidth]{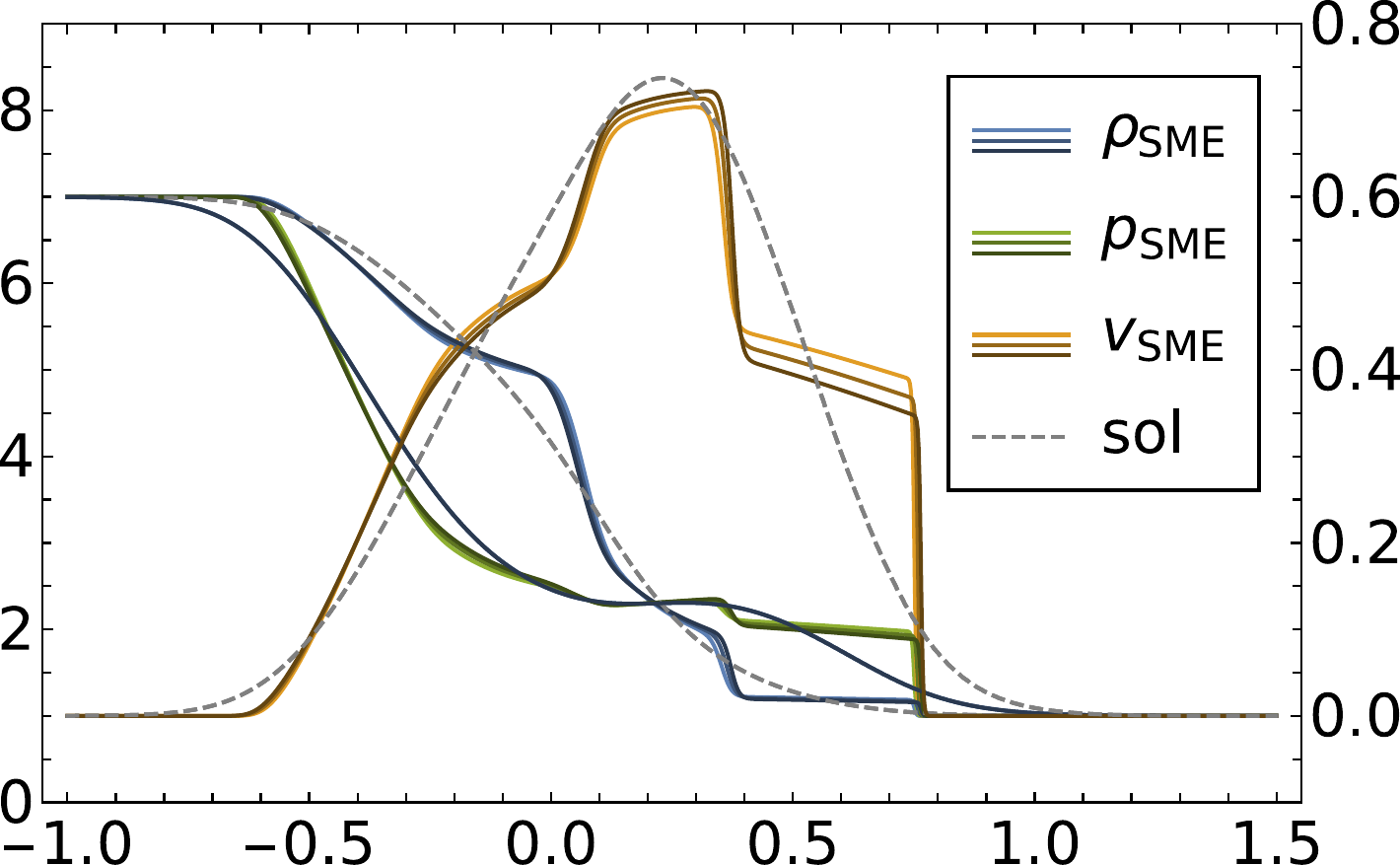}
 \caption{Full view.}
 \label{priforders2}
\end{subfigure}
\hfill
\begin{subfigure}[t]{0.5\textwidth}
 \centering
 \includegraphics[width=\textwidth]{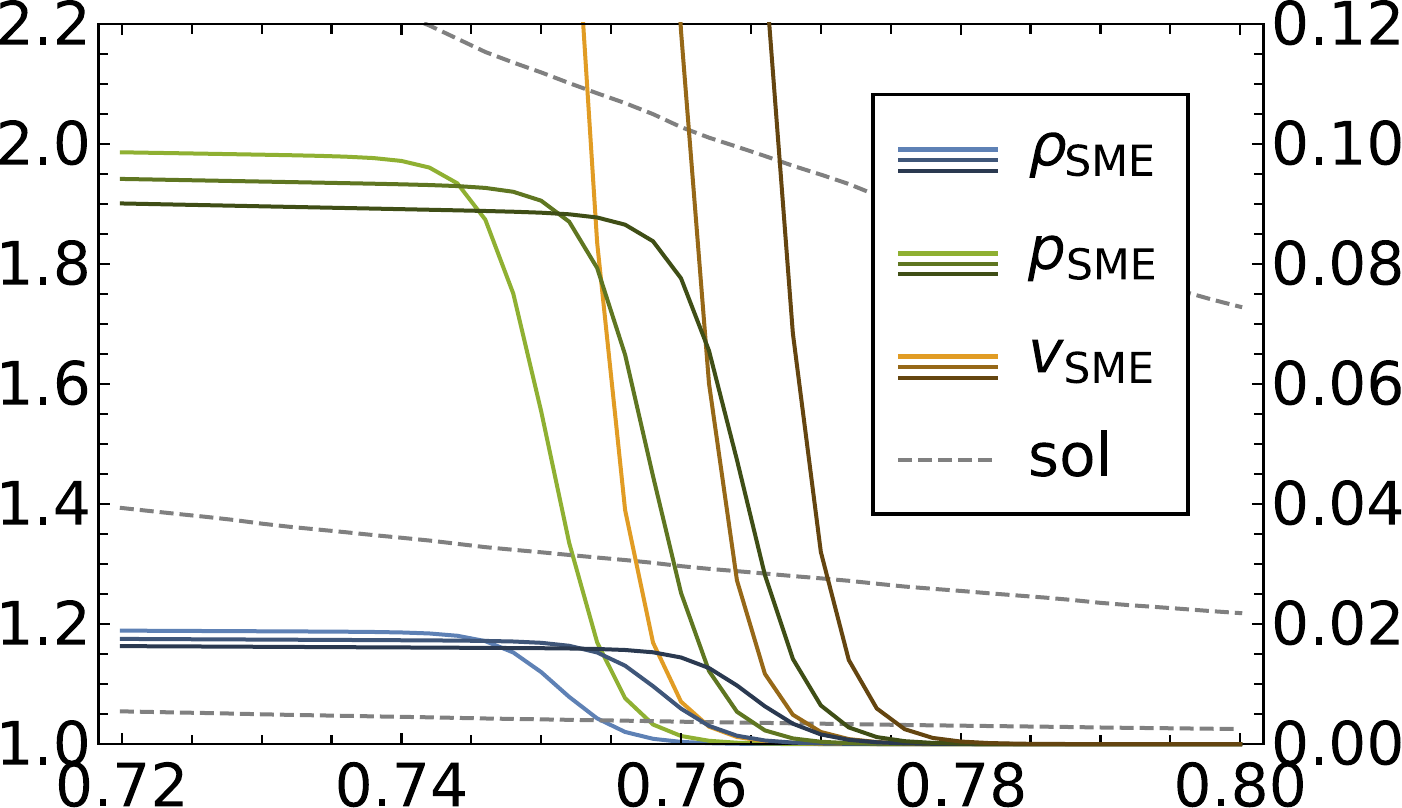}
 \caption{Enlarged view of (a).}
 \label{diforderszoomed}
 \end{subfigure}
\caption{Shock tube result for different spline order $k$ from bright to dark, $[\xi_{min},\xi_{max}]=[-2,2]$, $Kn=0.5$, $n=4$. Left axis for $\rho$ and $p$, right axis for $v$. Reference solution ``sol'' from DVM. (a) full solution, (b) zoom-in.}
\label{diforders}
\end{figure}
Most notably, we do not achieve better results for splines of higher order. A plausible explanation is the same effect that we have already seen when discussing the approximation properties in Figure \ref{unweightedb/plot}. When the number of splines is small in relation to the $\xi$-range, spline functions can only very roughly scan the original function's shape. Smaller details like rounded edges, that are the strength of higher-order splines, are not yet relevant here.

\begin{remark}
\label{remark_hyp}
    Figure \ref{diforders} clearly shows a subshock pattern, which is characteristic for moment models. The same can be observed for the standard Grad, HME or QBME models described in \cite{Koellermeier2017b}. The reason is the underlying set of bounded propagation speeds of the model. The subshocks typically become smaller when increasing $n$ and eventually vanish, as we will see later.
\end{remark}

\subsubsection{{Effect of the support $[\xi_{min},\xi_{max}]$}}
To compare different $\xi$-ranges in a reasonable way, the number of splines should remain unchanged when extending the $\xi$-range. The spline grid parameter $\Delta \xi=\frac{\xi_{max}-\xi_{min}}{(n-1)}$ is therefore kept constant.

We show simulations on the three $\xi$-ranges $[\xi_{min},\xi_{max}]=[-2, 2]$, $[-3, 3]$ and $[-4, 4]$ for all three orders. For $[-5,5]$ and $[-6,6]$ the simulation in the case $\textrm{Kn}=0.5$, $k=1$ was unstable for some particular numbers of splines. This may be caused by the occurrence of imaginary eigenvalues. However, the simulation was successful for the spline orders $k=2$ and $k=3$ with otherwise the same parameter setting.


Figure \ref{resultsDifferentBases} shows the results for a fixed spline distance of $\Delta \xi=\frac{2}{3}$.
It can be seen that the results are more accurate for a wider $\xi$-range. This effect is much stronger at the transition from $[-2,2]$ to $[-3,3]$ than at the transition from $[-3,3]$ to $[-4,4]$, probably due to the fact that the distribution function is decreasing for larger values of $\xi$. It is thus less problematic to assume $f(\xi)=0$ far away from the origin. For the largest range, the solution has already converged in the case $\textrm{Kn}=0.05$ and it is very accurate for $\textrm{Kn}=0.5$.
\begin{figure}[H]
 \centering
 \begin{subfigure}[t]{0.49\textwidth}
 \centering
 \includegraphics[width=\textwidth]{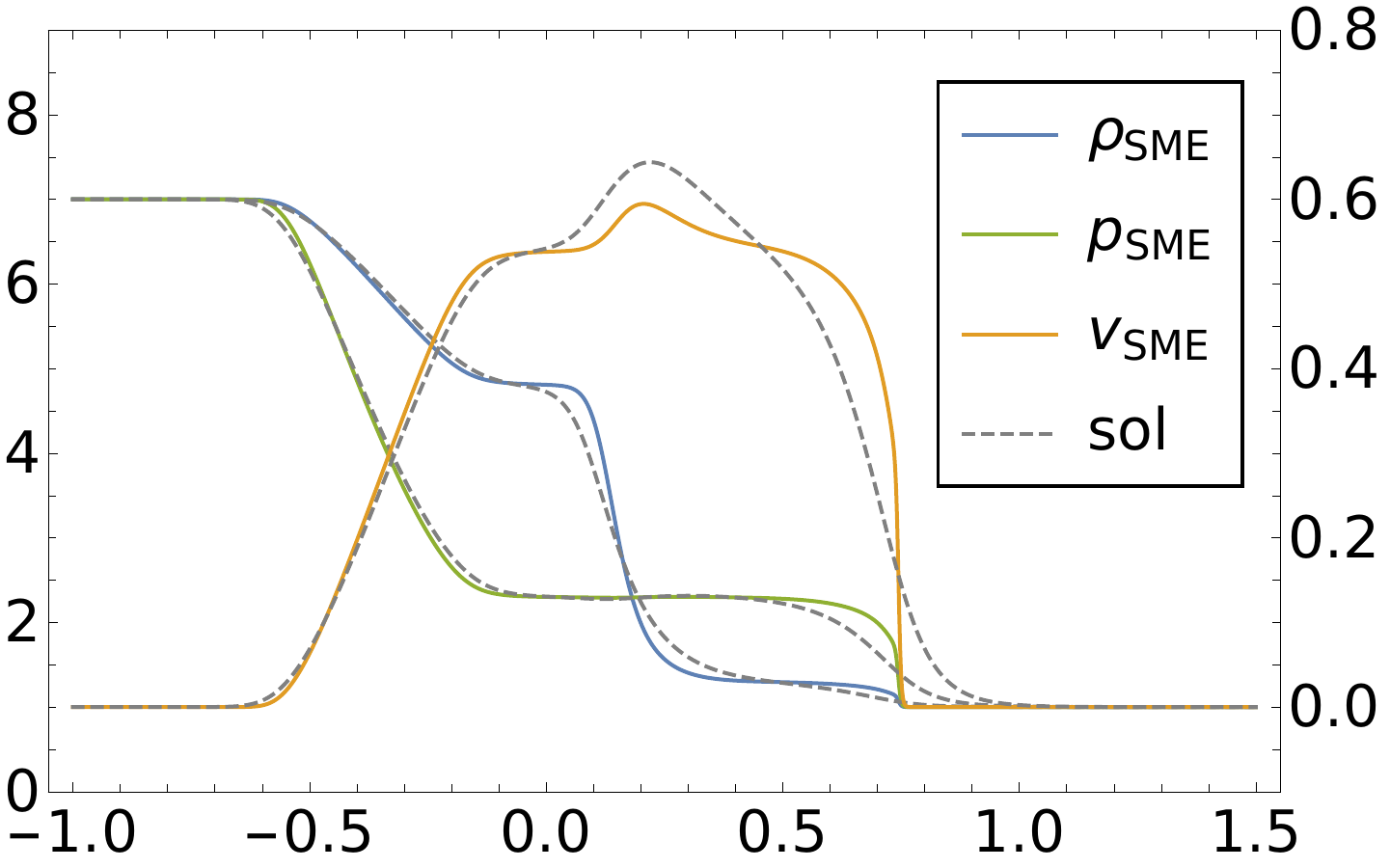}
 \caption{n=7 {on} $[-2,2]$, $\textrm{Kn}=0.05$.}
 \label{base1}
\end{subfigure}
\hfill
\begin{subfigure}[t]{0.49\textwidth}
 \centering
 \includegraphics[width=\textwidth]{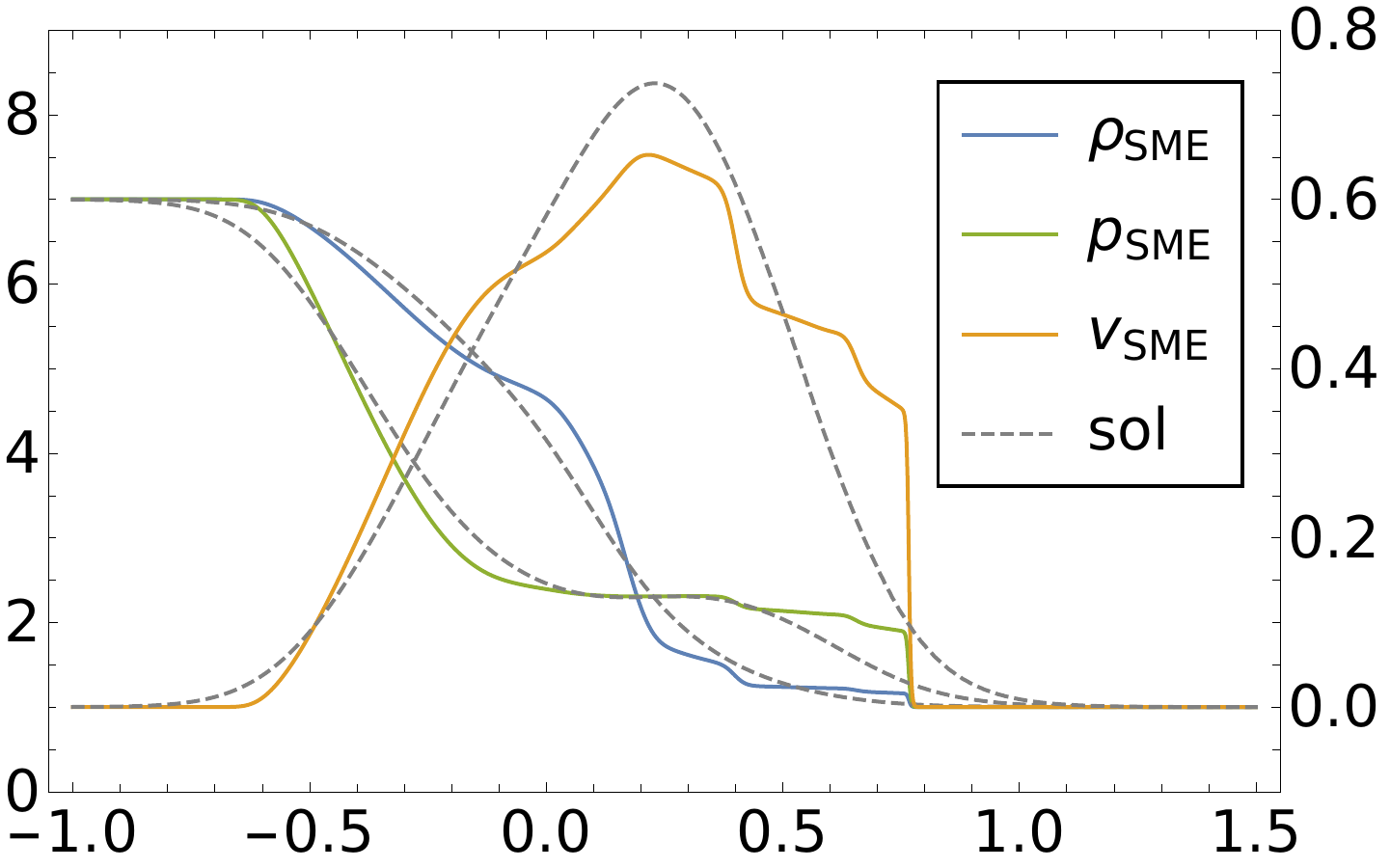}
 \caption{n=7 {on} $[-2,2]$, $\textrm{Kn}=0.5$.}
 \label{base2}
\end{subfigure}
\hfill
 \begin{subfigure}[t]{0.49\textwidth}
 \centering
 \includegraphics[width=\textwidth]{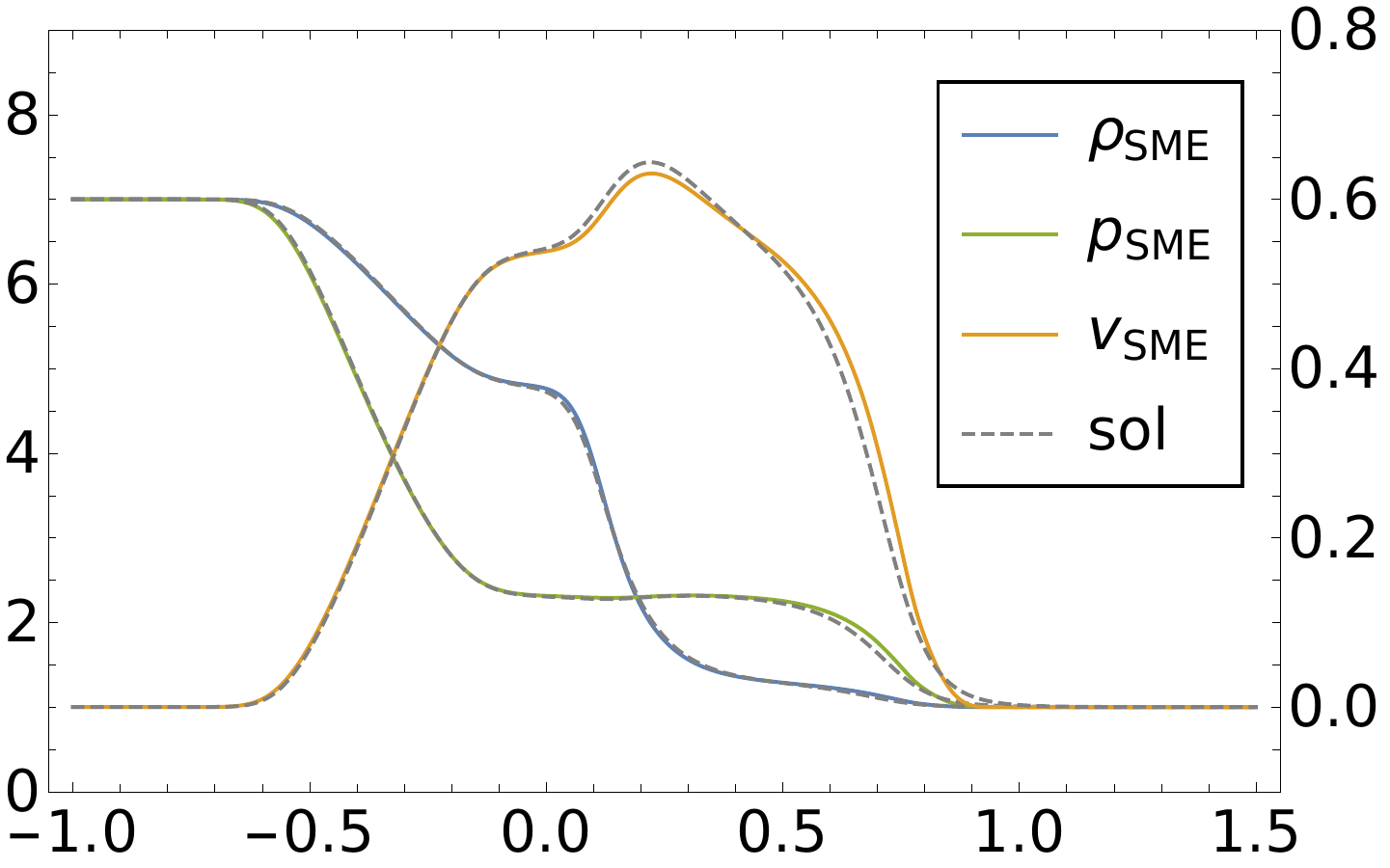}
 \caption{n=10 {on} $[-3,3]$, $\textrm{Kn}=0.05$.}
 \label{base3}
\end{subfigure}
\hfill
\begin{subfigure}[t]{0.49\textwidth}
 \centering
 \includegraphics[width=\textwidth]{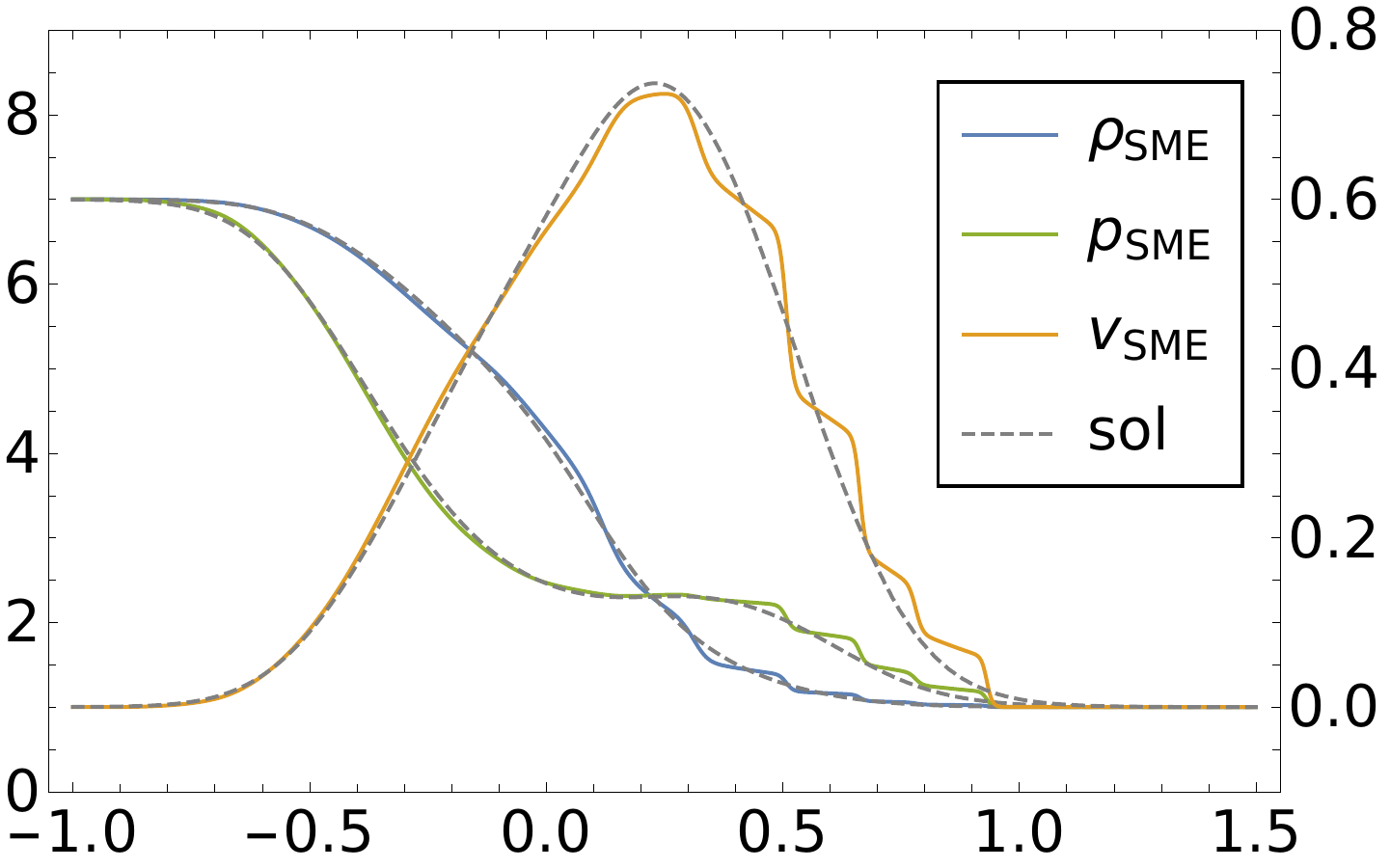}
 \caption{n=10 {on} $[-3,3]$, $\textrm{Kn}=0.5$.}
 \label{base4}
\end{subfigure}
\hfill
 \begin{subfigure}[t]{0.49\textwidth}
 \centering
 \includegraphics[width=\textwidth]{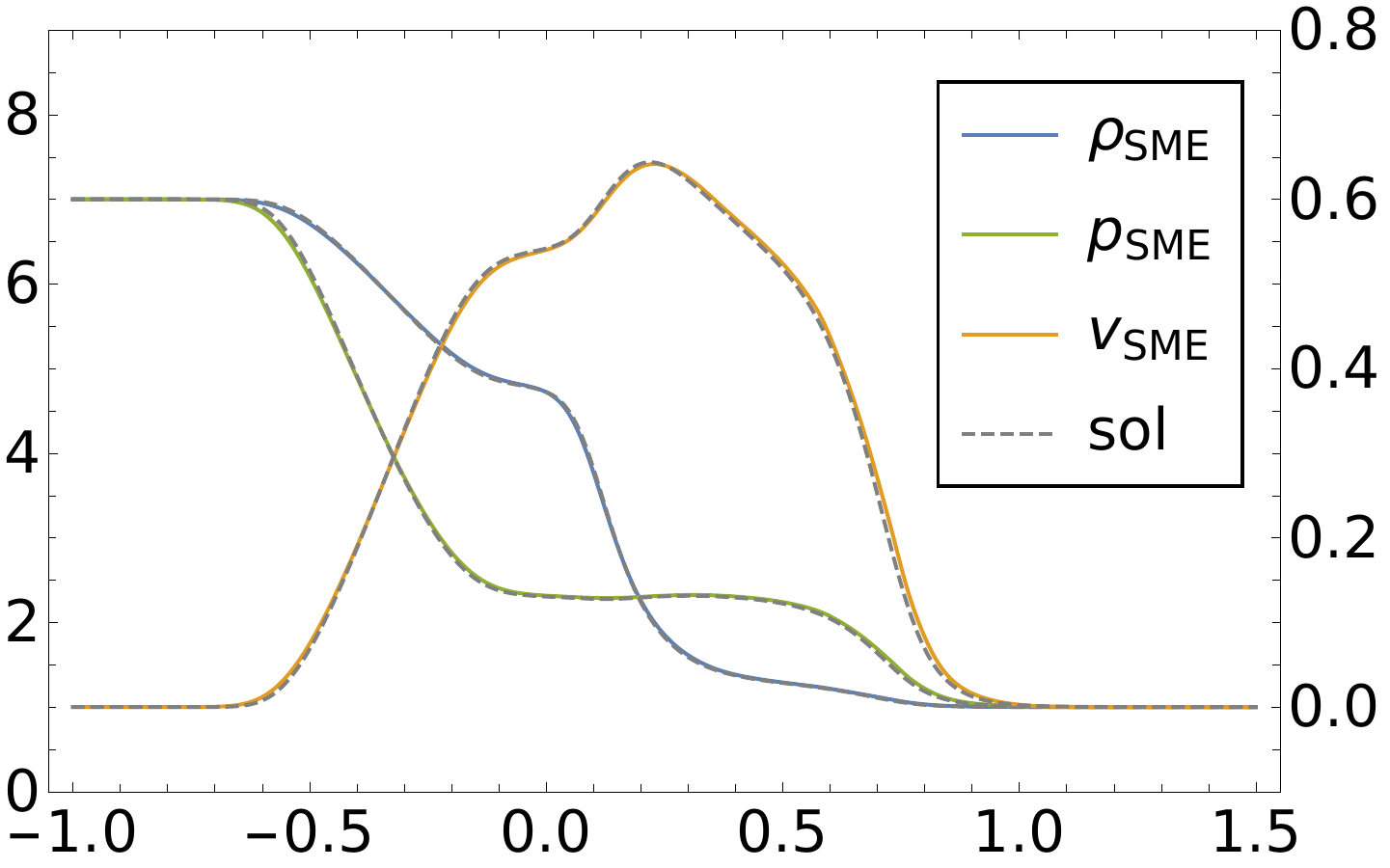}
 \caption{n=13 {on} $[-4,4]$, $\textrm{Kn}=0.05$.}
 \label{base5}
\end{subfigure}
\hfill
\begin{subfigure}[t]{0.49\textwidth}
 \centering
 \includegraphics[width=\textwidth]{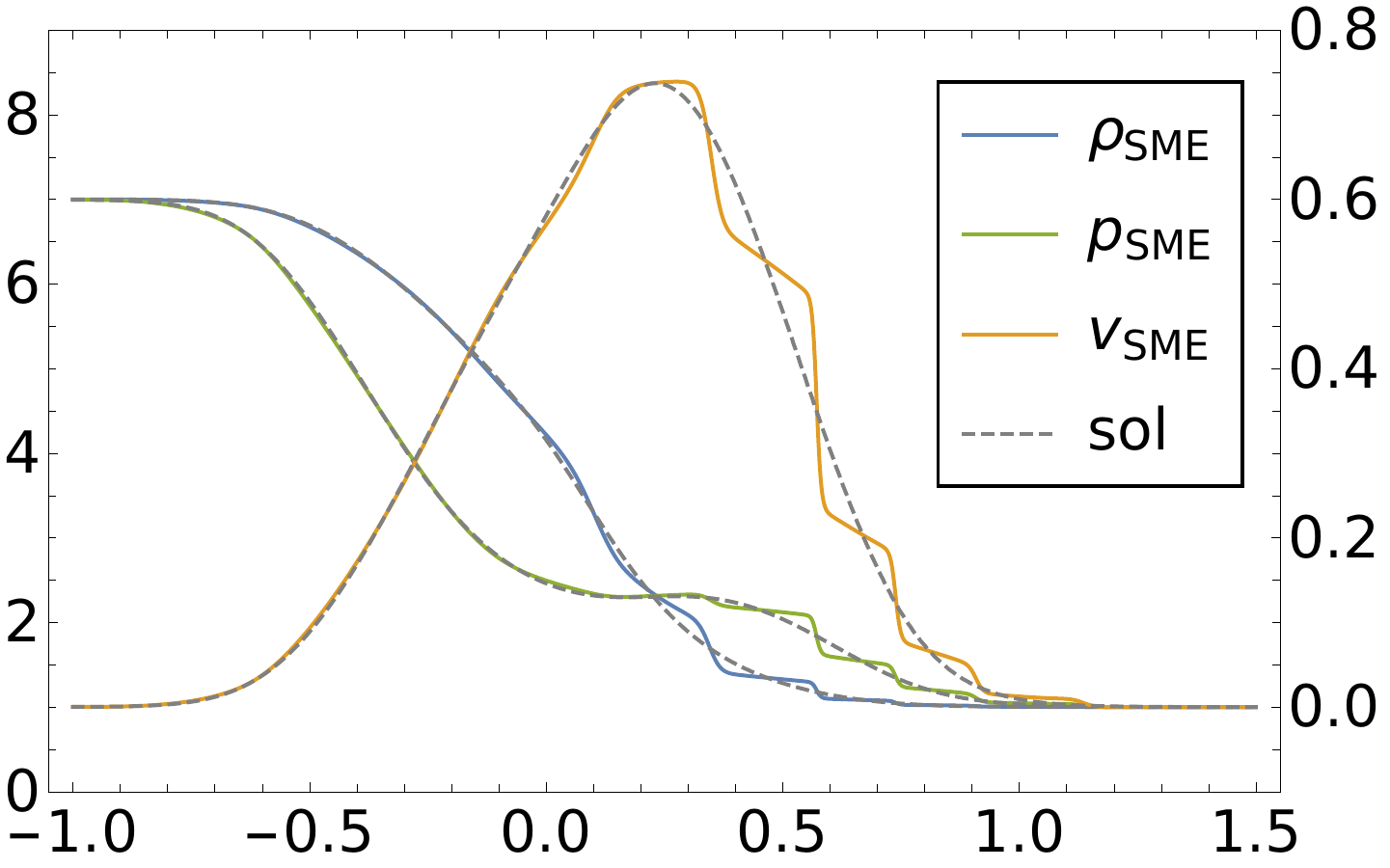}
 \caption{n=13 {on} $[-4,4]$, $\textrm{Kn}=0.5$.}
 \label{base6}
\end{subfigure}
\caption{Shock tube for different $\xi$-ranges with constant spline distance $\Delta \xi=\frac{2}{3}$. Left column (a),(c),(e) $\textrm{Kn}=0.05$, right column (b),(d),(f) $\textrm{Kn}=0.5$. (a),(b) $n=7$, (c),(d) $n=10$, (e),(f) $n=13$.}
\label{resultsDifferentBases}
\end{figure}

\subsubsection{{Effect of the number of splines}}
In Figures \ref{kase} and \ref{resultsDifferentN} the number of splines is increased while keeping the $\xi$-range constant. For the smaller Knudsen number in Figure \ref{kase}, the solution is already very accurate for $4$ splines (corresponding to $n=7$ equations) and does not change much for larger numbers of splines or equations. For the larger Knudsen number in Figure \ref{resultsDifferentN} however, we can see the improvement clearly by observing smaller subshocks and a better accuracy in all regions of the domain.


\begin{figure}[H]
\centering
\begin{subfigure}[t]{0.49\textwidth}
\centering
\includegraphics[width=\textwidth]{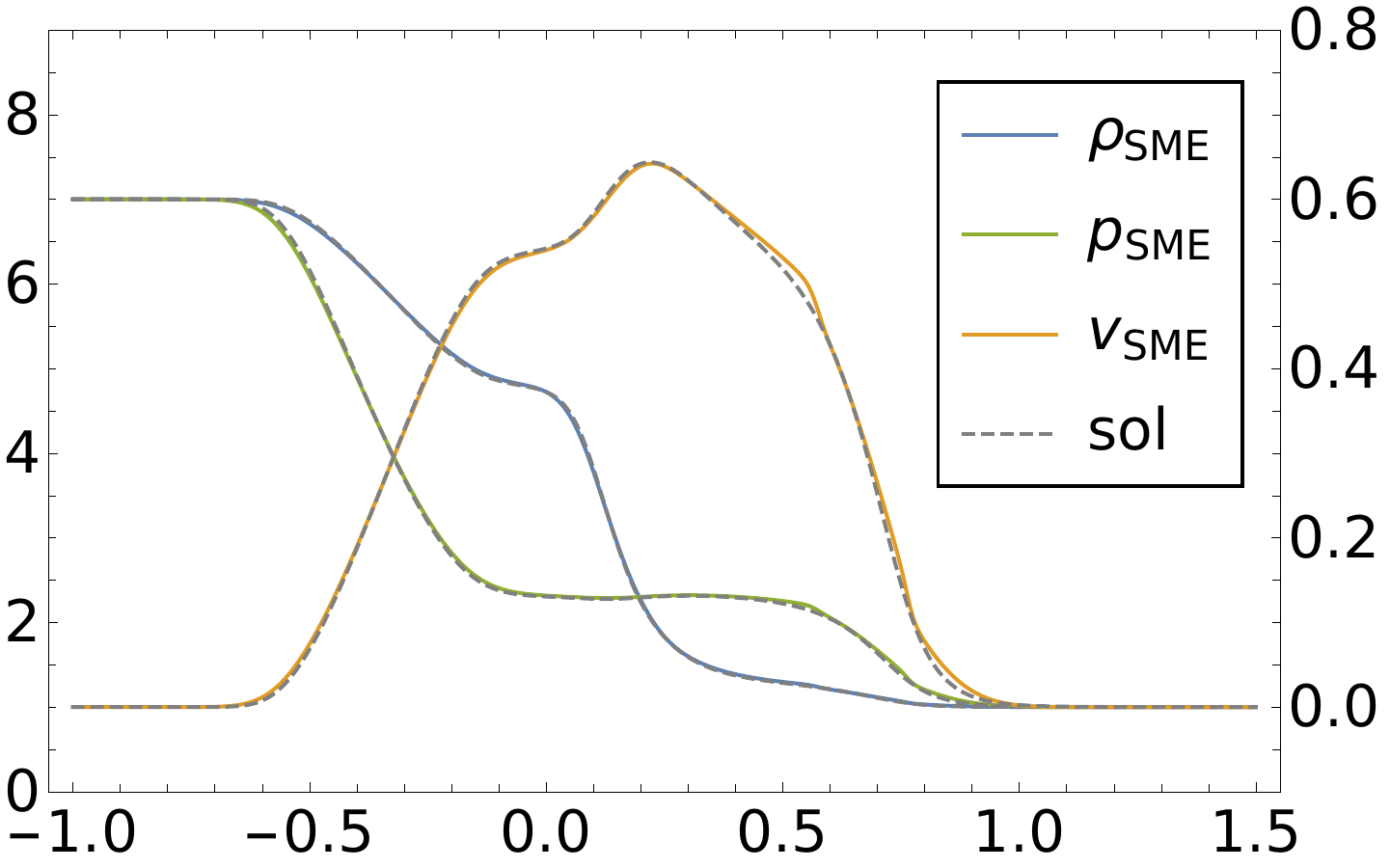}
\caption{$n=7$.}
\label{kase2}
\end{subfigure}
\hfill
\begin{subfigure}[t]{0.49\textwidth}
\centering
\includegraphics[width=\textwidth]{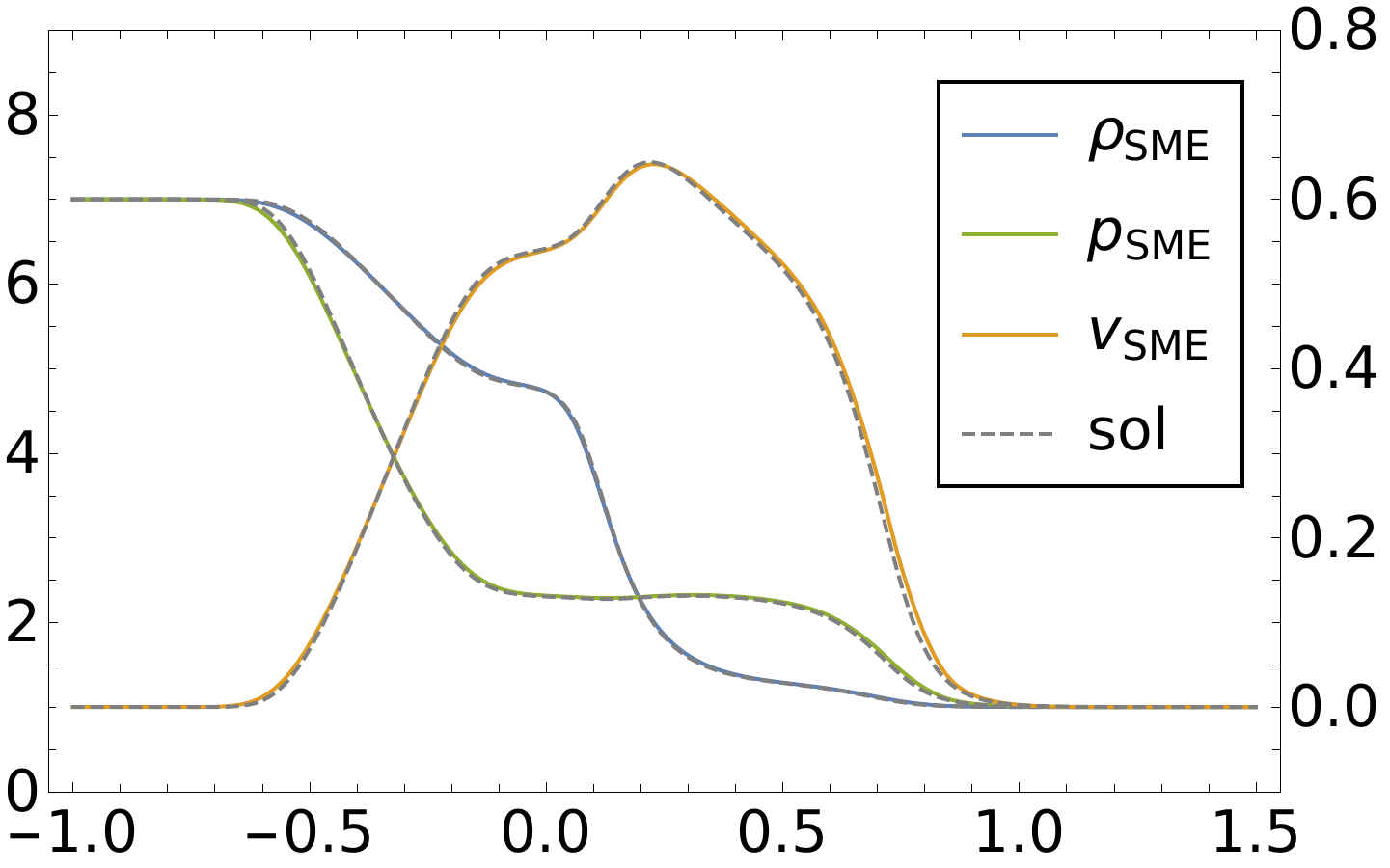}
\caption{$n=13$.}
\label{kase4}
\end{subfigure}
\caption{Shock tube for different $n$, $[\xi_{min},\xi_{max}]=[-4,4]$, $\textrm{Kn}=0.05$, $k=1$. (a) $n=7$, (b) $n=13$.}
\label{kase}
\end{figure}

\begin{figure}[H]
 \centering
\begin{subfigure}[t]{0.49\textwidth}
 \centering
 \includegraphics[width=\textwidth]{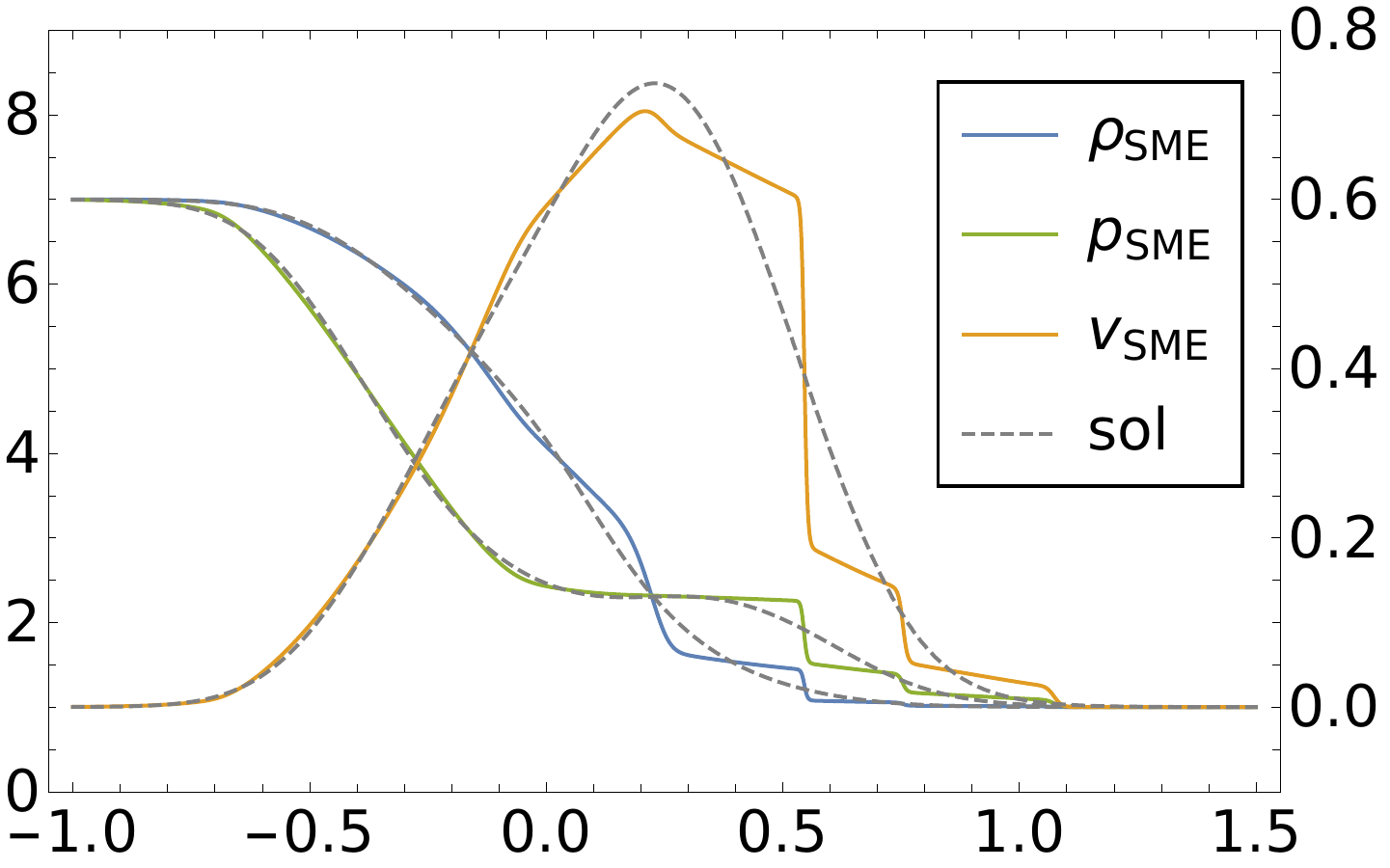}
 \caption{$n=7$.}
 \label{count2}
\end{subfigure}
\hfill
\begin{subfigure}[t]{0.49\textwidth}
 \centering
 \includegraphics[width=\textwidth]{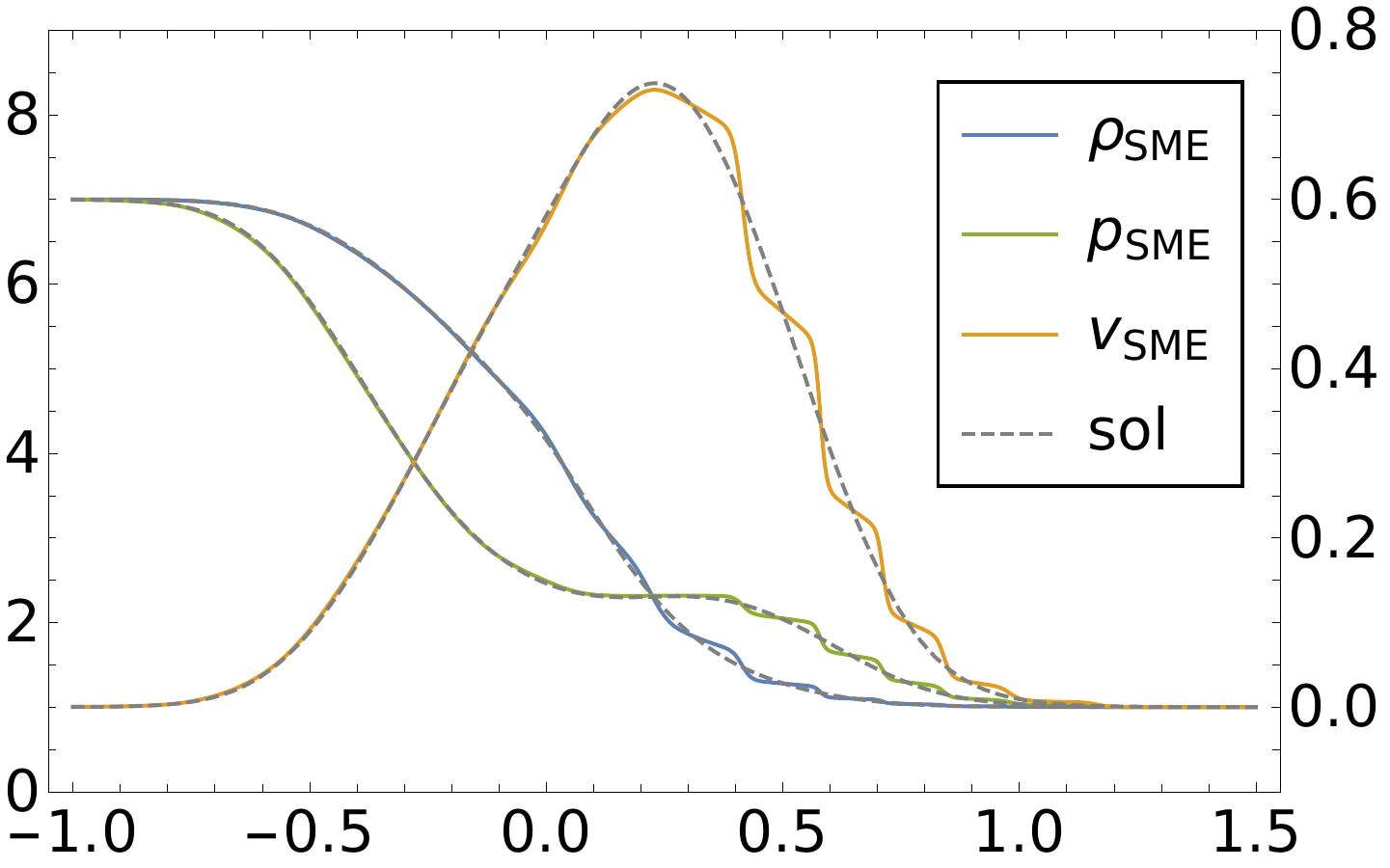}
 \caption{$n=13$.}
 \label{count4}
\end{subfigure}
\caption{Shock tube for different $n$, $[\xi_{min},\xi_{max}]=[-4,4]$, $\textrm{Kn}=0.5$, $k=1$. (a) $n=7$, (b) $n=13$.}
\label{resultsDifferentN}
\end{figure}

In the error plots in Figure \ref{errorPlotsDifferentOrders} the relative errors are defined as
\begin{equation}
  \err_{\rho}=\frac{\int | \rho_{SME}-\rho_{sol} | dx }{ \int | \rho_{sol} | dx},
\err_{p}=\frac{\int | p_{SME}-p_{sol} | dx }{ \int | p_{sol} | dx},
\textrm{{and}} \: \err_{v}=\frac{\int | v_{SME}-v_{sol} | dx }{ \int | v_{sol} | dx},
\end{equation}
and the error is plotted for different spline orders $k=1,2,3$ in each figure. Additionally, we consider the error evolution on different $\xi$-ranges $[\xi_{min},\xi_{max}]=[-3,3]$ and $[\xi_{min},\xi_{max}]=[-4,4]$.

In all settings the error for density and pressure is a factor 10 smaller than the velocity error in accordance with the literature \cite{Koellermeier2017b}. The errors clearly decrease with increasing number of splines. For about $10$ basis functions, the error can be reduced up to $0.2\%-0.4\%$ for $\rho$ and $p$ and to around $4\%$ for $v$, which is very accurate when considering the small number of equations involved. For $[\xi_{min},\xi_{max}]=[-2,2]$ the $\xi$-range is too small and no convergence is obtained (not shown here). A higher spline order leads to a slightly smaller error while showing the same trend as first order.

}
\begin{figure}[H]
 \centering
\begin{subfigure}[t]{0.49\textwidth}
 \centering
 \includegraphics[width=\textwidth]{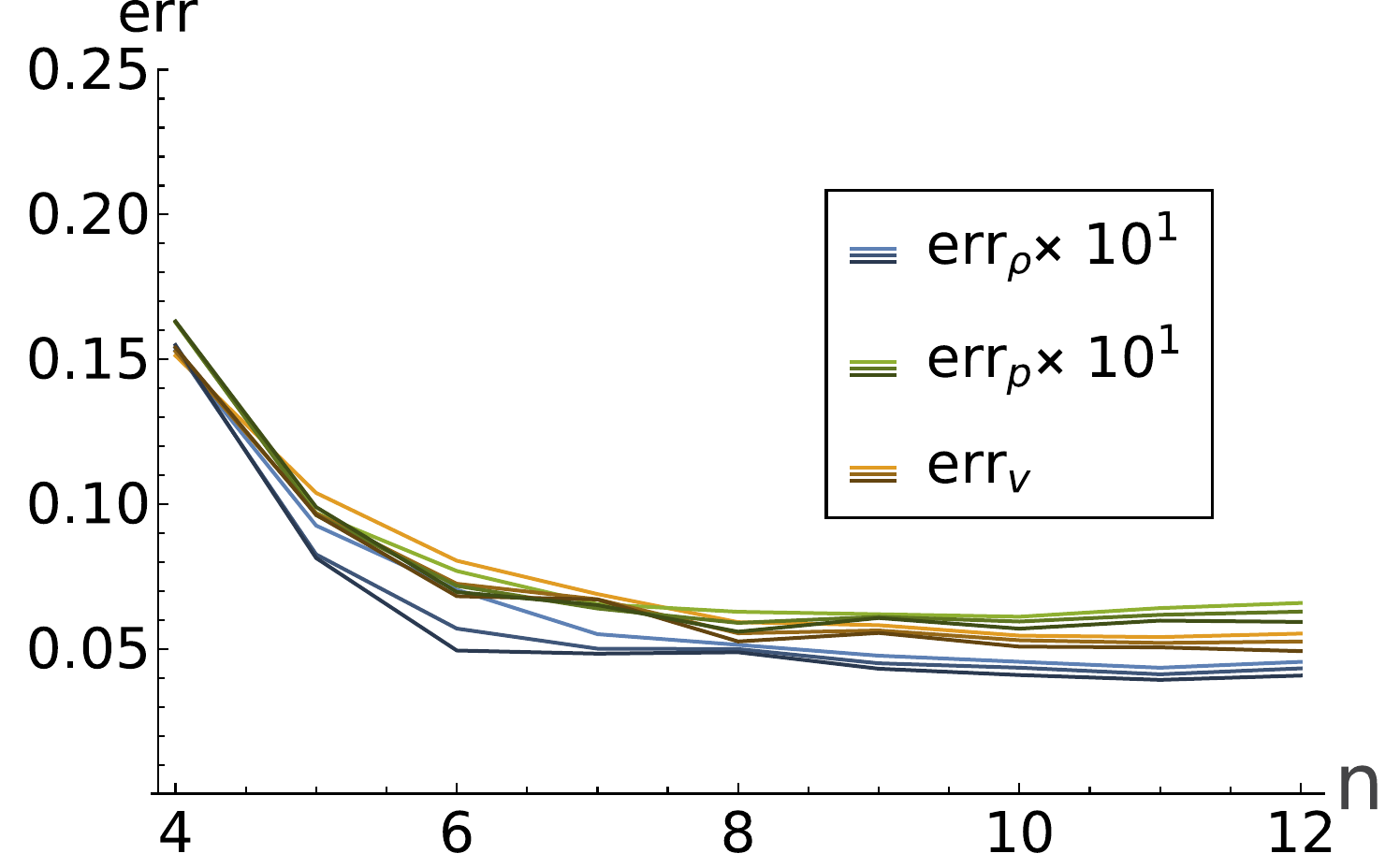}
 \caption{$[\xi_{min},\xi_{max}]=[-3,3]$, $k=1,2,3$ {(from bright to dark)}.}
 \label{error4}
\end{subfigure}
\hfill
\begin{subfigure}[t]{0.49\textwidth}
 \centering
 \includegraphics[width=\textwidth]{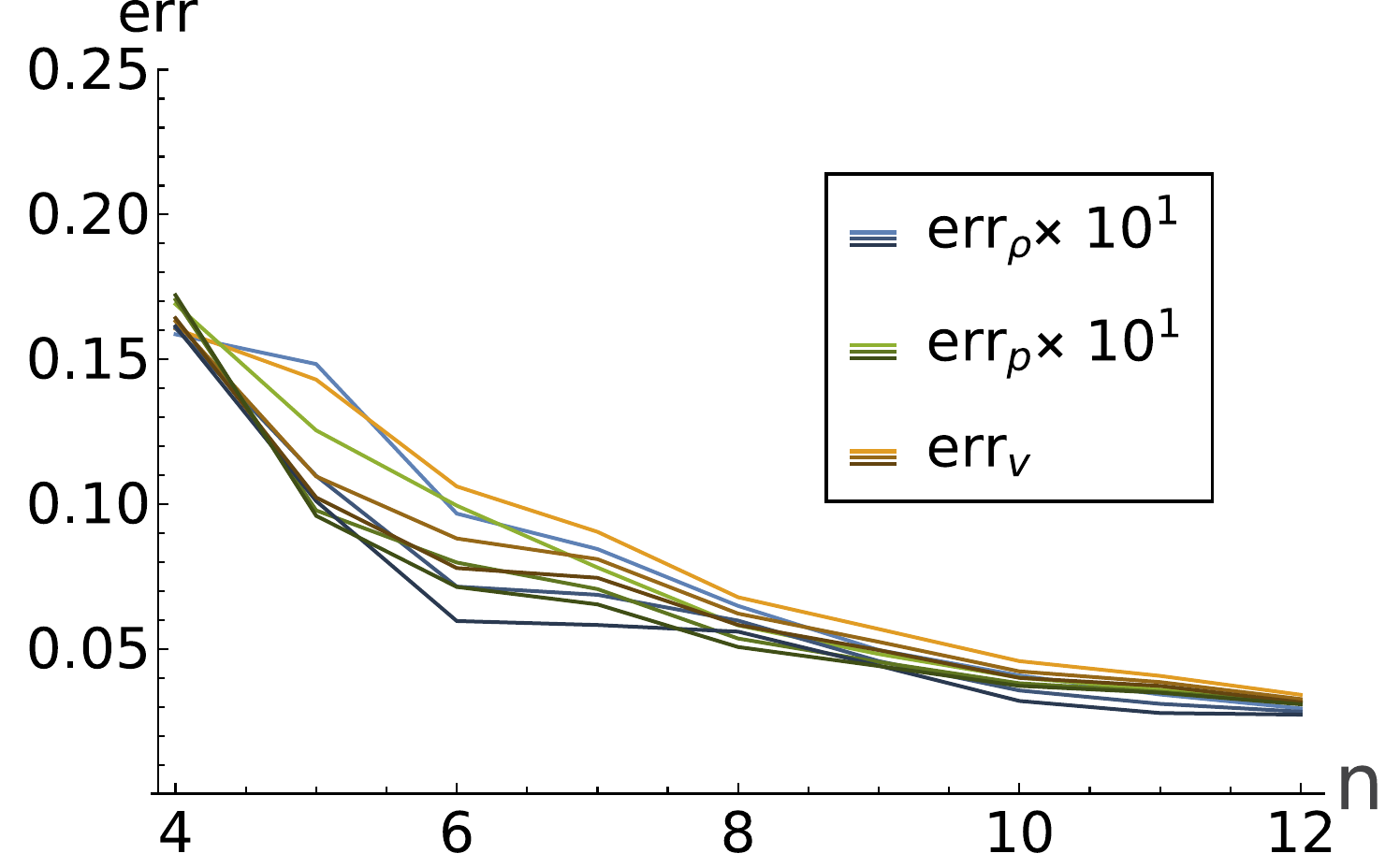}
 \caption{$[\xi_{min},\xi_{max}]=[-4,4]$, $k=1,2,3$ {(from bright to dark)}.}
 \label{error6}
\end{subfigure}
\caption{Simulation error depending on $n$ for different orders, $[\xi_{min},\xi_{max}]$, $\textrm{Kn}=0.5$. (a) $[\xi_{min},\xi_{max}]=[-3,3]$, (b) $[\xi_{min},\xi_{max}]=[-4,4]$.}
\label{errorPlotsDifferentOrders}
\end{figure}

\subsubsection{Computational efficiency}
To access the computational efficiency of the new spline models, we measured the run time for selected simulations on the interval $[\xi_{min},\xi_{max}]=[-4,4]$. Independent of the Knudsen number and the spline order, we have obtained the run time of approximately $t=19,28,44,77$ seconds for $n=4,7,10,13$, respectively. This yields an approximately linear increase in the runtime $t$ with increasing $n$. In Figure \ref{fig:errorAbhVonT} we plotted the error in the velocity variable for different Knudsen numbers depending on the run time that had been invested. This displays the computational efficiency of the models for different Knudsen numbers. As expected, for higher Knudsen numbers we need more equations and thus more computational effort to achieve accurate solutions, whereas relatively coarse solutions can be obtained in comparably short times.
There is a remaining error of less than $1 \%$ which is due to the spatial discretization within the numerical simulation, which is not the topic of this paper. Most importantly, the model error is reduced and is of the order of $1 \%$ for only $n=13$ equations.
The respective error plots for the density and pressure variables give the same result (not depicted here).

A convergence plot for the relative velocity error with respect to the Knudsen number is depicted in Fig \ref{fig:errorAbhVonKn}. Again, we observe that for larger Knudsen numbers more equations are necessary for an accurate simulation result. It is clearly visible that regardless of the number of equations $n$ the error converges to a small plateau of approximately $1\%$ for decreasing Knudsen number. As mentioned before, the remaining error is a result of numerical errors during the simulation and not a model error. A study of the other two macroscopic variables, density $\rho$ and pressure $p$, gives the same result.

\begin{figure}[H]
\centering
\begin{subfigure}[t]{0.49\textwidth}
 \centering
 \includegraphics[width=\textwidth]{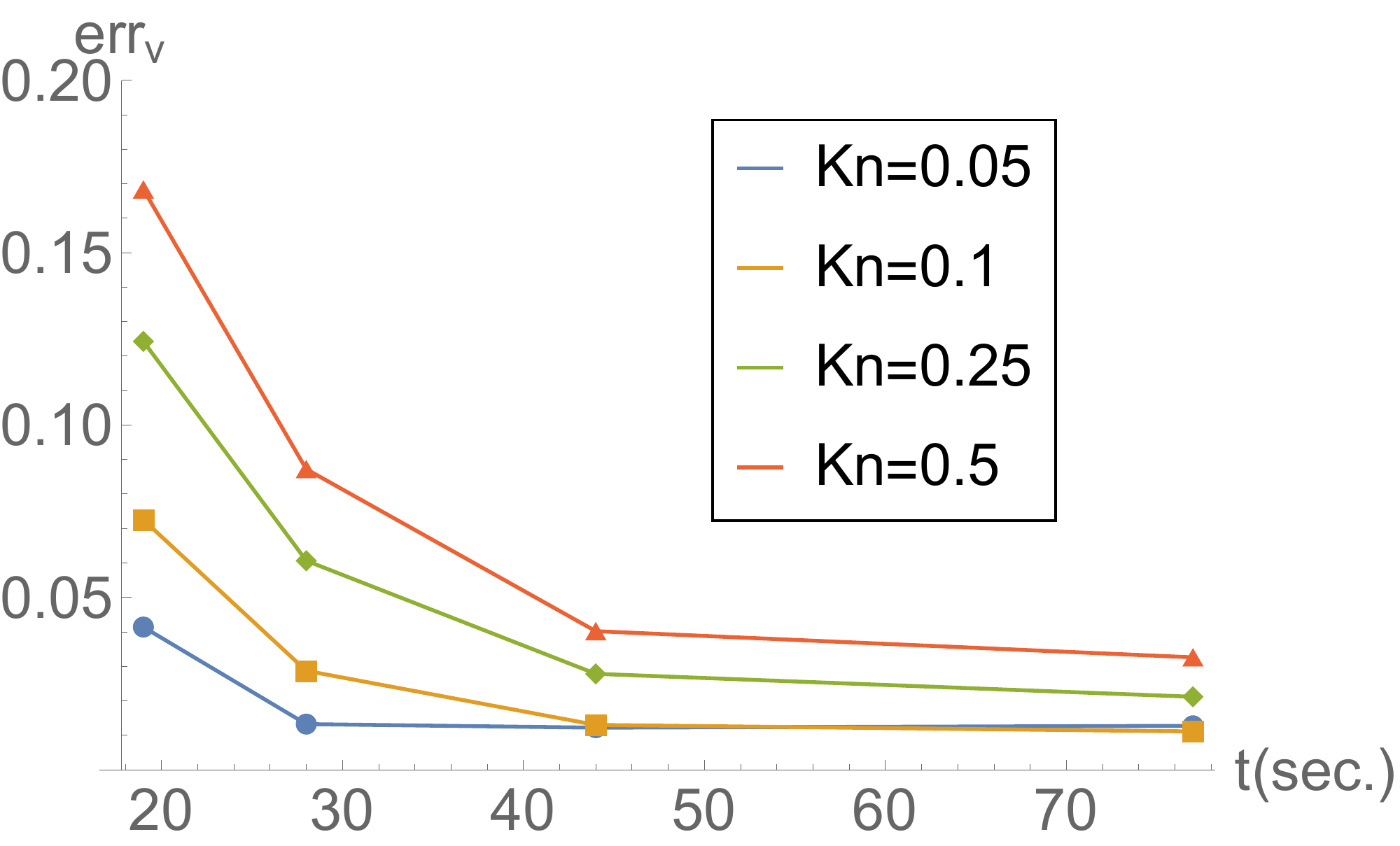}
 \caption{Error depending on $t$ for different Kn from~$0.05$ to~$0.5$.}
\label{fig:errorAbhVonT}
\end{subfigure}
\hfill
\begin{subfigure}[t]{0.49\textwidth}
 \centering
 \includegraphics[width=\textwidth]{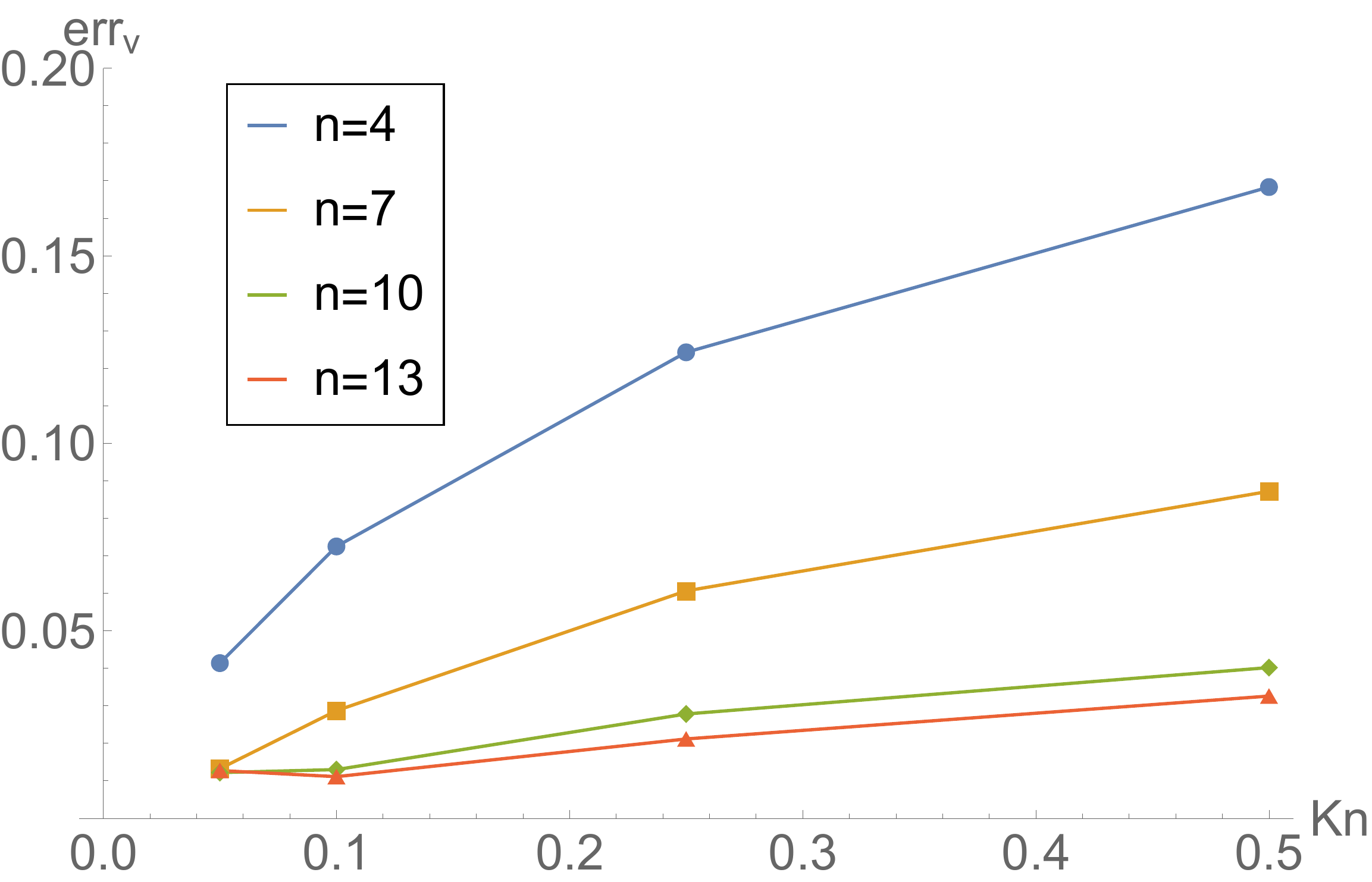}
 \caption{Error depending on Kn for~$n=3,7,10,13$.}
\label{fig:errorAbhVonKn}
\end{subfigure}
\caption{Simulation error in the velocity variable $v$, $k=1$, $[\xi_{min},\xi_{max}]=[-4,4]$. (a) depending on $t$, (b) depending on Kn.}
\end{figure}

\subsubsection{Comparison with existing QBME model}
In comparison with the reference solution obtained by a discrete velocity method, our method yields a fast and accurate solution. Now we compare the new SME with the existing Quadrature-Based Moment Equations (QBME), e.g. as described in \cite{Koellermeier2017b}. We compare using the same number of equations in Figure \ref{hermiteCmp} and we observe that the results look surprisingly similar and even deviations from the reference solution occur in similar positions. Both methods are very accurate for the case $\textrm{Kn}=0.05$. In Figure \ref{hermite4} the spline method appears to be more accurate. The SME model thus yields more accurate solutions than the existing moment model in this test case.
\begin{figure}[H]
 \centering
 \begin{subfigure}[t]{0.49\textwidth}
 \centering
 \includegraphics[width=\textwidth]{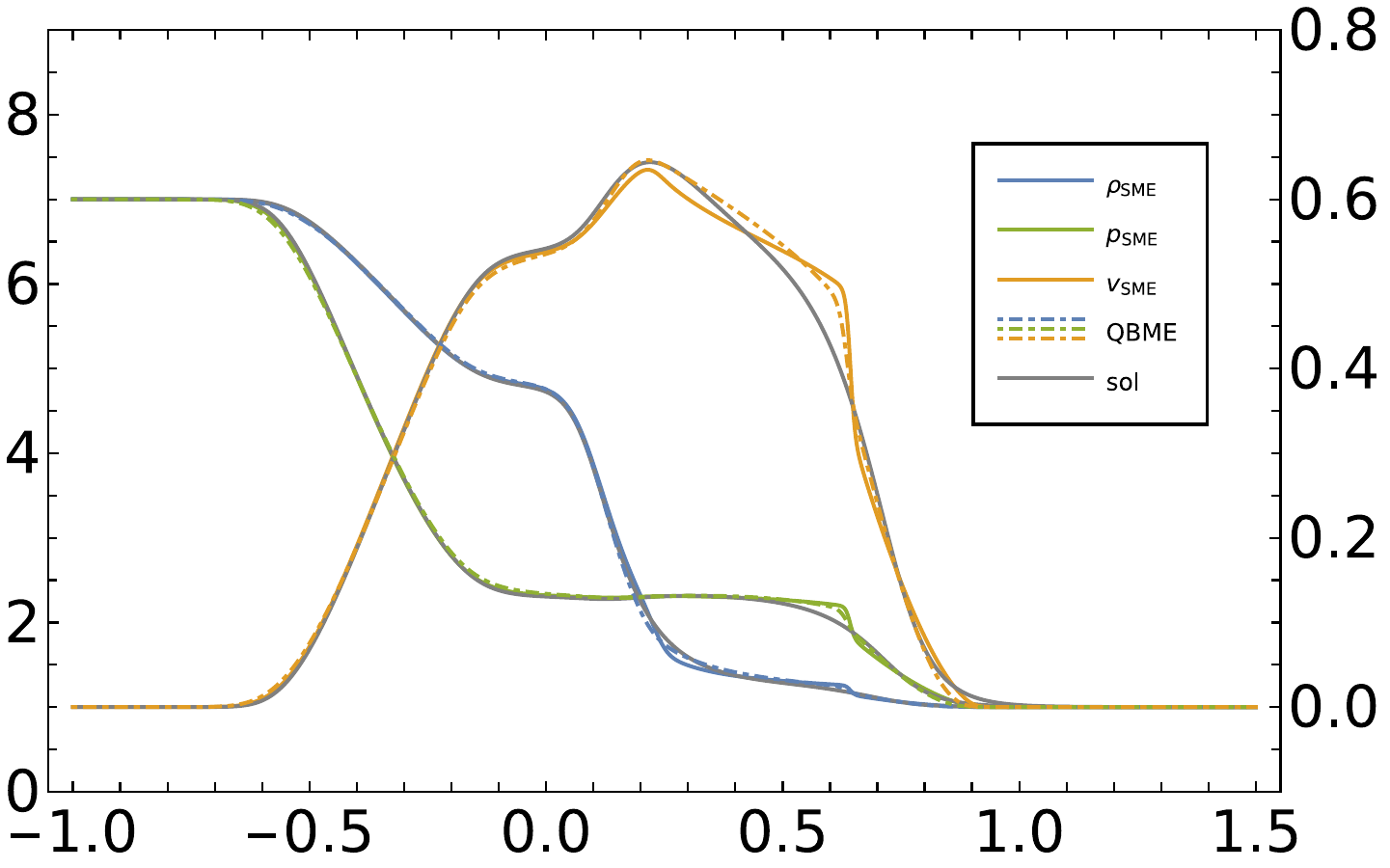}
 \caption{$n=5$, $\textrm{Kn}=0.05$.}
 \label{hermite1}
\end{subfigure}
\hfill
\begin{subfigure}[t]{0.49\textwidth}
 \centering
 \includegraphics[width=\textwidth]{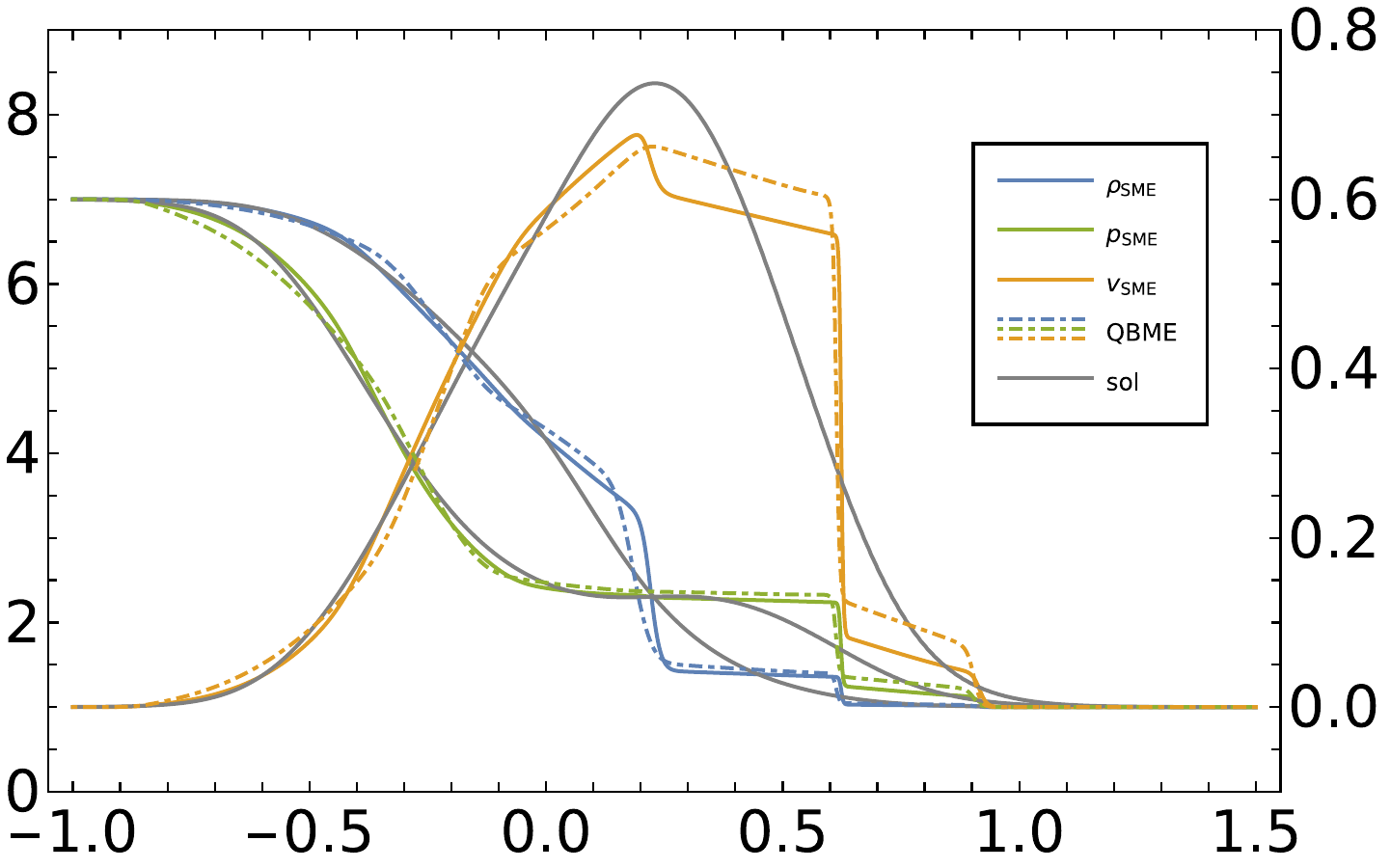}
 \caption{$n=5$, $\textrm{Kn}=0.5$.}
 \label{hermite2}
\end{subfigure}
\hfill
\hfill
\begin{subfigure}[t]{0.49\textwidth}
 \centering
 \includegraphics[width=\textwidth]{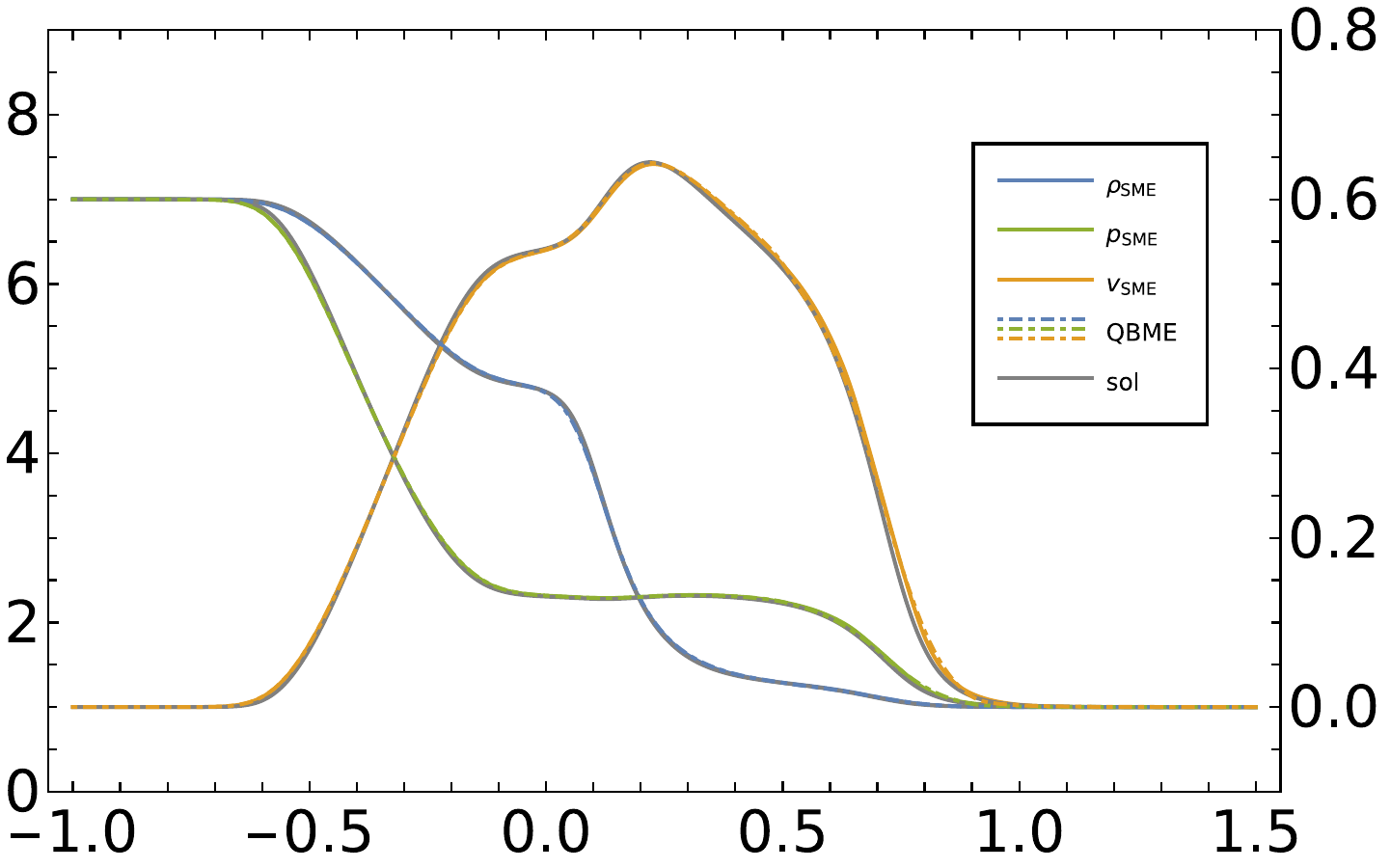}
 \caption{$n=11$, $\textrm{Kn}=0.05$.}
 \label{hermite3}
\end{subfigure}
\hfill
\begin{subfigure}[t]{0.49\textwidth}
 \centering
 \includegraphics[width=\textwidth]{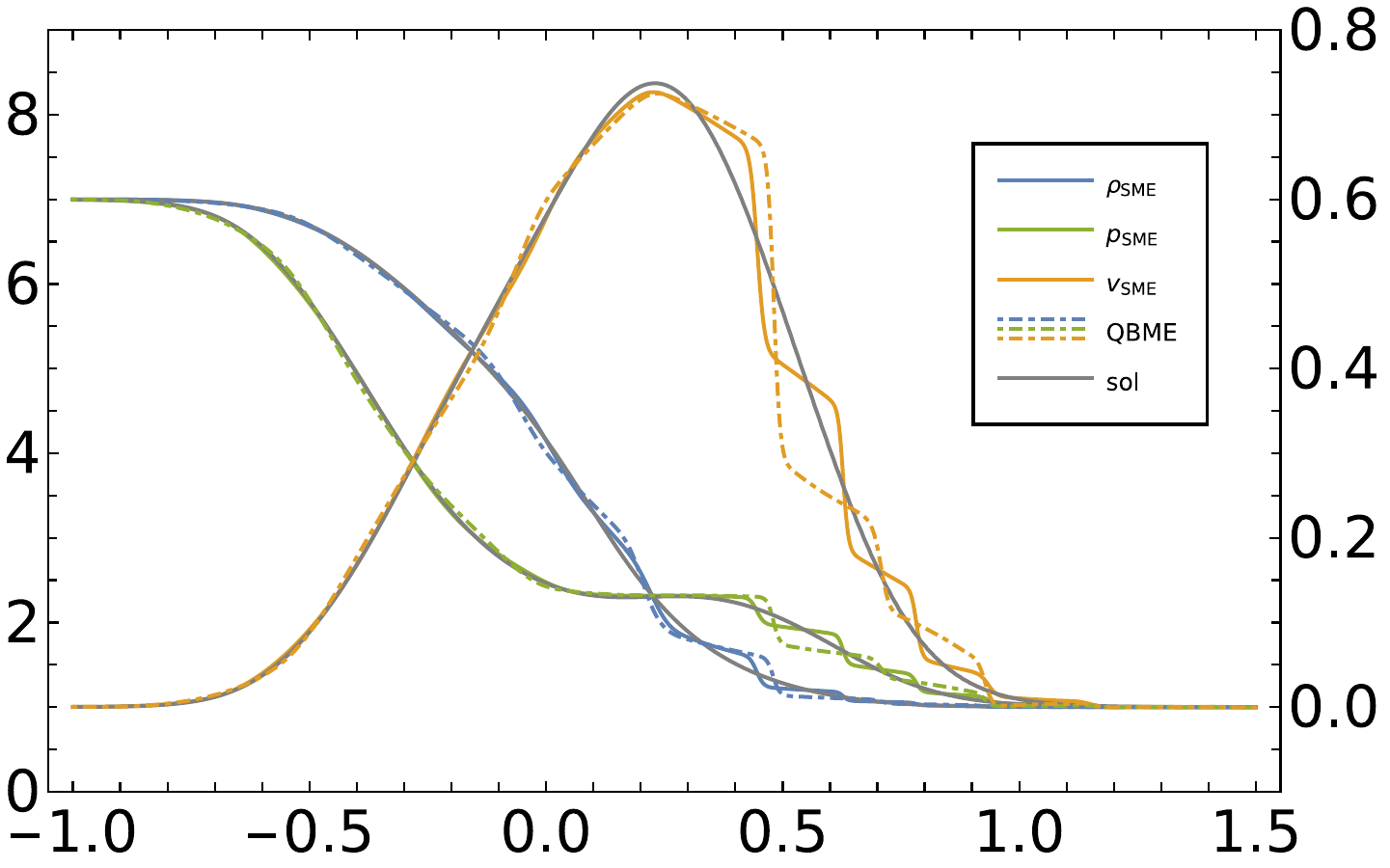}
 \caption{$n=11$, $\textrm{Kn}=0.5$.}
 \label{hermite4}
\end{subfigure}
\caption{Shock tube comparing SME to QBME model. $k=1$, $[\xi_{min},\xi_{max}]=[-4,4]$. Left column (a),(c) $\textrm{Kn}=0.05$, right column (b),(d) $\textrm{Kn}=0.5$. (a),(b) $n=5$, (c),(d) $n=11$.}
\label{hermiteCmp}
\end{figure}

Summarizing the shock tube test case, we saw that the new SME model approximates the solution with good accuracy and even outperforms the existing QBME model for the same number of equations.

\subsection{Symmetric two-beam problem}
In this test case from \cite{schaerer2015}, a BGK collision operator with a constant relaxation time $\tau = \textrm{Kn} \cdot t_{\textrm{END}}$ is used to model collisions similar to the shock tube test case.

The initial Riemann data is given by \eqref{e:2beam_IC}
\begin{equation}
    \vect{u}_M^L = \left( 1,0.5,1,0,\ldots,0\right)^T, \quad \quad \vect{u}_M^R = \left( 1,-0.5,1,0,\ldots,0\right)^T,
    \label{e:2beam_IC}
\end{equation}
and models two colliding Maxwellian distributed particle beams. We note that this test case is especially challenging for any type of polynomial ansatz as it is difficult to represent the analytical solution using a polynomial expansion. In the free streaming case $\textrm{Kn} = \infty$, the analytical solution is a sum of two Maxwellians distributions according to \cite{schaerer2015}.

The numerical tests are performed on the computational domain $[-10,10]$, discretized using $4000$ points. The end time is $t_{\textrm{END}}=0.3$ using a constant CFL number of approximately $0.5$ for all tests.

Tests are shown for $\textrm{Kn} = 0.1$ representing a small Knudsen number, $\textrm{Kn} = 1$ for a relatively large Knudsen number and $\textrm{Kn} = \infty$ leading to vanishing right-hand side and very sharp profiles. We show results for the SME model using $n=10$ with order $k=1$, and $[\xi_{min},\xi_{max}]=[-4,4]$, which was identified as accurate in the previous shock tube test case. The results are similar for other SME models using more equations, but we omit the results here for conciseness. We use the hyperbolic QBME model with $M+1=10$ equations for comparison. We note that the QBME model yields similar results as the HME model from \cite{Cai2013} for this test case. A discrete velocity method is used as reference solution. It was taken from \cite{schaerer2015} and computed using 2000 cells in physical space and 600 variables for the discretization of the microscopic velocity space. Note that the DVM method is computationally much more expensive in comparison to the lower-dimensional moment models. As in the literature, we display pressure $p=\rho\theta$ and the normalized heat flux $\bar{q}$ (see \ref{e:Q}), which can be computed for the moment models using
\begin{equation} \label{e6:normalized_heat_flux}
    \bar{q} = \frac{6 f_3}{\rho \sqrt{\theta}^3},
\end{equation}
while in the case of this specific spline model, we derive the heat flux to
\begin{equation}
    \begin{split}
        \label{e6:normalized_heat_flux10}
        \bar{q} \approx 0.000847777 \kappa_1 + 0.0129408 \kappa_2 +  0.0836094 \kappa_3 + 0.170286 \kappa_4 \\
         + 0.0836094 \kappa_5 + 0.0129408 \kappa_6 + 0.000847777 \kappa_7.
    \end{split}
\end{equation}

The results in Figure \ref{2beam} show a clear convergence of the new SME model and at least similar accuracy in comparison to the QBME model, which uses the same number of equations. For the collisionless test case with $\textrm{Kn} = \infty$, the step-like structure is due to the wave structure of the problem with at most $n=10$ propagation speeds, see remark \ref{remark_hyp}. We can see that the propagation speeds are slightly different than those of the QBME model, because of the different ansatz. For $\textrm{Kn} = 1$, the solution gets closer to the DVM reference, while still showing a similar accuracy to the QBME model. For the smallest Knudsen number $\textrm{Kn} = 0.1$, the new SME model has converged to the reference solution. Notably, also the heat flux is approximated with high precision.
\begin{figure}[H]
 \centering
 \begin{subfigure}[t]{0.49\textwidth}
 \centering
 \includegraphics[width=\textwidth]{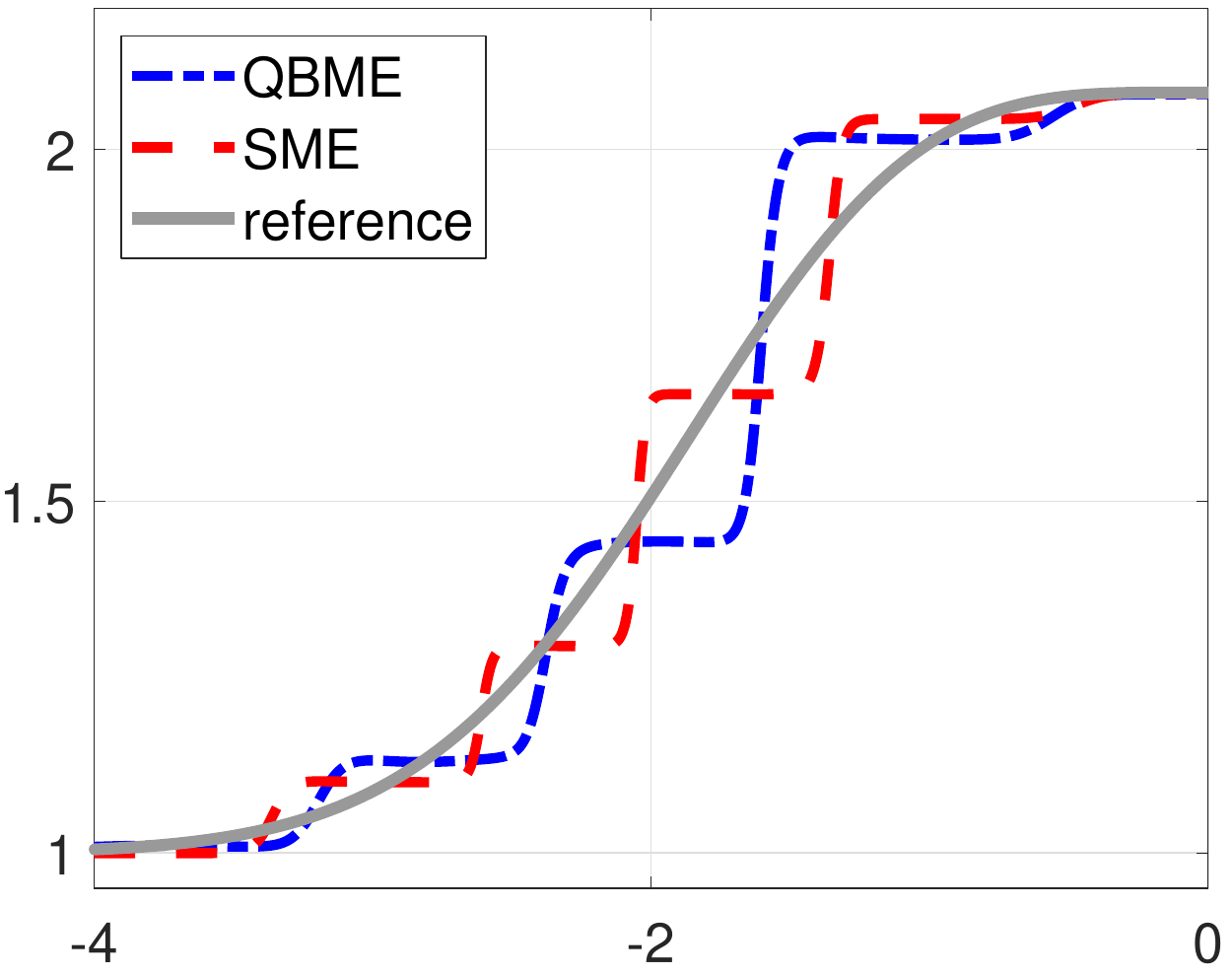}
 \caption{$Kn = \infty$, $p$}
\end{subfigure}
\hfill
 \begin{subfigure}[t]{0.49\textwidth}
 \centering
 \includegraphics[width=\textwidth]{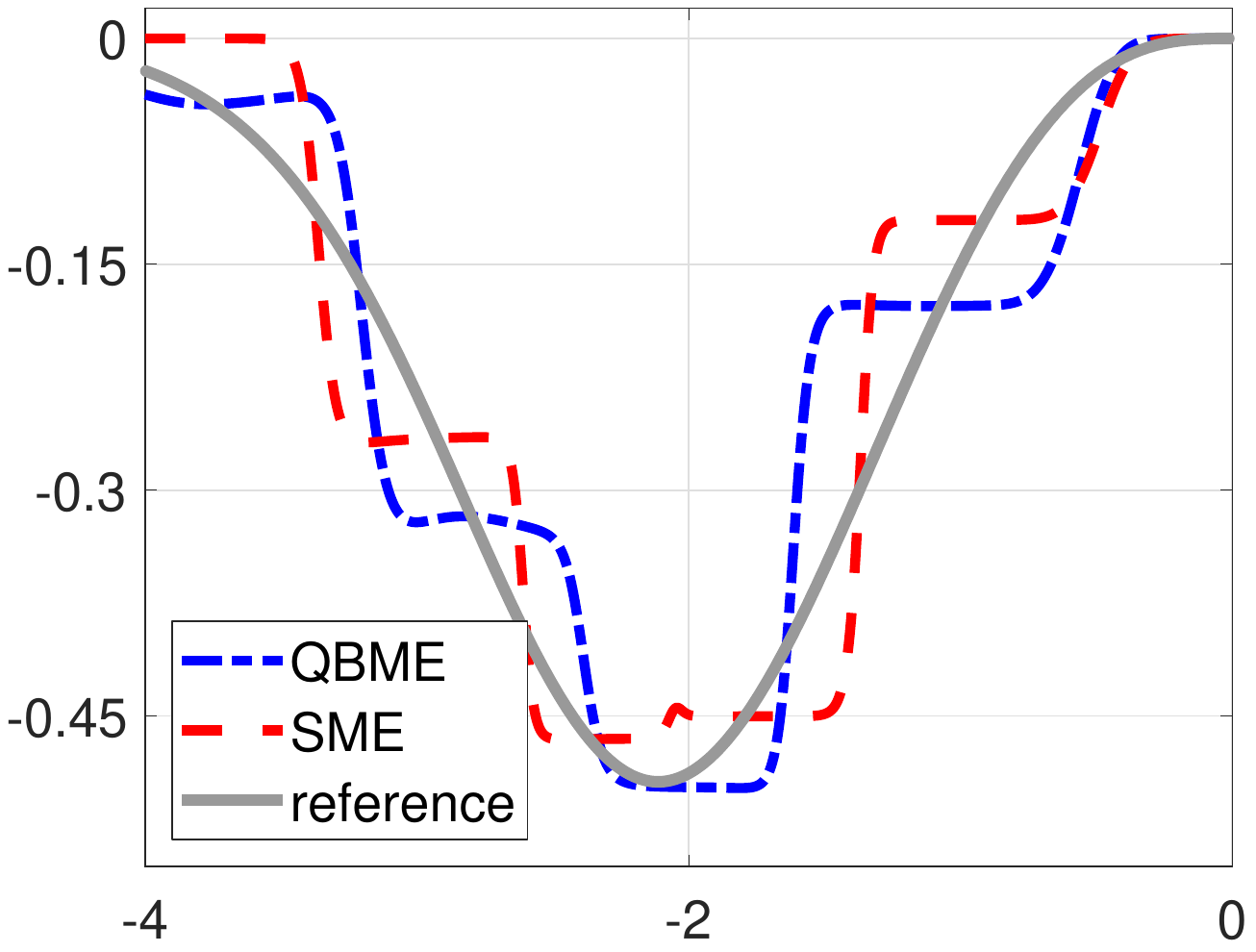}
 \caption{$Kn = \infty$, $\bar{q}$}
\end{subfigure}
\hfill
\begin{subfigure}[t]{0.49\textwidth}
 \centering
 \includegraphics[width=\textwidth]{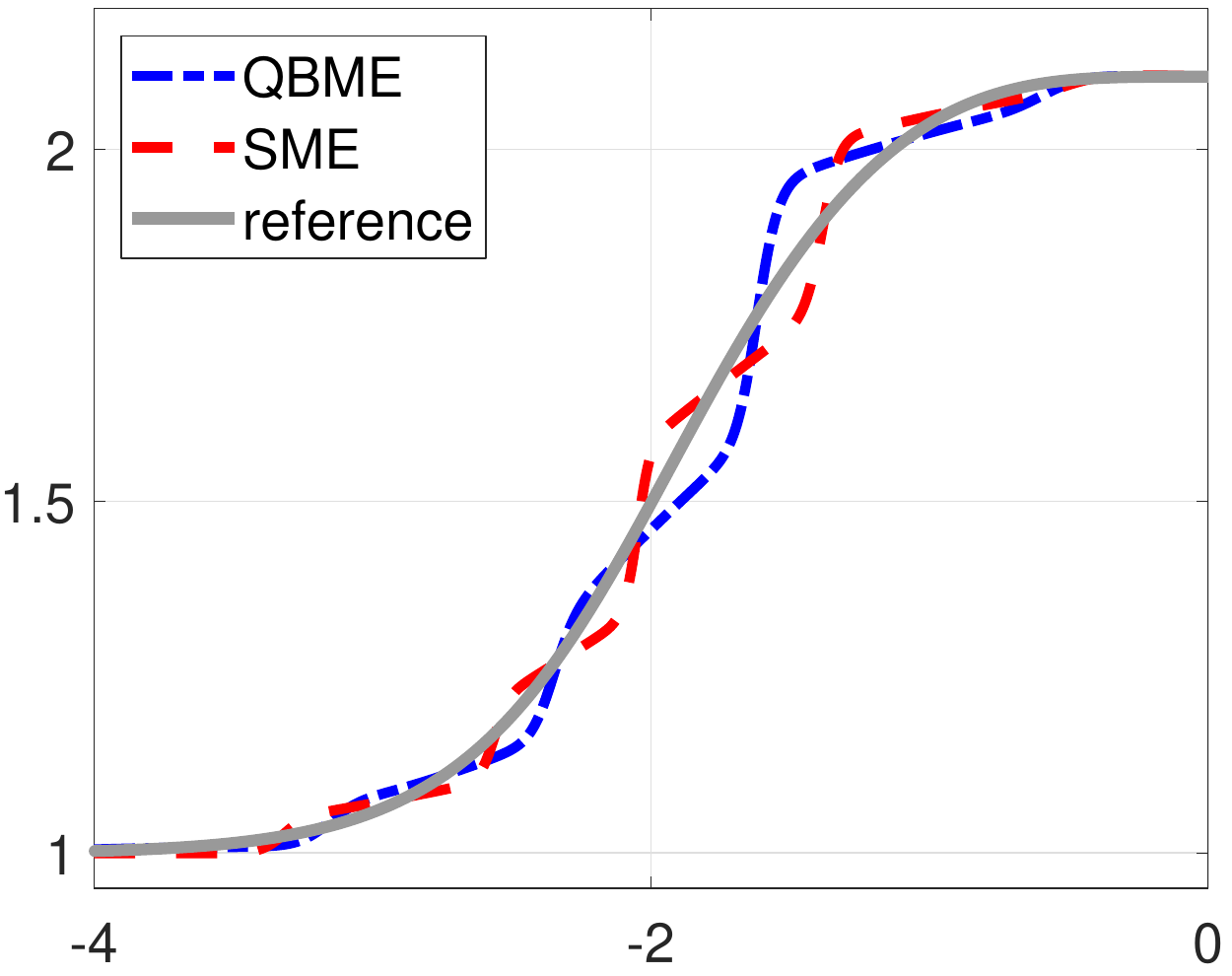}
 \caption{$Kn = 1$, $p$}
\end{subfigure}
\hfill
\begin{subfigure}[t]{0.49\textwidth}
 \centering
 \includegraphics[width=\textwidth]{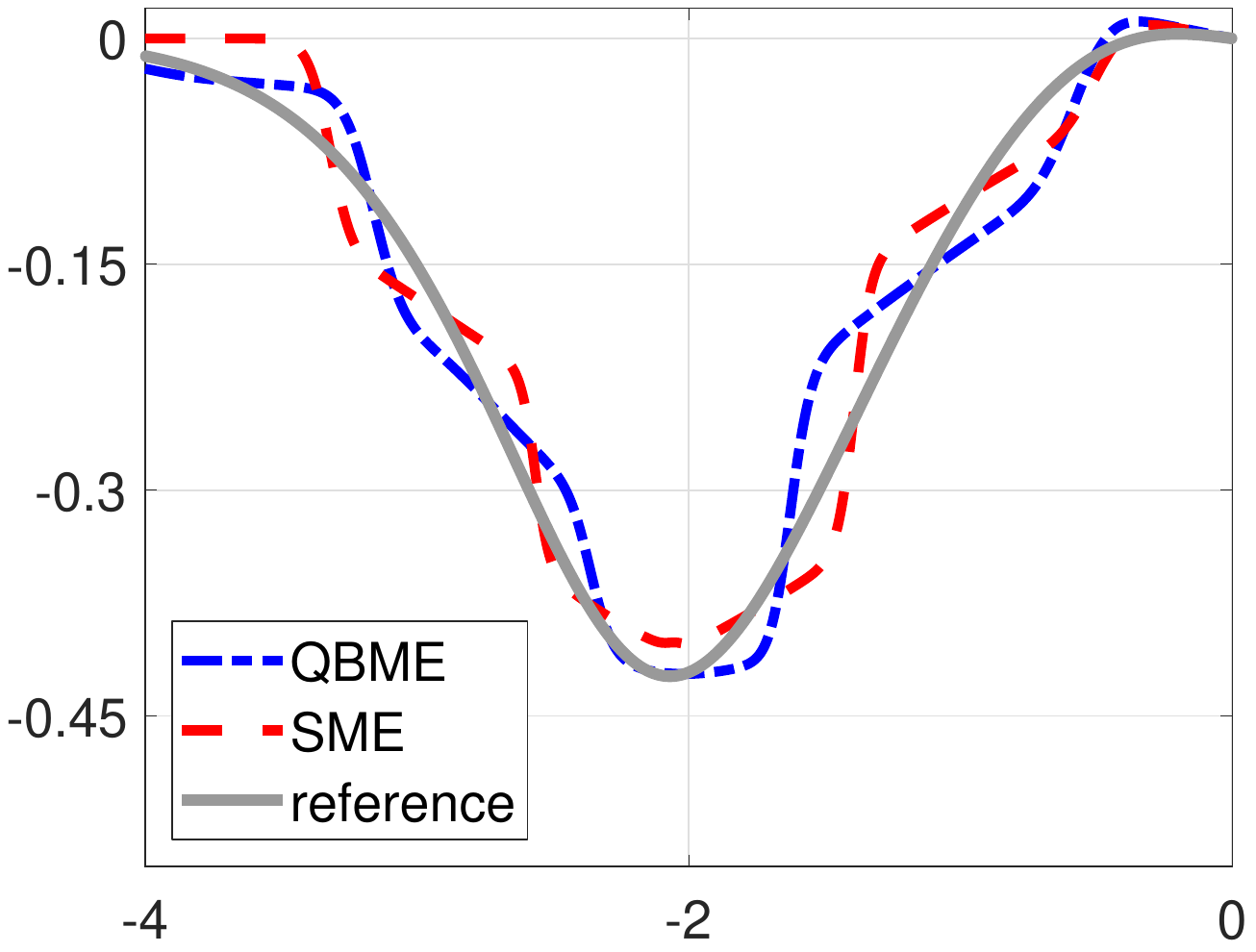}
 \caption{$Kn = 1$, $\bar{q}$}
\end{subfigure}
\hfill
 \begin{subfigure}[t]{0.49\textwidth}
 \centering
 \includegraphics[width=\textwidth]{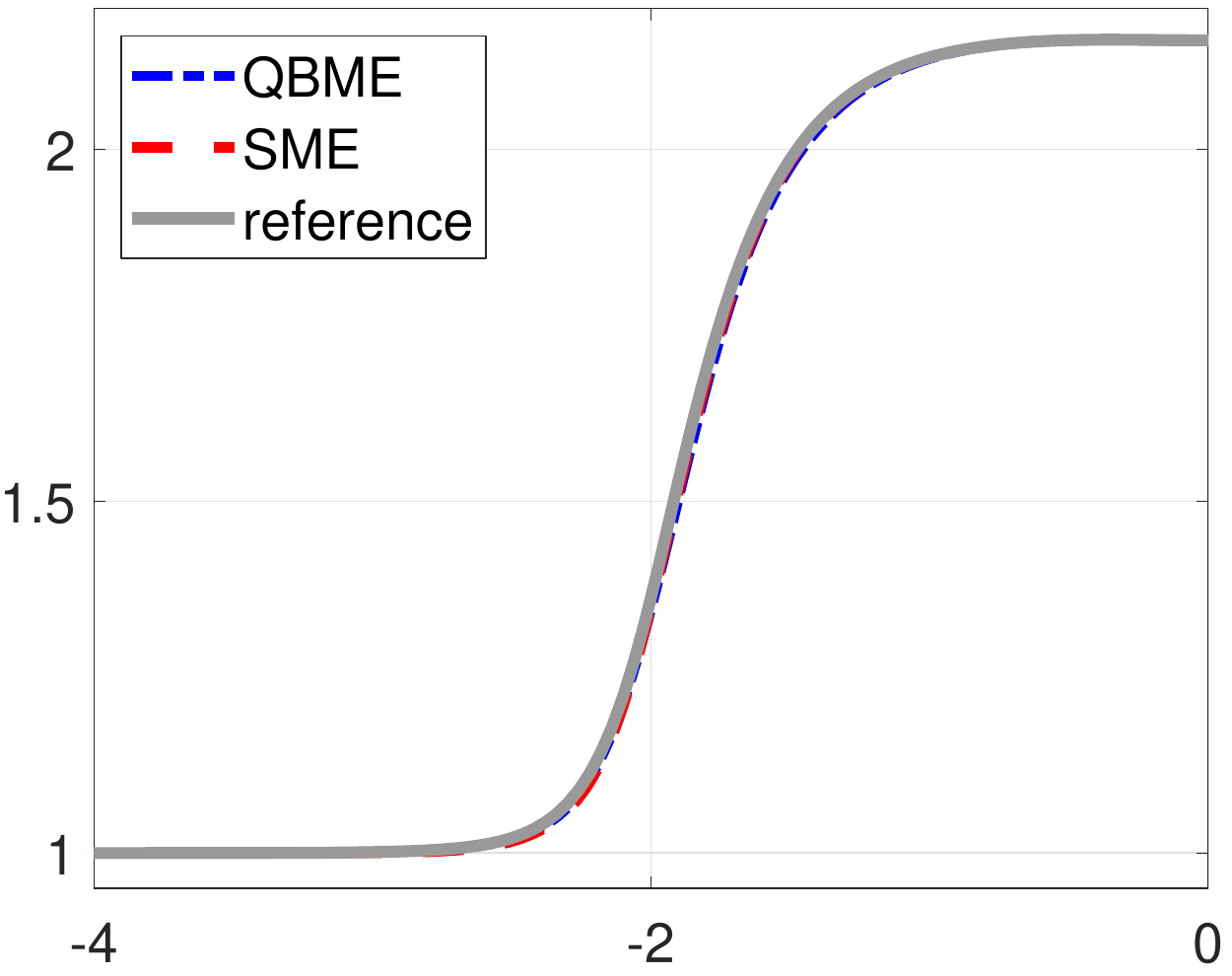}
 \caption{$Kn = 0.1$, $p$}
\end{subfigure}\hfill
 \begin{subfigure}[t]{0.49\textwidth}
 \centering
 \includegraphics[width=\textwidth]{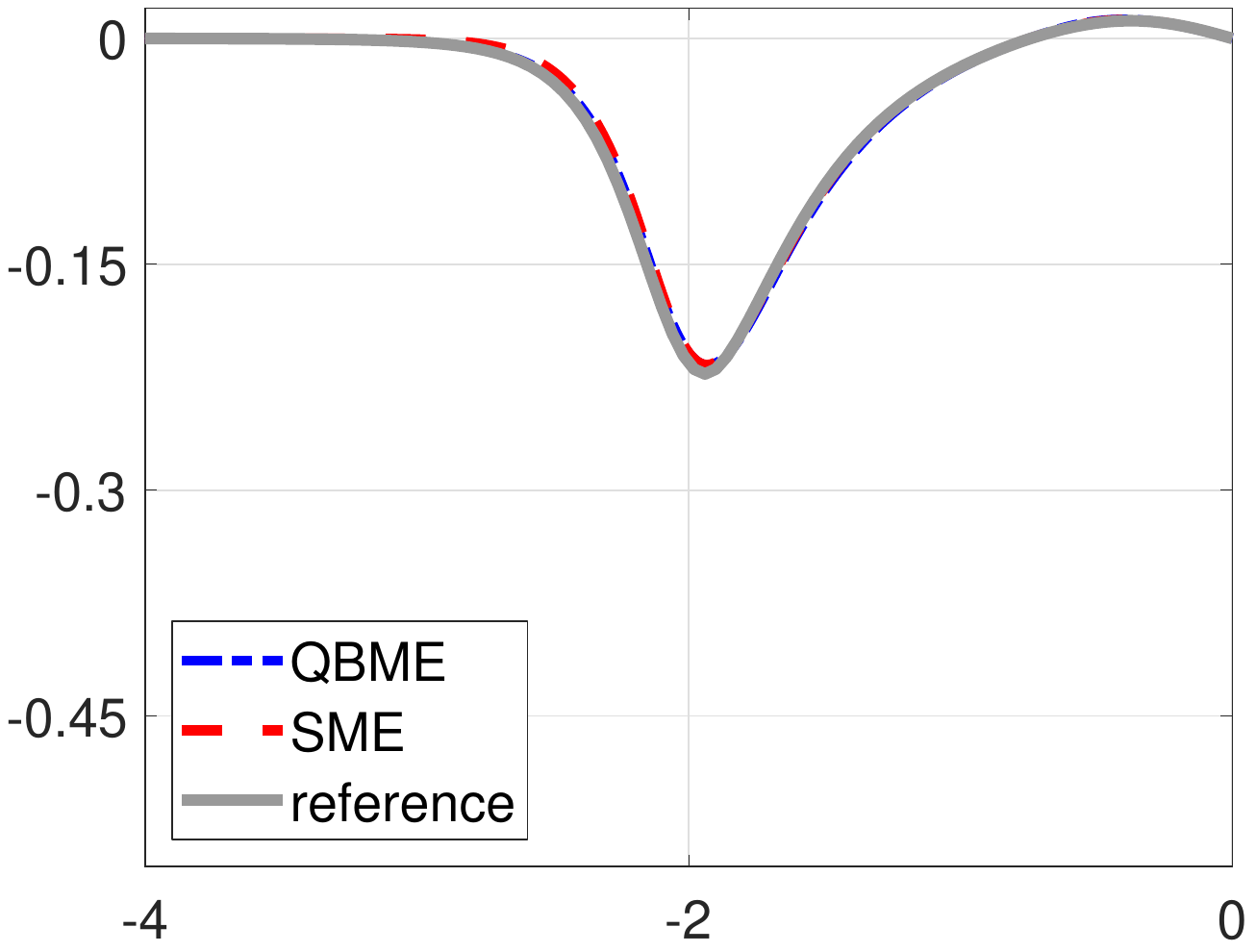}
 \caption{$Kn = 0.1$, $\bar{q}$}
\end{subfigure}
\caption{Symmetric two-beam result for SME and $n=10$. QBME for comparison and DVM for reference. Left column (a),(c), (e) $p$, right column (b), (d), (f) $\bar{q}$. (a),(b) $Kn = \infty$, (c),(d) $Kn = 1$, (e),(f) $Kn = 0.1$.}
\label{2beam}
\end{figure}

To better quantify the convergence of the new model, we plot the errors with decreasing Knudsen number in Figure \ref{2beamconv}. We clearly see that the error of the new SME model is smaller than the known QBME model for both the velocity $u$ and pressure $p$. For the heat flux, the error is of the same order, whereas the QBME model is slightly more accurate for smaller Knudsen number.

We summarize that the SME model successfully approximates the solution of the symmetric two-beam problem for different Knudsen numbers. For most cases, the model outperforms the, already accurate, existing QBME model.

\begin{figure}[H]
 \centering
 \begin{subfigure}[t]{0.3\textwidth}
 \centering
 \includegraphics[width=\textwidth]{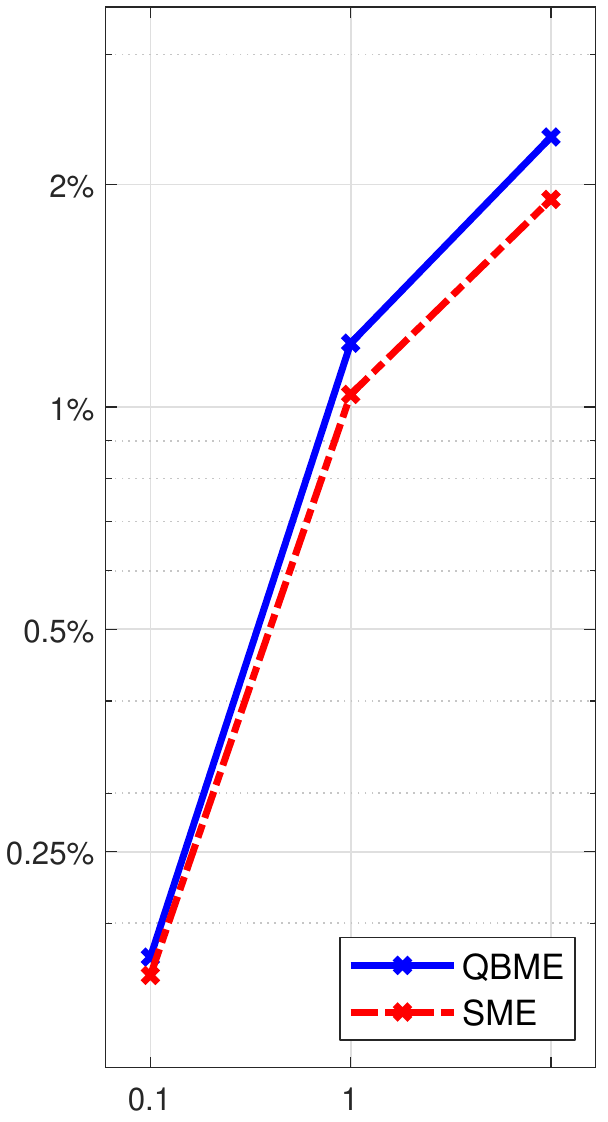}
 \caption{$u$}
\end{subfigure}
\begin{subfigure}[t]{0.3\textwidth}
 \centering
 \includegraphics[width=\textwidth]{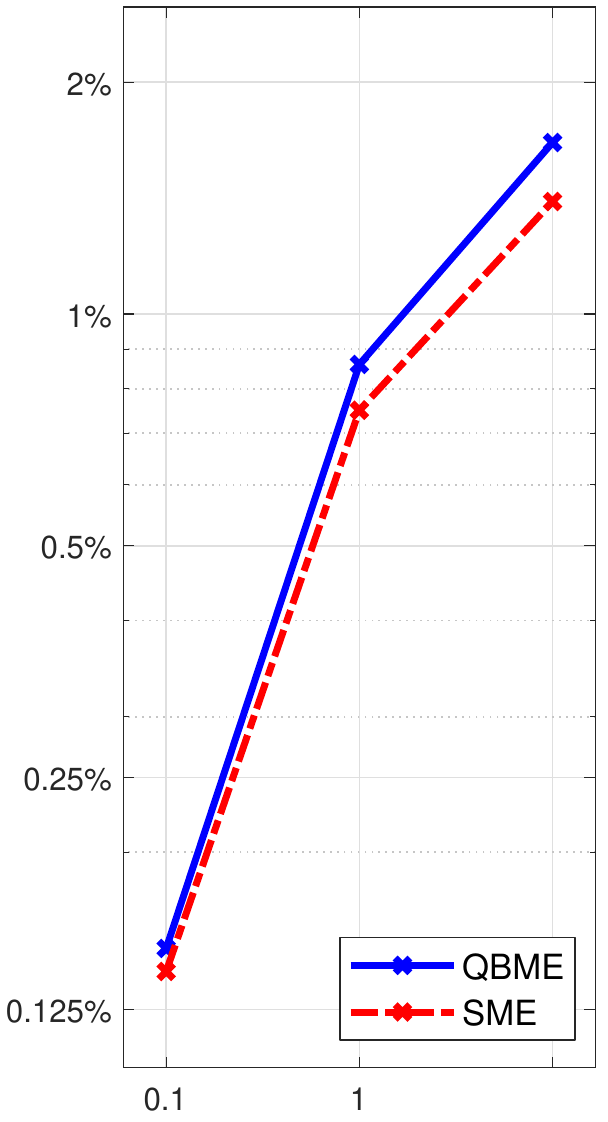}
 \caption{$p$}
\end{subfigure}
\begin{subfigure}[t]{0.3\textwidth}
 \centering
 \includegraphics[width=\textwidth]{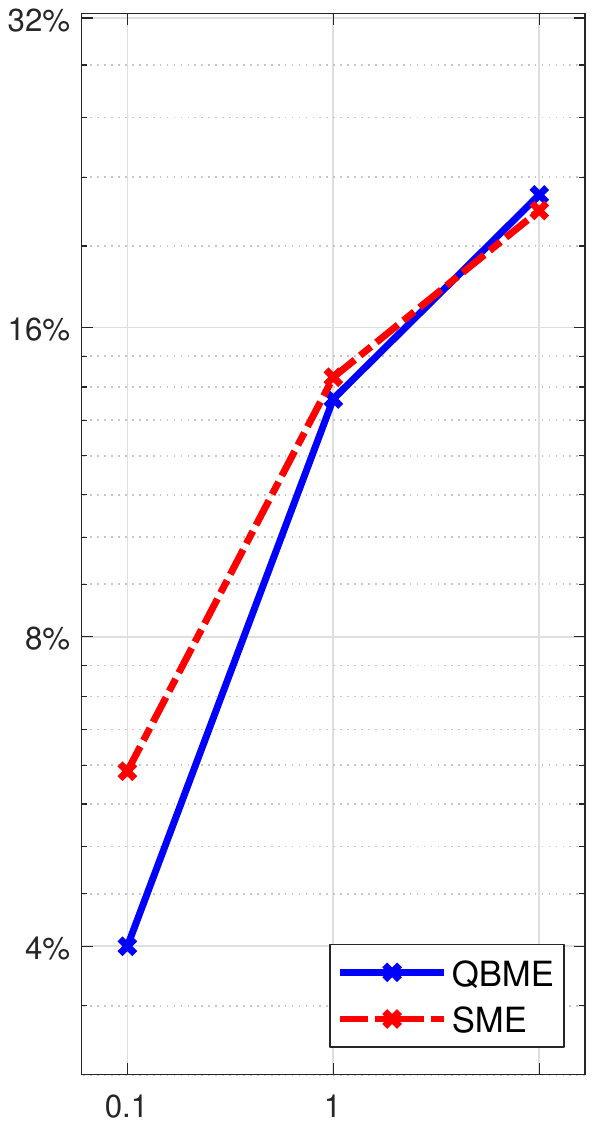}
 \caption{$\bar{q}$}
\end{subfigure}
\caption{Error convergence of symmetric two-beam simulations using $n=10$. (a) velocity $u$, (b) pressure $p$, (c) heat flux $\bar{q}$.}
\label{2beamconv}
\end{figure}

\subsection{Shock structure problem} \label{sec:Shock_structure}
For the last test case, the stationary shock structure problem, we again closely follow the descriptions in \cite{Koellermeier2017b,schaerer2015}. We choose given upstream and downstream boundary conditions from the Rankine-Hugoniot conditions \cite{Ruggeri1993}. Setting the Mach number $\mathrm{Ma} = 1.8$ leads to
\begin{equation}
    \begin{aligned}
        \rho_L &= 1, &
        \rho_R &= \rho_L \cdot \frac{2 \mathrm{Ma}^2}{\mathrm{Ma}^2+1} \approx 1.528,\\
        u_L &= \sqrt{3} \, \mathrm{Ma} \approx 3.118, \qquad&
        \vel_R &= \vel_L \cdot \frac{\mathrm{Ma}^2+1}{2 \mathrm{Ma}^2} \approx 2.040,\\
        \theta_L &= 1,&
        \theta_R &= \theta_L \cdot \frac{\left(1+\mathrm{Ma}^2\right)\left(3\mathrm{Ma}^2-1\right)}{4\mathrm{Ma}^2} \approx 2.853.
    \end{aligned}
\end{equation}

For the relaxation time on the right-hand side collision operator, a constant $\tau=0.01$ is used.
The steady-state solution is computed by time marching from the discontinuous initial data until
convergence. Afterwards the density is scaled using
\begin{equation}\label{e:BC}
     \widetilde{\rho} = \frac{\rho-\rho_L}{\rho_R-\rho_L} \Rightarrow \widetilde{\rho} \in [0,1]
\end{equation}
and the shock positions are aligned at $x=0$ by matching the values of $\widetilde{\rho} = 1/2$.
The computational grid is $[x_L,x_R]=[-78,78]$ using $N_x~=~7500$ cells with a CFL number of
approximately $0.5$. We will use the SME model with $n=7$ equations to match the number of equations from the test case done in \cite{Koellermeier2017b} and use the order $k=1$ and $[\xi_{min},\xi_{max}]=[-4,4]$, similar as before.
The reference DVM solution from \cite{Mieussens2000} includes $N_v=400$ velocity
points and $N_x = 4000$ spatial cells. The reference method is thus far more expensive to compute in comparison to the models used here. We furthermore compare with the QBME model mentioned before and the HME model \cite{Cai2013}, which is widely used in the literature. To allow for a fair comparison, we use $M+1=7$ variables for the other moment models QBME and HME as well. Note that the highly non-linear QBME and HME are based on global basis functions in comparison to the bounded support of the splines used for SME.

We consider density $\rho$, velocity $u$, pressure $p=\rho\theta$, and the normalized heat flux $\bar{q}$ from \ref{e6:normalized_heat_flux} and for our SME model here computed by
\begin{equation} \label{e6:normalized_heat_flux7}
    \bar{q} \approx 0.00699159 \kappa_1 + 0.179474 \kappa_2 + 0.179474 \kappa_3 + 0.00699159 \kappa_4.
\end{equation}

The results in Figure \ref{fig:shockStructure} show that there are only very small differences between the standard hyperbolic moment models QBME, HME and the new SME model. The density, velocity and pressure plots show a good approximation property of the SME model, despite the fact that only $n=7$ variables are used. For the normalized heat flux $\bar{q}$, slight differences with respect to QBME and HME can be seen. However, the differences are of the order of the differences between the hyperbolic QBME and HME themselves. There is no significant additional deviation from the reference solution for the SME model.

The relative error comparison in Table \ref{tab:shockStructureErrors} shows that SME model yields a good approximation quality in comparison with the other models. The error of the new SME model is indeed significantly smaller than the error of the HME model on the one hand. The SME results in approximately the same error as the QBME model on the other hand. The surprising result is that the use of only very few spline functions can compete and even outperform some of the hyperbolic moment models.

\begin{figure}[htb]
 \centering
 \begin{subfigure}[t]{0.47\textwidth}
 \centering
 \includegraphics[width=\textwidth]{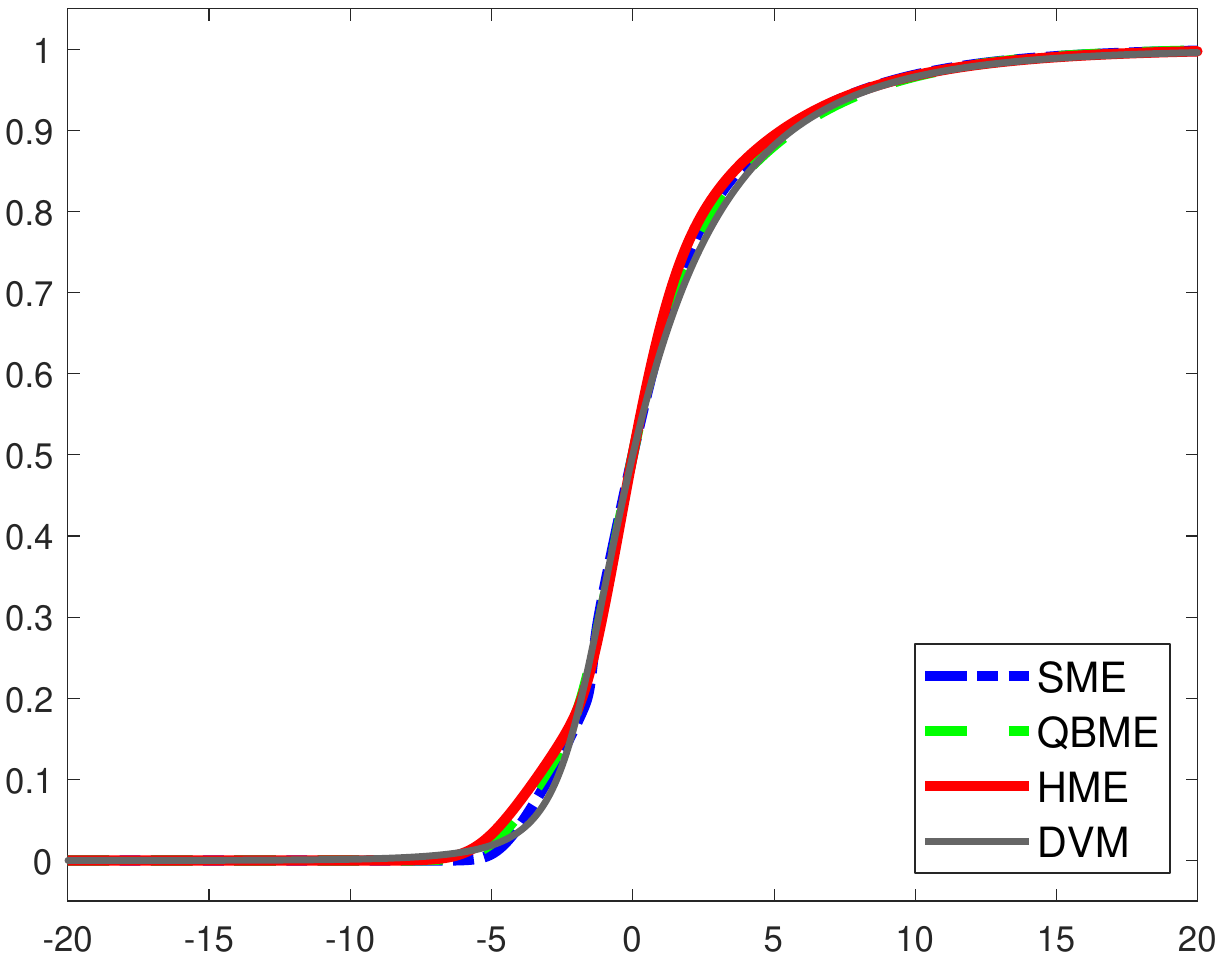}
 \caption{$\rho$}
 \label{fig:shockStructure_rho}
\end{subfigure}
\hfill
\begin{subfigure}[t]{0.47\textwidth}
 \centering
 \includegraphics[width=\textwidth]{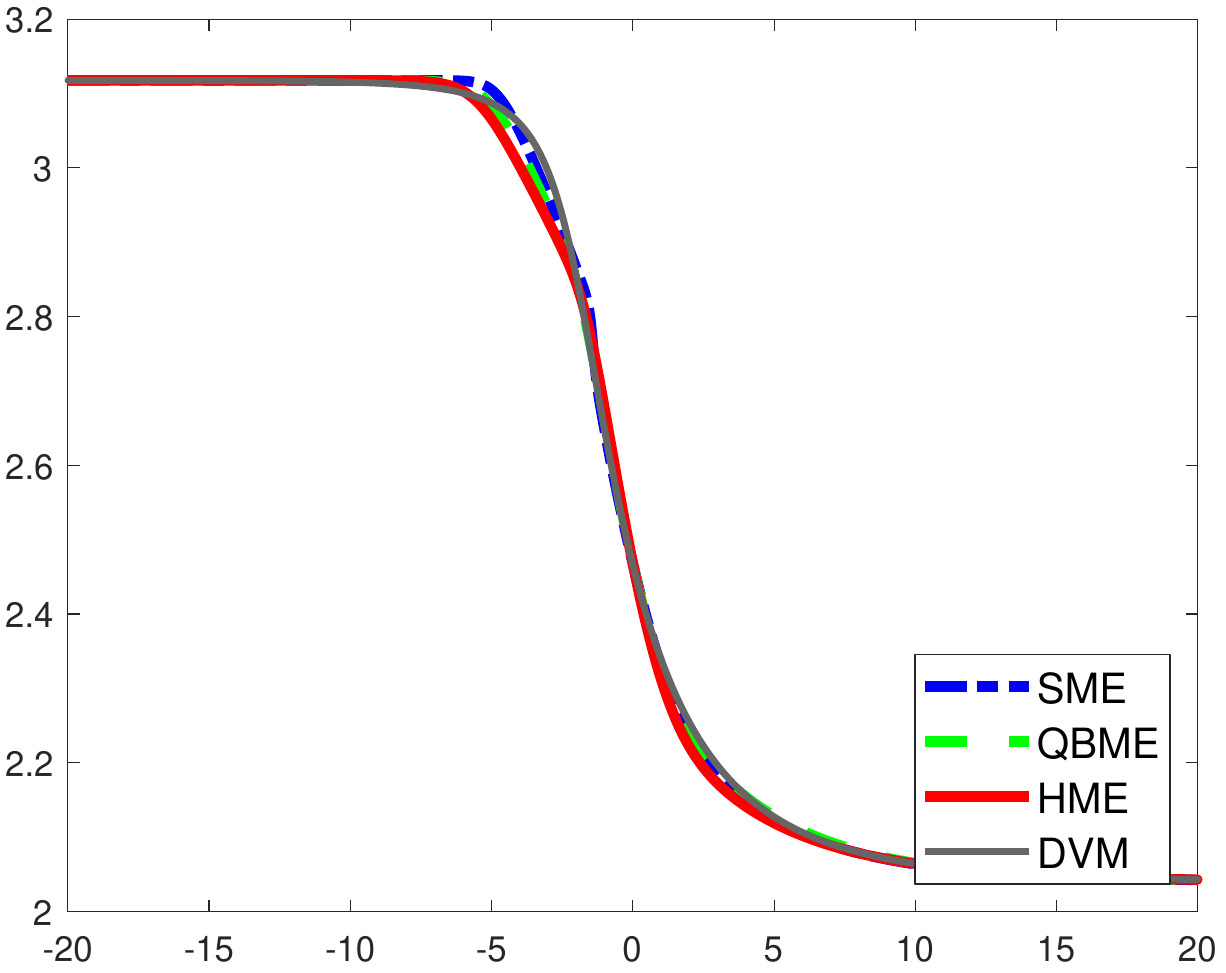}
 \caption{$u$}
 \label{fig:shockStructure_u}
\end{subfigure}
\hfill
 \begin{subfigure}[t]{0.47\textwidth}
 \centering
 \includegraphics[width=\textwidth]{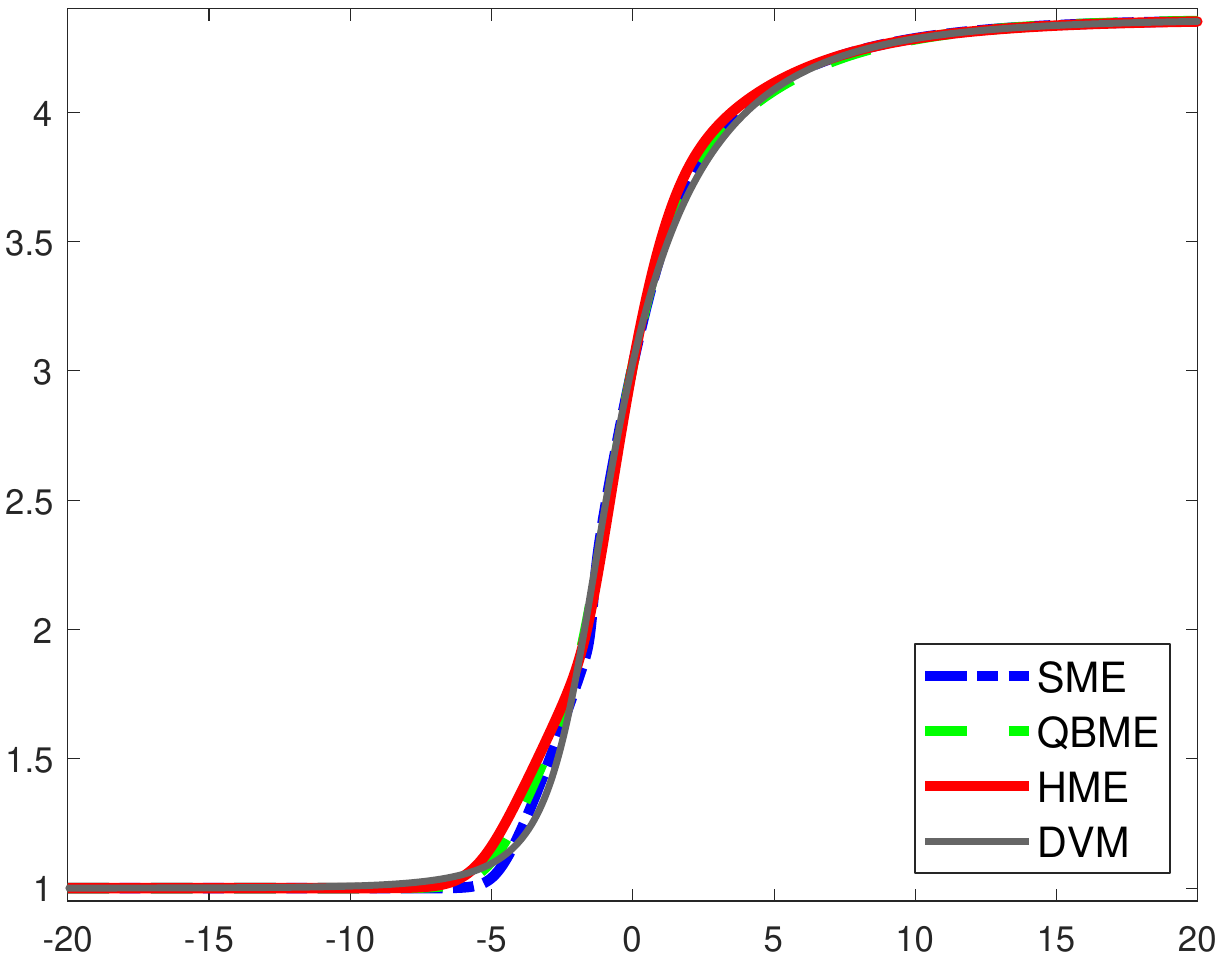}
 \caption{$p$}
 \label{fig:shockStructure_p}
\end{subfigure}
\hfill
\begin{subfigure}[t]{0.47\textwidth}
 \centering
 \includegraphics[width=\textwidth]{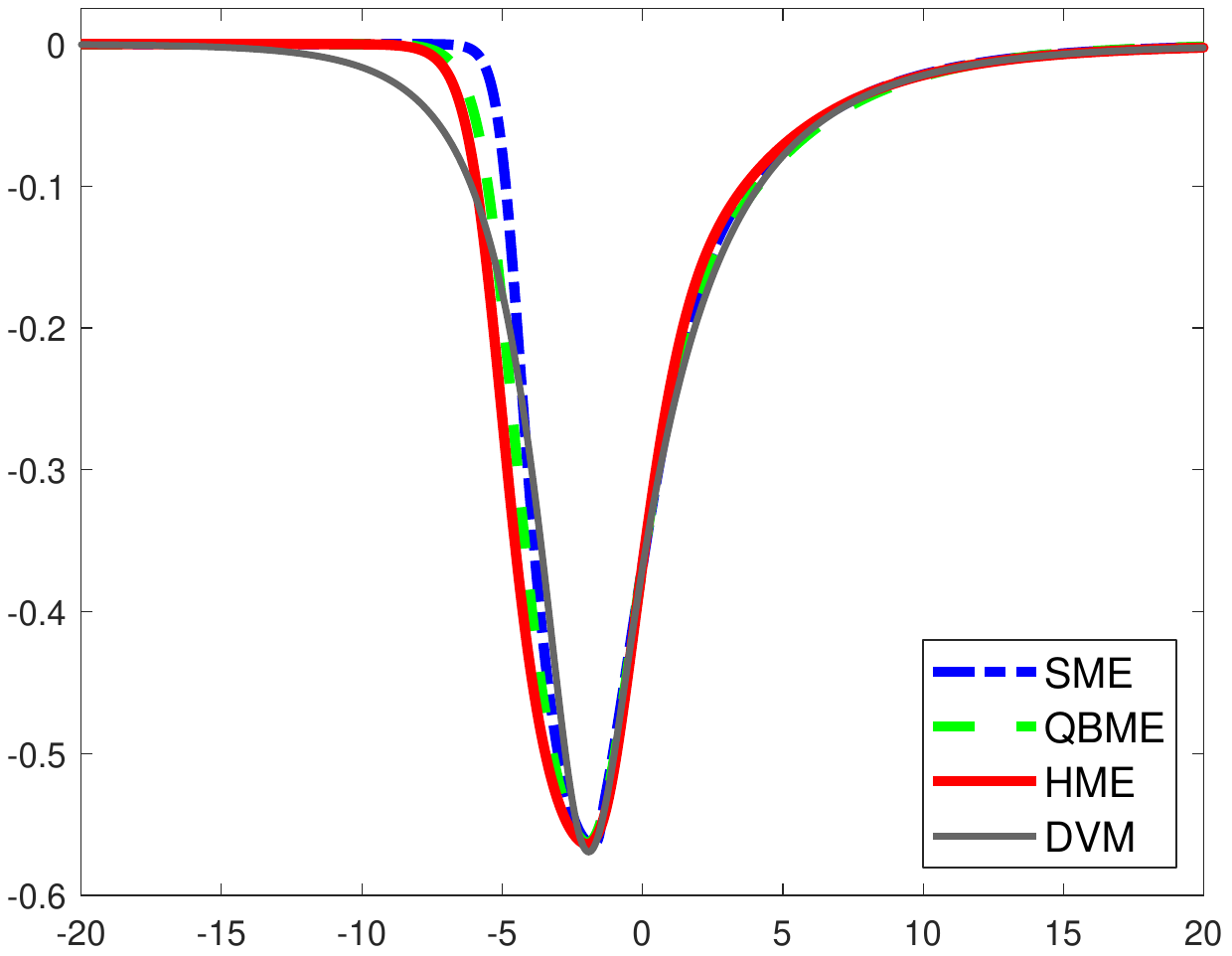}
 \caption{$\bar{q}$}
 \label{fig:shockStructure_Q}
\end{subfigure}
\caption{Shock structure result for SME and $n=7$. QBME, HME for comparison and DVM for reference. (a) density $\rho$, (b) velocity $u$, (c) pressure $p$, (d) heat flux $\bar{q}$.}
\label{fig:shockStructure}
\end{figure}

\begin{table}
\centering
\caption{Shock structure error comparison for different models.}
\label{tab:shockStructureErrors}
\tabcolsep7pt\begin{tabular}{c|c|c|c|c}
    model & $\rho$ & $\vel$ & $p$ & $Q$ \\
\hline
    SME & 1.09\% & 0.20\% & 0.59\% & 14.03\% \\
    QBME & 0.91\% & 0.19\% & 0.58\% & 12.26\% \\
    HME & 1.65\% & 0.32\% & 0.97\% & 16.22\%
\end{tabular}
\end{table}

\section{Conclusion and further work}
\label{sec:7}

In this paper we introduced the first moment model for kinetic theory based on spline basis functions, called Spline Moment Equations (SME).

After a concise definition of the spline ansatz space, we investigated the approximation properties of three classes of splines, namely unweighted splines, weighted splines and weighted fundamental constrained splines (FCS). We saw that the FCS resulted in good accuracy while conserving mass, momentum and energy of the approximated distribution function.

The FCS were used to derive moment equations from a shifted and scaled version of the one-dimensional Boltzmann-BGK equation to allow for an accurate and efficient discretization. The equations could be given in compact analytical form and a subsequent investigation of the hyperbolicity yielded a similar hyperbolicity domain as the more complex Grad model, despite the simplicity of the spline ansatz.

The resulting new SME model was systematically tested using a shock tube test case and the SME model yielded accurate solutions with decreasing error when using more and more equations. In additional two-beam and stationary shock structure test cases, the new SME outperformed several moment models and resulted in better accuracy.

Future work on Spline models for kinetic equations should consider the multi-dimensional extension and a more detailed investigation of the hyperbolicity loss as well as a consistent hyperbolic regularization. An application of the spline ansatz to the depth-averaged shallow water models similar as in \cite{Koellermeier2020c,Kowalski2019} can be promising.

\section*{Acknowledgments}
The research of J. Koellermeier was funded by a joint postdoctoral scholarship from Freie Universit\"at Berlin and Peking University under grand reference number 0503241821.

\appendix
\section{Appendix: Hyperbolic linearized SME }
\label{app}
In the shock tube test, the SME model also yields accurate solutions for the test case using $\textrm{Kn}=0.5$. This is already an improvement with respect to Grad's method, which fails to be stable in this test case \cite{Cai2013}. However, stability problems due to the loss of hyperbolicity may occur in more extreme test cases. We therefore test a simple linearized SME model LSME, in which the system matrix is linearized around equilibrium. A similar strategy (though in a different set of variables) was used in \cite{Cai2013} to achieve hyperbolicity using a modified version of Grad's equations called Hyperbolic Moment Equations (HME). The same strategy was used for shallow flow models in \cite{Koellermeier2020c}

In Figure \ref{linear} the results of the hyperbolic LSME model are presented. The results do not differ much from the SME. However, it is obvious that there remains a systematic error coming from the linearization, especially for the pressure $p$ and velocity $v$. One possible reason is that the linearization has a larger effect on the respective equations for $v,\theta$, while the first equation for $\rho$ is not changed at all.

When increasing the number of splines in Figure \ref{lincmp}, it becomes clear that the linearized SME model does not converge to the solution of the Boltzmann-BGK equation. The error for the pressure remains at a high plateau despite some reductions for velocity and density.

In summary, it must be said that the linearized system does not converge to the exact solution, the linearization seems too much of a simplification. Future work on the investigation of other hyperbolic spline-based models is necessary, see \cite{Fan2016}.


\begin{figure}[H]
 \centering
 \begin{subfigure}[t]{0.49\textwidth}
 \centering
 \includegraphics[width=\textwidth]{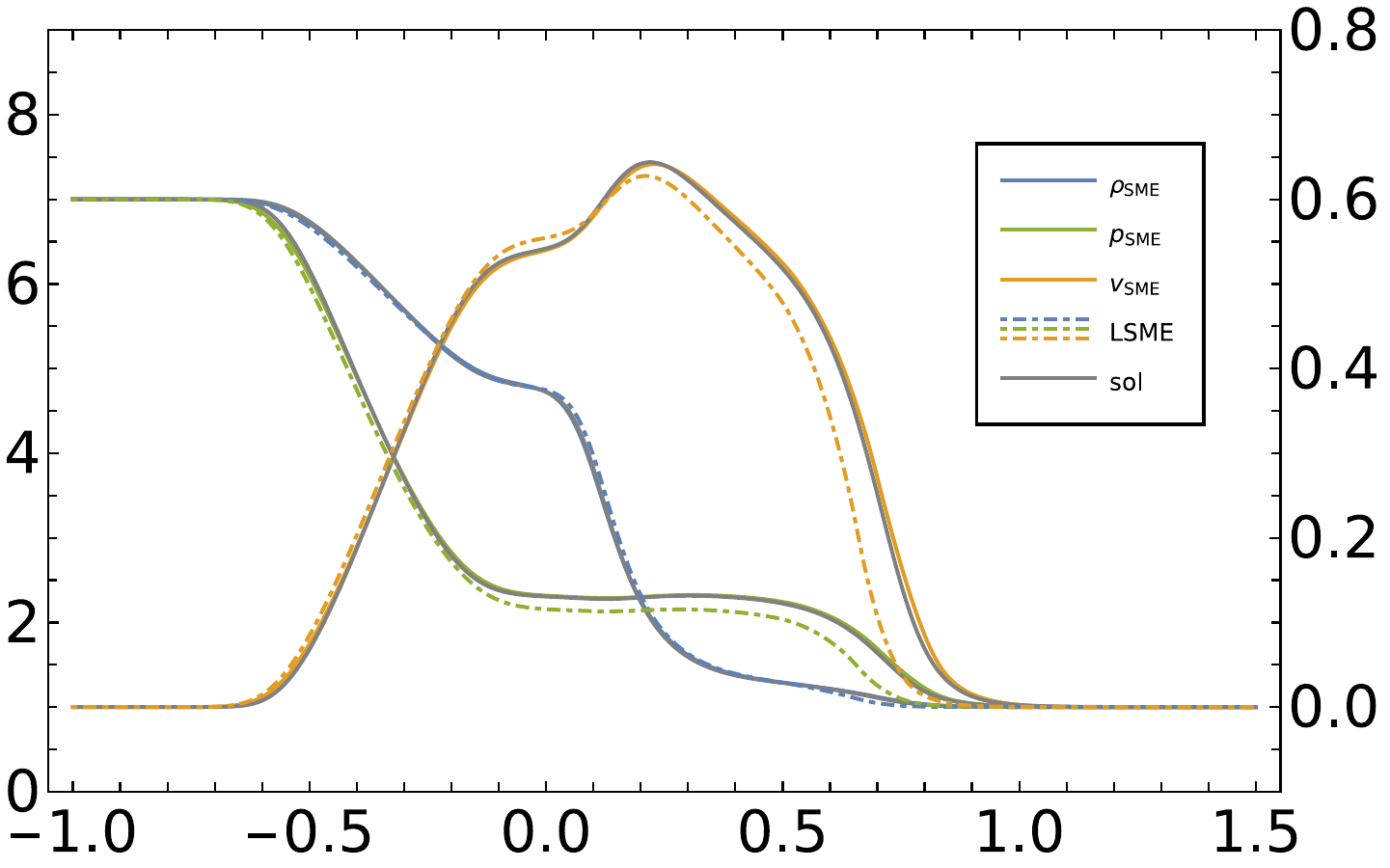}
 \caption{$\textrm{Kn}=0.05$}
 \label{linear5}
\end{subfigure}
\hfill
\begin{subfigure}[t]{0.49\textwidth}
 \centering
 \includegraphics[width=\textwidth]{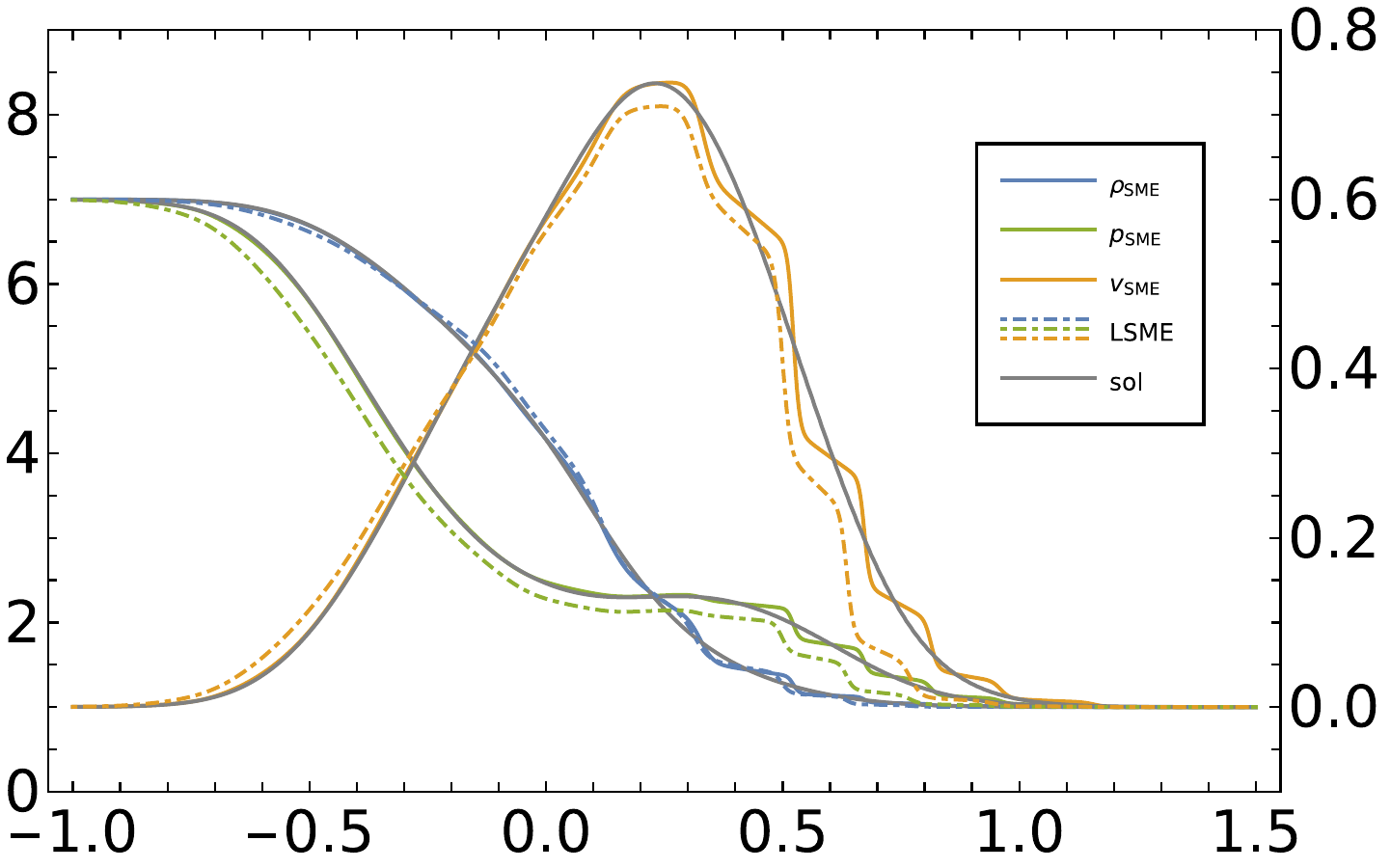}
 \caption{$\textrm{Kn}=0.5$}
 \label{linear6}
\end{subfigure}
\caption{Shock tube linearized SME comparison. $k=1$, $n=12$, $[\xi_{min},\xi_{max}]=[-4,4]$. (a) $\textrm{Kn}=0.05$, (b) $\textrm{Kn}=0.5$.}
\label{linear}
\end{figure}

\begin{figure}[H]
 \centering
 \begin{subfigure}[t]{0.49\textwidth}
 \centering
 \includegraphics[width=\textwidth]{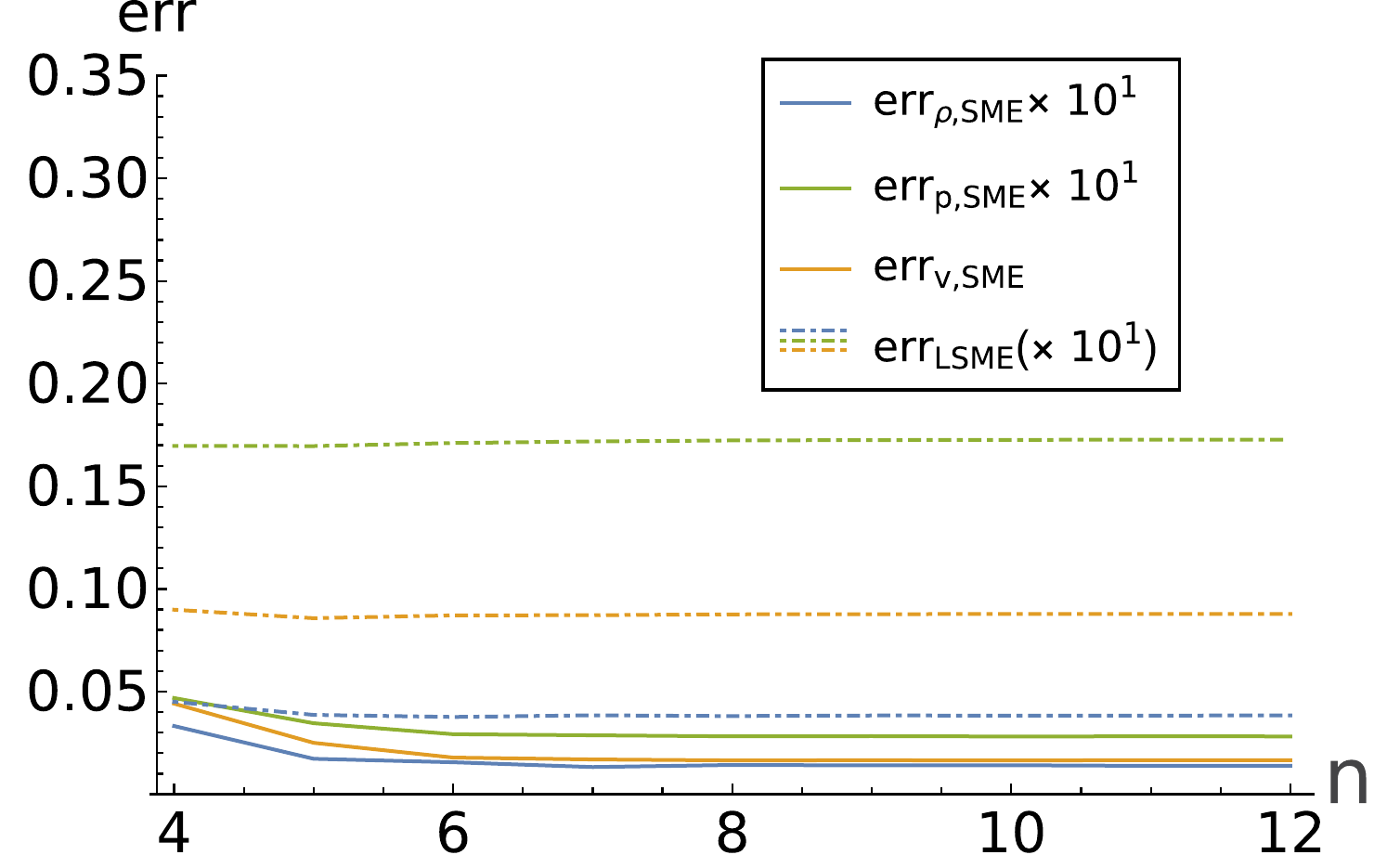}
 \caption{$\textrm{Kn}=0.05$.}
 \label{lincmp5}
\end{subfigure}
\hfill
\begin{subfigure}[t]{0.49\textwidth}
 \centering
 \includegraphics[width=\textwidth]{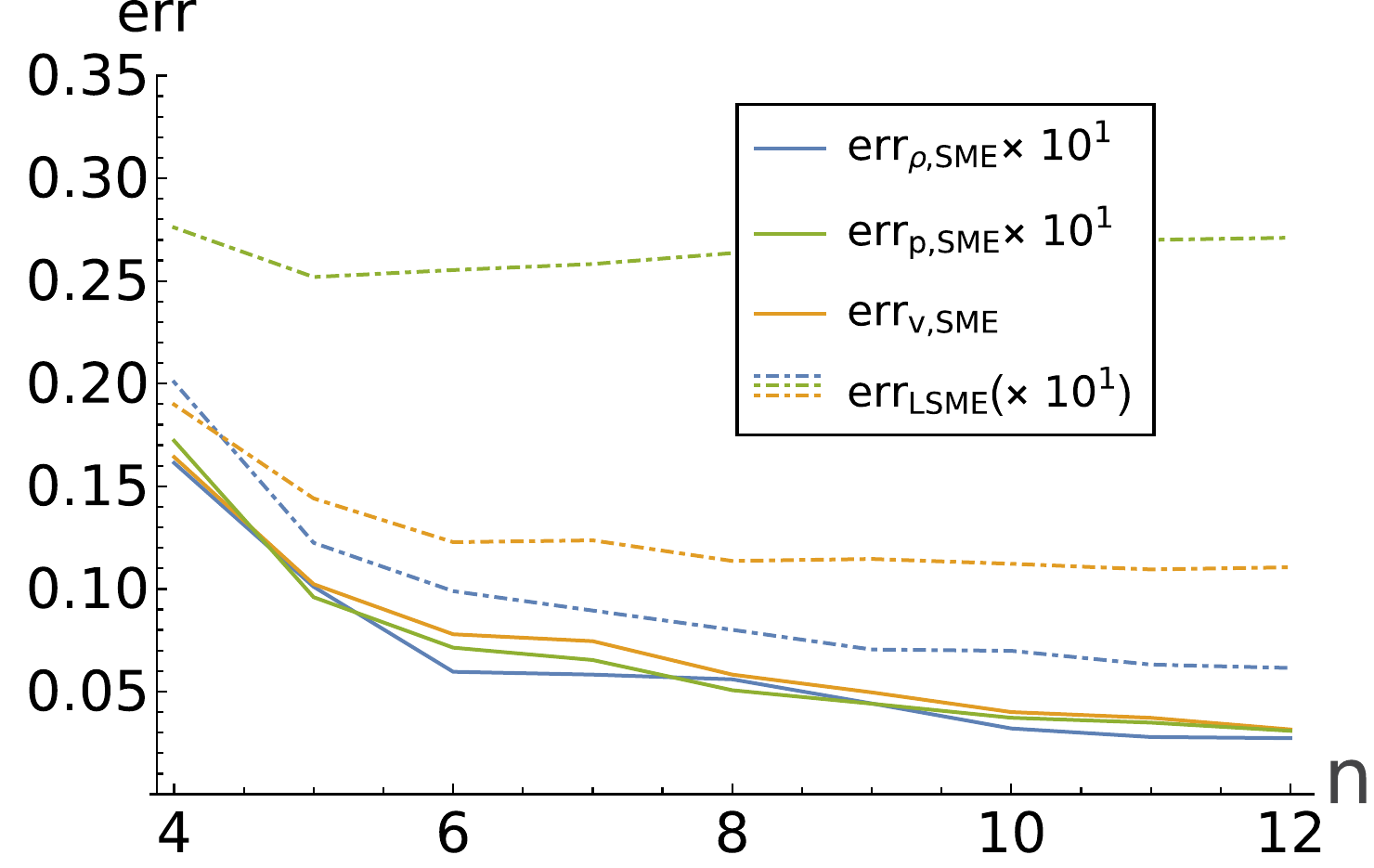}
 \caption{$\textrm{Kn}=0.5$.}
 \label{lincmp6}
\end{subfigure}
\caption{Simulation error depending on number of splines for linearized SME, $k=3$, $[\xi_{min}, \xi_{max}]=[-4,4]$. (a) $\textrm{Kn}=0.05$, (b) $\textrm{Kn}=0.5$.}
\label{lincmp}
\end{figure}

\section*{Data Availability Statement}
The datasets generated during or analysed during the current study are available from the corresponding author on reasonable request.
\bibliographystyle{abbrv}

\end{document}